\documentclass[11pt]{article}

\usepackage[french]{babel}
\usepackage[T1]{fontenc}
\usepackage{comment}
\usepackage{amssymb,epsfig} 
\usepackage{amsmath}
\usepackage[dvips]{color}
\usepackage{mathrsfs}

\newtheorem{defi}{D\'efinition}[section]
\newtheorem{prop}[defi]{Proposition}
\newtheorem{theo}[defi]{Th\'eor\`eme}
\newtheorem{conj}[defi]{Conjecture}
\newtheorem{lemm}[defi]{Lemme}
\newtheorem{coro}[defi]{Corollaire}
\newtheorem{rema}[defi]{Remarque}
\newtheorem{exem}[defi]{Exemple}
\newtheorem{exer}{Exercice}

\newcommand{\bdefi}{\begin{defi}}
\newcommand{\edefi}{\end{defi}}
\newcommand{\bprop}{\begin{prop}}
\newcommand{\eprop}{\end{prop}}
\newcommand{\btheo}{\begin{theo}}
\newcommand{\etheo}{\end{theo}}
\newcommand{\blemm}{\begin{lemm}}
\newcommand{\brema}{\begin{rema}}
\newcommand{\erema}{\end{rema}}
\newcommand{\bexer}{\begin{exer}}
\newcommand{\eexer}{\end{exer}}
\newcommand{\bconj}{\begin{conj}}
\newcommand{\econj}{\end{conj}}
\newcommand{\elemm}{\end{lemm}}
\newcommand{\bcoro}{\begin{coro}}
\newcommand{\ecoro}{\end{coro}}
\newcommand{\bexem}{\begin{exem}}
\newcommand{\eexem}{\end{exem}}
\newcommand{\dem}{\noindent{\bf D\'emonstration. }}
\newcommand{\rem}{\noindent{\bf Remarque. }}

\newcommand{\A}{{\cal A}}
\newcommand{\R}{{\cal R}}

\newcommand{\M}{{\cal M}}
\renewcommand{\P}{{\cal P}}
\newcommand{\G}{{\cal G}}

\newcommand{\D}{{\cal D}}
\newcommand{\E}{{\cal E}}

\renewcommand{\H}{{\cal H}}
\newcommand{\OOO}{{\cal O}}
\newcommand{\C}{{\cal C}}

\newcommand{\F}{{\cal F}}
\newcommand{\V}{{\cal V}}

\newcommand{\N}{{\cal N}}
\renewcommand{\P}{{\cal P}}

\usepackage{amssymb}
\usepackage{amsmath}

\newcommand{\maths}[1]{{\mathbb #1}}  

\newcommand{\RR}{\maths{R}}
\newcommand{\NN}{\maths{N}}
\newcommand{\CC}{\maths{C}}
\newcommand{\QQ}{\maths{Q}}

\newcommand{\SSS}{\maths{S}}

\newcommand{\HH}{\maths{H}}

\newcommand{\ZZ}{\maths{Z}}

\newcommand{\PP}{\maths{P}}

\newcommand{\aaa}{{\mathfrak a}}

\newcommand{\ppp}{{\mathfrak p}}

\newcommand{\ra}{\rightarrow}

\newcommand{\ov}[1]{{\overline{#1}}} 
\newcommand{\wt}[1]{{\widetilde{#1}}}
\newcommand{\wh}[1]{{\widehat{#1}}}

\newcommand{\ga}{\gamma}
\newcommand{\Ga}{\Gamma}

\newcommand{\cqfd}{\hfill$\Box$}

\newcommand{\PSL}[2]{{\operatorname{PSL}_{#1}(#2)}}
\newcommand{\SL}[2]{\operatorname{SL}_{#1}(#2)}

\newcommand{\hdr}{{\HH}^2_\RR}
\newcommand{\htr}{{\HH}^3_\RR}

\renewcommand\mathcal{\mathscr}

\newcommand{\note}[1]{\marginpar{\scriptsize #1}}
\newcommand{\bs}{\backslash}
\newcommand{\id}{\operatorname{id}}

\newcommand{\Vol}{\operatorname{Vol}}
\newcommand{\vol}{\operatorname{vol}}

\newcommand{\arcosh}{{\operatorname{argcosh}}}
\newcommand{\artanh}{{\operatorname{argtanh}}}

\renewcommand{\S}{\textsection}

\newcommand{\card}{\operatorname{Card}}
\newcommand{\autom}{\operatorname{SO}}
\newcommand{\PSO}{\operatorname{PSO}(Q,\ZZ)}
\newcommand{\discr}{\operatorname{D}}
\newcommand{\Stab}{\operatorname{Stab}}
\newcommand{\Tr}{\operatorname{Tr}}

\newcounter{fig}

\usepackage{epsfig}
\usepackage[dvips]{color}

\title{\'Equidistribution, comptage et approximation par 
irrationnels quadratiques}
\author{Jouni Parkkonen \and Fr\'ed\'eric Paulin}
\date{\today} 

\begin{document} 
\bibliographystyle{../alphanum}
\maketitle 

\begin{abstract} 
\noindent   
Soit $M$ une vari\'et\'e hyperbolique de volume fini, nous montrons
que les hypersurfaces \'equidistantes \`a une sous-vari\'et\'e $C$ de
volume fini totalement g\'eod\'esique s'\'equidistribuent dans
$M$. Nous donnons un asymptotique pr\'ecis du nombre de segments
g\'eod\'esiques de longueur au plus $t$, perpendiculaires communs \`a
$C$ et au bord d'un voisinage cuspidal de $M$. Nous en d\'eduisons des
r\'esultats sur le comptage d'irrationnels quadratiques sur $\QQ$ ou
sur une extension quadratique imaginaire de $\QQ$, dans des orbites
donn\'ees des sous-groupes de congruence des groupes modulaires.

\begin{center}
{\bf Abstract}
\end{center}

Let $M$ be a finite volume hyperbolic manifold, we show the
equidistribution in $M$ of the equidistant hypersurfaces to a finite
volume totally geodesic submanifold $C$. We prove a precise
asymptotic on the number of geodesic arcs of lengths at most $t$, that
are perpendicular to $C$ and to the boundary of a cuspidal
neighbourhood of $M$. We deduce from it counting results of quadratic
irrationals over $\QQ$ or over imaginary quadratic extensions of
$\QQ$, in given orbits of congruence subgroups of the modular groups.
\end{abstract}

\section{Introduction}

Soit $M$ une vari\'et\'e hyperbolique (c'est-\`a-dire riemannienne
lisse, compl\`ete et \`a courbure sectionnelle cons\-tan\-te $-1$)
connexe, de volume fini, de dimension $n$ au moins $2$.  Soit $C$ une
sous-vari\'et\'e immerg\'ee de volume fini de $M$, connexe, de
dimension $k<n$, totalement g\'eod\'esique (par exemple une
g\'eod\'esique ferm\'ee). Soit $\Sigma$ une composante connexe du
fibr\'e normal unitaire de $C$. Pour tout $t>0$, notons $\Sigma(t)=
\{\exp tv\;:\; v\in\Sigma\}$ l'hypersurface strictement convexe
immerg\'ee de $M$, pouss\'ee au temps $t$ le long des rayons
g\'eod\'esiques orthogonaux \`a $C$ dans la direction $\Sigma$.

Notre premier r\'esultat dit que la moyenne riemannienne de $\Sigma(t)$
s'\'equidistribue quand $t$ tend vers $+\infty$ vers la moyenne
riemannienne de $M$. 

\btheo\label{theo:introequidistrib} 
Pour toute fonction $\varphi$ continue \`a support compact dans $M$, 
nous avons
$$
\lim_{t\ra+\infty}\oint_M \varphi(x)\; d\vol_{\Sigma(t)}(x) = 
\oint_M \;\varphi(x)\; d\operatorname{vol}_{M}(x)\;.
$$
\etheo

Lorsque $n=2$ et $\Sigma$ est le fibr\'e normal unitaire tout entier
de $C$, ce r\'esultat est d\^u \`a \cite{EskMcMul93}, sur des id\'ees
remontant \`a la th\`ese de Margulis (voir par exemple
\cite{Margulis04}). Dans celle-ci, l'\'equidistribution de grandes
horosph\`eres, que l'on peut voir comme les $\Sigma(t)$ o\`u $\Sigma$
est le fibr\'e normal unitaire sortant du bord (non totalement
g\'eod\'esique) d'un voisinage cuspidal (ou de Margulis, voir par
exemple \cite{BusKar91}) d'un bout de $M$, est d\'emontr\'ee. Nous
renvoyons \`a \cite{Babillot02b,Roblin03} pour des r\'esultats
g\'en\'eraux d'\'equidistribution d'horosph\`eres, et par exemple \`a
\cite{EskMcMul93} (dont les m\'ethodes ne peuvent s'appliquer dans
notre situation si $n\geq 3$ et $k>0$, l'espace homog\`ene $\PSL
2\CC/\Big\{\Big(\begin{array}{cc}a &0\\0&
a^{-1}
\end{array}\Big)\;:\;a>0 \Big\}$ n'\'etant par exemple
pas sym\'etrique affine) et aux travaux de Benoist, Duke, Einsiedler,
Eskin, Gorodnick, Lindenstrauss, Margulis, Mozes, Oh, Ratner, Quint,
Rudnick, Sarnak, Shah, Tomanov, Ullmo, Venkatesh, Weiss pour d'autres
r\'esultats d'\'equidistribution dans des espaces homog\`enes.

\medskip Nous utilisons ensuite ce th\'eor\`eme pour montrer un
r\'esultat de comptage de rayons g\'eod\'esiques partant
perpendiculairement de $C$ et convergeant dans un bout de $M$. Pour
normaliser les longueurs, notons $\H$ un voisinage cuspidal d'un bout
de $M$.  Pour tout $t\geq 0$, notons $\N(t)$ le
nombre (compt\'e avec multiplicit\'e) de segments (localement)
g\'eod\'esiques perpendiculaires commun \`a $C$ et \`a $\partial\H$,
de longueur au plus $t$.

\btheo\label{theo:introcomptage} Quand $t$ tend vers $+\infty$, nous
avons l'\'equivalent
$$
\N(t)\sim 
\frac{\Vol(\SSS_{n-k-1})\Vol(\H)\Vol(C)}{\Vol(\SSS_{n-1})\Vol(M)}
\;e^{(n-1)t}\;.
$$
\etheo

L'origine de ces r\'esultats asymptotiques pr\'ecis remonte \`a Huber
\cite{Huber59} en dimension $2$, et nous renvoyons particuli\`erement
\`a \cite{Hermann62} pour d'autres r\'esultats de comptages d'arcs
perpendiculaires \`a des sous-vari\'et\'es totalement
g\'eod\'esiques. G\'en\'eralisant les r\'esultats de Huber en courbure
variable et en dimension quelconque, Margulis (voir par exemple
\cite{Margulis04}) donne un \'equivalent du nombre d'arcs (localement)
g\'eod\'esiques entre deux points de $M$ de longueur au plus $t$. Voir
par exemple \cite{BelHerPau01,Roblin03,HerPau04} pour les comptages
d'arcs perpendiculaires \`a des horosph\`eres.

\bigskip Nous donnons dans cette introduction trois exemples
d'applications arithm\'etiques du th\'eor\`eme
\ref{theo:introcomptage} (et de ses variantes, expos\'ees \`a la fin
de la partie \ref{sect:comptage}).

Soit $Q(X,Y)=AX^2+BXY+CY^2$ une forme quadratique binaire enti\`ere,
et $\discr=\discr_Q=B^2-4AC$ son discriminant. Notons $\P$ l'ensemble
des \'el\'ements de $\ZZ^2$ \`a coordonn\'ees premi\`eres entre
elles. Une {\it repr\'esentation} (resp.~{\it repr\'esentation propre})
d'un entier $n\in\ZZ$ par $Q$ est un \'el\'ement $x\in\ZZ^2$
(resp.~$x\in\P$) tel que $Q(x)=n$. Lorsque $Q$ est d\'efinie positive
(en particulier, $\discr<0$), Gauss a donn\'e un \'equivalent quand
$t$ tend vers $+\infty$ du nombre de repr\'esentations propres par $Q$
d'entiers inf\'erieurs \`a $t$ (voir par exemple
\cite{Landau58}, \cite[Chap.~10]{Cohn62})~:
$$
\card\big(\{x\in\P\;:\; Q(x)\leq t\}\big)\sim
\frac{12}{\pi\sqrt{-\discr_Q}}\;t\;.
$$
Supposons maintenant que $Q$ soit primitive ($A,B,C$ sont premiers
entre eux), ind\'efinie ($\discr>0$), et ne soit pas, à multiple
rationnel près, un produit de deux formes lin\'eaires enti\`eres
($\discr$ n'est pas un carr\'e). Au lieu de compter des points entiers
dans des ellipses, il s'agit maintenant de compter des points entiers
entre des hyperboles (voir par exemple \cite{Cohn62} et la partie
\ref{subsec:comptrepquadbin} pour des dessins). Notons $\autom(Q,\ZZ)$
le groupe sp\'ecial orthogonal entier de $Q$ (appel\'e le {\it groupe
  des automorphes} de $Q$), qui est infini. Notons
$$
\Psi_Q(t)=\card\big(\autom(Q,\ZZ)\bs\{x\in\P\;:\; |Q(x)|\le t\}\big)
$$
le nombre de classes modulo $\autom(Q,\ZZ)$ de repr\'esentations
propres par $Q$ d'entiers inf\'erieurs \`a $t$ en valeur absolue. Soit
$(t_Q,u_Q)$ la solution fondamentale (c'est-\`a-dire $t_Q,u_Q>0$ et
$u_Q$ est minimal) de l'\'equation de Pell-Fermat $t^2-\discr\,
u^2=4$, notons $ \epsilon_{Q}=\frac{t_Q+u_Q\sqrt{\discr}}{2}$ et
$R_Q=\log \epsilon_Q$.

\bcoro\label{coro:introformes} Quand $s$ tend vers $+\infty$, nous
avons $\Psi_Q(t)\sim \frac{12\, R_Q}{\pi^2\sqrt{\discr_Q}}\; t$.
\ecoro

Ce r\'esultat est connu (voir par exemple \cite[page 164]{Cohn62} pour
le comptage des repr\'esentations impropres). Mais notre preuve permet
de nombreuses variations, voir en particulier la partie
\ref{subsec:comptrepquadbin} lorsque l'on impose de plus des
conditions de congruences sur les coordonn\'ees des points entiers.

\medskip Passons au second exemple d'application arithm\'etique. Les
r\'esultats de comptage, et de comptage asymptotique, de nombres
alg\'ebriques, de complexit\'e born\'ee (degr\'e, hauteur, etc),
satisfaisant des contraintes alg\'ebriques donn\'ees, sont au coeur de
la th\'eorie des nombres. Si $K$ est un corps de nombres, d'anneau des
entiers $\OOO_K$, le groupe $\operatorname{PSL}_2(\OOO_K)$ agit sur les
irrationnels quadratiques sur $K$, mais non transitivement, et il est
int\'eressant de donner une asymptotique pr\'ecise du nombre
d'\'el\'ements dans une orbite donn\'ee ayant une complexit\'e
born\'ee. Pour tout irrationnel quadratique $\alpha$ sur $K$, nous
noterons $\alpha^\sigma$ son conjugu\'e de Galois sur $K$. Les
r\'esultats de ce papier sont motiv\'es par des probl\`emes,
introduits et d\'evelopp\'es dans \cite{ParPau10MZ}, d'approximation
diophantienne d'\'el\'ements de $\RR$ ou $\CC$ par des irrationnels
quadratiques sur $\QQ$ (ou sur une extension imaginaire quadratique de
$\QQ$) qui appartiennent \`a une orbite donn\'ee du groupe modulaire (ou
du groupe de Bianchi). Dans ce cadre, l'inverse de la distance entre
$\alpha$ et $\alpha^\sigma$ est une complexit\'e pertinente, lorsque
$\alpha$ est consid\'er\'e modulo translations enti\`eres (voir la
partie \ref{subsec:justifcomplexiteh} pour des justifications). Nous
renvoyons par exemple \`a \cite{Bugeaud04,Waldschmidt00} pour des
monographies sur l'approximation diophantienne par des nombres
alg\'ebriques.

Notre second exemple d'application arithm\'etique est un r\'esultat de
comptage (pour la complexit\'e ci-dessus) d'irrationnels quadratiques
sur $\QQ$ dans une orbite du groupe modulaire
$\operatorname{PSL}_2(\ZZ)$ (\'etendu par l'automorphisme de
Galois). Fixons $\alpha_0$, un entier quadratique sur $\QQ$
irrationnel r\'eel. Notons $\OOO=\ZZ+\alpha_0\ZZ$ l'ordre qu'il
d\'efinit dans $\QQ(\alpha_0)$.  Soient $D=\Tr(\alpha_0)^2
-4N(\alpha_0)$ le discriminant de $\OOO$, $(t_0,u_0)$ la solution
fondamentale de l'\'equation de Pell-Fermat $t^2-Du^2=4$, $\epsilon_0=
\frac{t_0+u_0\sqrt D}{2}$ l'unit\'e fondamentale de $\OOO$, et
$R_0=\log\epsilon_0$ le r\'egulateur de $\OOO$. Notons $n_0=2$ si $\OOO$
contient une unit\'e de norme $-1$, et $1$ sinon.

\bcoro\label{coro:introcompirrquad}
Quand $t$ tend vers $+\infty$, nous avons 
$$
\operatorname{Card} \big(\{\alpha\in
\operatorname{PSL}_2(\ZZ)\cdot\{\alpha_0,\alpha_0^\sigma\}\!\!\!
\mod \ZZ\;:\; |\alpha-\alpha^\sigma|^{-1}\leq t\}\big) \sim
\frac{24\;R_0}{\pi^2n_0}\; \;t\;.
$$ 
\ecoro

Nous renvoyons aux parties \ref{subsec:comptirraquad} et
\ref{subsec:comptirraquadimag} pour des r\'esultats plus g\'en\'eraux
(en particulier concernant le comptage d'irrationnels quadratiques
dans les orbites de sous-groupes de congruence de
$\operatorname{PSL}_2(\ZZ)$), et des r\'esultats de comptage
d'irrationnels quadratiques sur une extension quadratique imaginaire
de $\QQ$ dans des orbites par les groupes de Bianchi
$\operatorname{PSL}_2(\OOO_{\QQ(i\sqrt{d})})$, dont un cas
tr\`es particulier est le suivant.

\bcoro\label{coro:introorbianchi} Soient $\phi=\frac{1+\sqrt{5}}{2}$
le nombre d'or et $K$ un corps de nombres quadratique imaginaire,
d'anneau des entiers $\OOO_K$, de discriminant $D_K$ et de fonction
z\'eta $\zeta_K$. Quand $t$ tend vers $+\infty$, nous avons
$$
\operatorname{Card} \big(\{\alpha\in
\operatorname{PSL}_2(\OOO_K)\cdot\phi \!\!\!\mod \OOO_K\;:\; 
|\alpha-\alpha^\sigma|^{-1}\leq t\}\big) \sim
\frac{8\,\pi^2\;\log\phi}{|D_K|\;\zeta_K(2)}\; \;t^2
$$ 
si $D_K\neq -4$, la moiti\'e sinon. 
\ecoro

Nous renvoyons \`a \cite{ParPau11CJM} pour d'autres applications
arithm\'etiques, o\`u nous donnons un \'equivalent pr\'ecis du nombre de
repr\'esentations propres non \'equivalentes d'entieres rationnels de
valeur absolue au plus $t$ par une forme hermitienne binaire enti\`ere
ind\'efinie fix\'ee.

\medskip Nous commen\c{c}ons cet article (partie
\ref{sect:comptageconvcocomp}) par un r\'esultat de comptage en
courbure n\'egative, dont les hypoth\`eses sur l'espace ambiant sont
bien plus g\'en\'erales que celles du th\'eor\`eme
\ref{theo:introcomptage}, mais qui ne donne qu'un encadrement
asymptotique au lieu d'un \'equivalent (voir le th\'eor\`eme
\ref{theo:geomasymp}).  Nous d\'emontrons une version un peu plus
g\'en\'erale du th\'eor\`eme d'\'equidistribution
\ref{theo:introequidistrib} dans la partie \ref{sect:equidistrib},
ainsi que comment d\'eduire, d'un r\'esultat d'\'equidistribution au
niveau du fibr\'e tangent unitaire, un r\'esultat d'\'equidistribution
dans la vari\'et\'e, y compris en courbure variable~: les courbures
principales des horosph\`eres n'\'etant pas forc\'ement constantes,
des termes correcteurs remarquables apparaissent (voir le th\'eor\`eme
\ref{theo:equidistributionenbas}).  Dans la partie
\ref{sect:comptage}, nous d\'emontrons, en utilisant les r\'esultats
d'\'equidistribution pr\'ec\'edents, le r\'esultat de comptage d'arcs
g\'eod\'esiques, et diverses variantes pour les orbifolds et les
\'enonc\'es de cas particuliers utiles pour les applications
arithm\'etiques.

Bien s\^ur, les r\'esultats \ref{coro:introformes} et
\ref{coro:introcompirrquad} ne sont pas ind\'ependants: il existe une
correspondance bien connue entre les r\'eels irrationnels quadratiques
et les formes quadratiques binaires ind\'efinies primitives ne
repr\'esentant pas $0$. Voir la partie \ref{subsec:justifcomplexiteh}
pour des rappels sur cette correspondance, ainsi que pour diverses
interpr\'etations alg\'ebriques de notre complexit\'e
$|\alpha-\alpha^\sigma|^{-1}$.

\bigskip {\small {\it Remerciements :} Nous remercions Pierre Pansu
pour son aide dans la preuve des propositions \ref{prop:pansu} et
\ref{prop:pansujacobi}, et K.~Belabas pour des dicussions fructueuses
concernant les applications arithm\'etiques.}

\section{Sur le comptage des classes \`a droite d'un 
sous-groupe con\-vexe-cocompact}
\label{sect:comptageconvcocomp}

Nous renvoyons par exemple \`a \cite{Bourdon95,BriHae98} pour tout
rappel concernant les espaces ${\rm CAT}(-1)$ et leurs groupes
d'isom\'etries. Nous reprenons pour le th\'eor\`eme
\ref{theo:geomasymp} le cadre de travail de
\cite[Sect.~3]{ParPau10MZ}, en modifiant quelques hypoth\`eses, et en
nous inspirant de \cite[Theo.~1.1]{BelHerPau01} et
\cite[Theo.~3.4]{HerPau04} pour les preuves. Si $f$ et $g$ sont deux
applications d'un ensemble $E$ dans $\RR$, alors nous noterons
$f\asymp g$ s'il existe une constante $c\geq 1$ telle que
$\frac{1}{c}\, f(t)\leq g(t)\leq c\,f(t)$ pour tout $t$ dans $E$.

Soient $X$ un espace m\'etrique g\'eod\'esique ${\rm CAT}(-1)$ propre,
de bord \`a l'infini $\partial_\infty X$, et $\Ga$ un groupe discret
d'isom\'etries de $X$, d'exposant critique not\'e $\delta=\delta_\Ga$.
Rappelons que, pour tous $x,y$ fix\'es dans $X$, si $f_{\Ga,x,y}(t)=
\operatorname{Card}\{\ga\in\Ga\;:\;d(x,\ga y)\leq t\}$, alors
$$
\delta_\Ga=\limsup_{t\ra +\infty} \frac{1}{t} \log f_{\Ga,x,y}(t)\;.
$$ 
Rappelons que la condition $f_{\Ga,x,y}(t)\asymp e^{\delta t}$ ne
d\'epend pas des points $x,y$, et qu'elle est v\'erifi\'ee par exemple
si $X/\Ga$ est une variété riemanienne de volume fini
\cite{Margulis69,Bourdon95}, voir
\cite{Patterson76,Sullivan79,Roblin03} pour de nombreux cas plus
g\'en\'eraux.

Notons $\Lambda\Ga\subset \partial_\infty X$ l'ensemble limite de
$\Ga$, et, lorsque $\Lambda\Ga$ contient au moins deux points, notons
$\C\Ga\subset X$ l'enveloppe convexe de $\Lambda\Ga$.  Soit $\Ga_0$ un
sous-groupe {\it convexe-cocompact} (c'est-\`a-dire tel que
$\Lambda\Ga_0$ contienne au moins deux points et $\C\Ga_0/\Ga_0$ soit
compact), d'indice infini dans $\Ga$. Notons $\delta_0$ son exposant
critique, qui est strictement inf\'erieur \`a $\delta$ (ceci est bien
connu, voir par exemple \cite[page 116]{DalOtaPei00}), et
$C_0=\C\Gamma_0$.

Soit $C_\infty$ un convexe ferm\'e non vide de $X$, de stabilisateur
dans $\Ga$ not\'e $\Ga_\infty$. Notons $\delta_\infty$ l'exposant
critique de $\Ga_\infty$.

Supposons que $C_\infty$ ne contienne pas $\C\Ga$, que
$\Ga_\infty\backslash \partial C_\infty$ soit compact, et que
l'intersection de $C_\infty$ et de $\ga C_0$ soit non vide seulement
pour un ensemble fini de classes $[\ga]$ dans $\Ga_\infty\backslash
\Ga/\Ga_0$. Pour toute double classe $[\ga]$ dans
$\Ga_\infty\backslash \Ga/\Ga_0$, notons
$$
D([\ga])=d(C_\infty,\ga C_0)\;,
$$
ce qui ne d\'epend pas du repr\'esentant choisi.

Pour un exemple, on pourra prendre pour $X$ un espace sym\'etrique de
courbure sectionnelle au plus $-1$, pour $\Ga$ un sous-groupe
d'isom\'etries de $X$ de covolume fini, pour $\Ga_0$ le stabilisateur
dans $\Ga$ de l'axe de translation dans $X$ d'un \'el\'ement
hyperbolique de $\Ga$, et pour $C_\infty$ un point de $X$ ou une
horoboule de $X$ centr\'ee en un point fixe parabolique de
$\Ga$. Rappelons qu'alors $\delta_\infty < \delta$, et nous renvoyons
\`a \cite{Roblin03} pour de nombreux autres exemples o\`u cette
in\'egalit\'e est v\'erifi\'ee.

\btheo \label{theo:geomasymp} Soient $(X,\Ga,\Ga_0,C_\infty)$ comme
ci-dessus, si $f_{\Ga,x,y}(t)\asymp e^{\delta t}$ et si
$\delta_\infty<\delta$, alors
$$
\operatorname{Card} \{r\in \Ga_\infty\backslash
\Ga/\Ga_0\;:\;D(r)\leq t\} \;\asymp \;e^{\delta t}\;.
$$ 
\etheo

\dem 
Soient $x_\infty\in \partial C_\infty$, $x_0\in C_0$, $F$ le
sous-ensemble fini des $[\ga]$ dans $\Ga_\infty\backslash \Ga/\Ga_0$
tels que $C_\infty\cap\ga C_0$ soit non vide, et $\R=\Ga_\infty
\backslash \Ga/\Ga_0-F$.  Pour tout $r$ dans $\R$, soit $\ga_r$ un
repr\'esentant de la double classe $r$ tel que
$$
d(x_\infty, \ga_r x_0)=
\min_{\alpha\in \Ga_\infty, \beta\in \Ga_0} 
d(x_\infty, \alpha\ga_r \beta x_0)\;.
$$ 
Ce minimum existe, ainsi que 
$$
\eta=\min_{r\in \R} d(x_\infty, \ga_r x_0)>0\;,
$$ 
car l'action de $\Ga$ sur l'espace propre $X$ est proprement
discontinue.

\blemm \label{lem:perpcomm} Il existe $c>0$ tel que pour tout
$r\in\R$, le segment g\'eod\'esique $[x_\infty, \ga_r x_0]$ est \`a
distance de Hausdorff au plus $c$ du segment perpendiculaire commun
entre $C_\infty$ et $\ga_r C_0$.  
\elemm

\dem 
Les espaces m\'etriques quotients $\Ga_\infty\backslash \partial
C_\infty$ et $\Ga_0 \backslash C_0$ sont compacts, par les
hypoth\`eses. Notons $c'>0$ le maximum de leurs diam\`etres.

Pour tout $\ga\in\Ga$ tel que $C_\infty\cap\ga C_0=\emptyset$, notons
$[u_\ga,v_\ga]$ le segment perpendiculaire commun entre $C_\infty$ et
$\ga C_0$, avec $u_\ga\in \partial C_\infty$. Remarquons que, pour
tous $\alpha\in\Ga_\infty$ et $\beta\in\Ga_0$, nous avons alors
$u_{\alpha\ga\beta}=\alpha u_\ga$, $v_{\alpha\ga\beta}=\alpha v_\ga$,
et
$$
d(x_\infty,\alpha\ga\beta x_0)=d(\alpha^{-1}x_\infty,\ga\beta x_0)
\leq
d(\alpha^{-1}x_\infty,u_\ga)+d(u_\ga,v_\ga)+d(v_\ga,\ga\beta x_0)\;.
$$ 
Par d\'efinition de $c'$, pour tout $[\ga]\in\R$, quitte \`a
multiplier $\ga$ \`a gauche par un \'el\'ement de $\Ga_\infty$ et \`a
droite par un \'el\'ement de $\Ga_0$, ce qui ne change pas
$d(u_{\ga_r},v_{\ga_r})$, les distances $d(x_\infty,u_\ga)$ et
$d(v_\ga,\ga x_0)$ peuvent \^etre rendues au plus $c'$. Donc
$d(x_\infty, \ga_r x_0) \leq 2c'+ d(u_{\ga_r},v_{\ga_r})$ pour tout
$r\in\R$.

Puisque $X$ est CAT$(-1)$, par convexit\'e et disjoinction de
$C_\infty$ et de $C_0$, il existe une constante $c''=c''(\eta)$ (ne
d\'ependant que de $\eta$) telle que pour tout $[\ga]\in \R$, nous
ayons
\begin{equation}\label{eq:quasigeod}
d(x_\infty, \ga x_0)\geq d(x_\infty, u_{\ga})+
d(u_{\ga},v_{\ga})+ d(v_{\ga},\ga x_0)-c''\;.
\end{equation}
En posant $c=2c'+c''$, nous avons $d(x_\infty, u_{\ga_r})\leq c$ et
$d(v_{\ga_r}, \ga_r x_0)\leq c$ pour tout $r\in\R$. Le r\'esultat en
d\'ecoule, par convexit\'e.  
\cqfd

\medskip
Par l'in\'egalit\'e triangulaire, il d\'ecoule de ce lemme que
$$
D(r)=d(C_\infty,\ga_r C_0)\geq d(x_\infty,\ga_r x_0)-2c
$$ 
pour tout $r\in\R$. Puisque $f_{\Ga,x_\infty,x_0}(n) \asymp e^{\delta
  n}$, il existe donc une constante $c_+>0$ telle que
$$
\operatorname{Card} \{r\in \R\;:\;D(r)\leq n\} \leq c_+e^{\delta n}
$$ 
pour tout $n$ dans $\NN$. Pour d\'emontrer une minoration analogue
du terme de gauche, nous commen\c{c}ons par un lemme \'el\'ementaire.

\blemm \label{lem:taupin} Pour tous $M\in\NN$ et $A,\delta,\delta'>0$
o\`u $\delta>\delta'$, il existe $B>0$ tel que pour toutes suites
strictement positives $(a_n)_{n\in\NN}, (b_n)_{n\in\NN},
(c_n)_{n\in\NN}$ v\'erifiant, pour tout $n$ dans $\NN$, que
$b_n,c_n\leq Ae^{\delta'n}$, $a_n\leq Ae^{\delta n}$ et
$\sum_{k+\ell+m=n}^{n+M} a_kb_\ell c_{m}\geq \frac{1}{A} e^{\delta
  n}$, on ait $\sum_{k=0}^{n}a_{k}\geq B e^{\delta n}$.  
\elemm

\dem Les s\'eries positives $\sum e^{(\delta'-\delta)k}$ et $\sum
k\;e^{(\delta'-\delta)k}$ sont convergentes, notons $S$ et $S'$ leur
somme. Soient $n,N\geq 1$ des entiers. Alors
\begin{align*}
\sum_{k=0}^n\;\sum_{\ell+m=n+N-k}^{n+N-k+M} a_k\,b_\ell \,c_{m}&\leq 
\;\sum_{k=0}^n\;\sum_{\ell+m=n+N-k}^{n+N-k+M} A\,e^{\delta
  k}A\,e^{\delta'\ell}A\,e^{\delta'm}\\&
\leq (M+1)A^3\sum_{k=0}^n\;(n+N-k+M+1)\,e^{\delta
  k+\delta'(n+N-k+M)}\\&
\leq (M+1)A^3 e^{\delta' (N+M)}(S(N+M+1)+S')\;e^{\delta n}\;.
\end{align*}
De plus, en posant $B'=(M+1)(N+M)A^2e^{\delta'(N+M)}>0$, nous avons
\begin{align*}
B'\sum_{k=n+1}^{n+N+M} a_k&\geq
\sum_{k=n+1}^{n+N+M}\;\sum_{\ell+m=n+N-k}^{n+N-k+M} a_k\,b_\ell \,c_{m}\\&
=\sum_{k+\ell+m=n+N}^{n+N+M} a_kb_\ell c_{m}-\sum_{k=0}^n\;
\sum_{\ell+m=n+N-k}^{n+N-k+M}a_k\,b_\ell \,c_{m}\\
&\geq\frac{1}{A}\,e^{\delta(n+N)}-(M+1)A^3 
e^{\delta' (N+M)}(S(N+M+1)+S')\;e^{\delta n}\;.
\end{align*}
Fixons $N$ assez grand tel que $B''=\frac{1}{A}e^{\delta N}-(M+1)A^3
e^{\delta' (N+M)}(S(N+M+1)+S')>0$, ce qui est possible car
$\delta>\delta'$. Alors pour tout $n\geq N+M$,
$$ 
\sum_{k=0}^{n}a_{k}\geq \frac{B'' e^{\delta (n-N-M)}}{B'}\;.
$$
Le r\'esultat en d\'ecoule.
\cqfd

\medskip
Montrons maintenant qu'il existe $c_->0$ tel que 
$$
\operatorname{Card} \{r\in \R\;:\;D(r)\leq n\} \geq c_-e^{\delta n}
$$ 
pour tout $n$ dans $\NN$. Soit $\delta'>0$ tel que
$\max\{\delta_\infty,\delta_0\} <\delta' <\delta$. En particulier,
$$
b_n=\operatorname{Card} \{\alpha\in\Ga_\infty\;:\;n\leq
d(x_\infty,\alpha x_\infty)< n+1\}
$$ 
et 
$$
c_n=\operatorname{Card}
\{\beta\in\Ga_0\;:\;n\leq d(x_0,\beta x_0)< n+1\}
$$ 
sont des $o(e^{\delta'n})$, par d\'efinition d'un exposant critique.

Puisque $f_{\Ga,x_\infty,x_0}(n) \asymp e^{\delta n}$ et puisque $F$
est fini, il existe donc, par un argument de s\'erie g\'eom\'etrique,
$A'>0$ et $T>0$ tels que, pour tout $t$ dans $\RR$,
$$
g(t)=\operatorname{Card} \{(\alpha,\beta,r)\in \Ga_\infty \times
\Ga_0\times \R\;:\;t\leq d(x_\infty,\alpha\ga_r\beta x_0)\leq t+T\}
\geq A' e^{\delta t}\;.
$$
Soit $(\alpha,\beta,r)\in\Ga_\infty\times\Ga_0\times \R$. Par le lemme
\ref{lem:perpcomm}, nous avons
\begin{align*}
d(x_\infty,\alpha\ga_r\beta x_0)&\leq
d(\alpha^{-1}x_\infty,x_\infty)+d(x_\infty,\ga_r x_0)+d(x_0,\beta
x_0)\\&\leq d(x_\infty,\alpha x_\infty)+D(r)+d(x_0,\beta
x_0)+2c
\end{align*}
et, par la formule \eqref{eq:quasigeod} puis par le lemme
\ref{lem:perpcomm},
\begin{align*}
d(x_\infty,\alpha\ga_r\beta x_0)&\geq
d(x_\infty,\alpha u_{\ga_r})+d(u_{\ga_r},v_{\ga_r})+d(v_{\ga_r},
\ga_r\beta x_0)-c''\\&\geq d(x_\infty,\alpha x_\infty)+D(r)+
d(x_0,\beta x_0)-2c-c''\;.
\end{align*}
Posons 
$$
a_n=\operatorname{Card} \{r\in\R\;: \;n\leq D(r)< n+1\}\leq 
c_+e^{\delta(n+1)}\;.
$$ 
Si $n+2c+3\leq d(x_\infty,\alpha\ga_r\beta x_0)\leq n +2c+3 +T$, nous
avons donc, en notant les parties enti\`eres $k=E[D(r)]$,
$\ell=E[d(x_\infty,\alpha x_\infty)]$, $m=E[d(x_0,\beta x_0)]$,
$$
k+l+m\geq D(r)+d(x_\infty,\alpha x_\infty)+d(x_0,\beta
x_0)-3\geq d(x_\infty,\alpha\ga_r\beta x_0)-2c-3\geq n
$$
et
$$
k+l+m\leq D(r)+d(x_\infty,\alpha x_\infty)+d(x_0,\beta
x_0)\leq d(x_\infty,\alpha\ga_r\beta x_0)+2c+c''\leq n +M\;,
$$ 
o\`u $M=E[4c+c''+3+T]+1$. Donc par un argument de comptage, pour tout
$n\in\NN$, nous avons
$$
g(n+2c+3)\leq \sum_{k+\ell+m=n}^{n+M} a_kb_\ell c_{m}\;.
$$ 
Ceci entra\^ine que $\sum_{k+\ell+m=n}^{n+M} a_k b_\ell c_{m}\geq
A'e^{\delta (n+2c+3)}$ pour tout $n\in\NN$. Le r\'esultat d\'ecoule
alors du lemme \ref{lem:taupin}, par d\'efinition de $a_n$. 
 \cqfd

 \medskip Sous des hypoth\`eses bien plus restrictives sur
 $(X,\Gamma)$, le th\'eor\`eme \ref{theo:geomasymp} peut \^etre
 am\'elior\'e. Nous raisonnerons comme dans \cite{EskMcMul93} en trois
 \'etapes, en montrant un r\'esultat d'\'equidistribution, un
 r\'esultat de comptage asymptotique, et l'am\'elioration r\'esultera
 alors de calculs de volumes de voisinages tubulaires.

\section{\'Equidistribution des hypersurfaces \'equidistantes}
\label{sect:equidistrib}

Soit $M$ une vari\'et\'e riemannienne lisse compl\`ete de dimension
$n\geq 2$. Nous noterons $\V_\epsilon A$ 
le $\epsilon$-voisinage ouvert 
d'une partie $A$ de $M$.

Nous noterons $d\vol_M$ la mesure riemannienne de $M$ et $\Vol(M)$ son
volume total. Lorsque $\Vol(M)$ est fini, nous appellerons {\it
  moyenne riemannienne } sur $M$ la mesure de probabilit\'e sur $M$
\'egale \`a $\frac{1}{\Vol(M)} \;d\vol_M$. Nous noterons
$\oint\varphi\; d\vol_M$ la moyenne riemannienne d'une application
$\varphi\in \C_c(M)$ (\`a valeurs complexes, continue \`a support
compact).  Pour toute application mesurable $f:M\ra \,[0,+\infty[\,$,
lorsque l'int\'egrale $\int_{M} f \;d\vol_{M}$ est finie non nulle,
nous appellerons {\it $f$-moyenne} de $M$ la mesure de probabilit\'e
$\frac{f}{\int_{M}f\, d\vol_{M}}\;d\vol_{M}$ sur $M$.

Notons $\pi: T^1M \ra M$ le fibr\'e tangent unitaire de $M$ et
$(g^t)_{t\in\RR}$ le flot g\'eod\'esique sur $T^1M$.  Pour tout $v\in
T^1M$, notons $\ga_v(t)$ la géodésique de $M$ telle que
$\dot{\ga}_v(0)=v$. Munissons $T^1M$ de la m\'etrique riemannienne
usuelle (dite de Sasaki), qui fait de $\pi$ une submersion
riemannienne.  Rappelons que la mesure riemannienne $d\vol_{T^1M}$,
appel\'ee {\it mesure de Liouville}, se d\'esint\`egre par $\pi$ en
$$
d\vol_{T^1M}=\int_{x\in M} d\vol_{T^1_xM}\;d\vol_{M}(x)\;,
$$
o\`u $d\vol_{T^1_xM}$ est la mesure sph\'erique sur la fibre
$T^1_xM$ de $\pi$ au-dessus de $x\in M$. En particulier,
\begin{equation}\label{eq:volTunMM}
\Vol(T^1M)=\Vol(\SSS_{n-1})\Vol(M)\;.
\end{equation} 
Nous appellerons {\it moyenne de Liouville} la mesure de Liouville
normalis\'ee pour \^etre de probabilit\'e.

\medskip Soit $\iota:M'\ra M$ une immersion lisse d'une vari\'et\'e
lisse $M'$ dans $M$. Nous munirons $M'$ de sa structure de vari\'et\'e
riemannienne lisse image r\'eciproque par $\iota$. Nous noterons
encore par abus $d\vol_{M'}$ la mesure sur $M$ image par l'immersion
$\iota$ de la mesure riemannienne de $M'$ (et de même
pour les moyennes et $f'$-moyennes o\`u $f':M'\ra \,[0,+\infty[\,$ est
mesurable).

Nous noterons $\pi_\perp:\nu^1M'\ra M'$ le fibr\'e normal unitaire de
$M'$ dans $M$ (pour tout $x\in M'$, nous avons $\nu_x^1M'=
T^1_{\iota(x)}M \cap \big(T\iota(T_xM')\big)^\perp$) et $\iota_\perp:
\nu^1M'\ra T^1M$ l'immersion naturelle. Si $M'$ est connexe, de
codimension $p$ dans la composante connexe de $M$ contenant $M'$,
alors $\nu^1M'$ est connexe si et seulement si $p\neq 1$. De plus, si
$p=1$, alors $\nu^1M'$ admet exactement deux composantes
connexes. Pour tout ouvert $\Sigma$ de $\nu^1M'$, nous noterons par
abus $g^t\Sigma$ la sous-vari\'et\'e lisse immerg\'ee de $T^1M$,
pouss\'ee, par le flot g\'eod\'esique au temps $t$, de l'image de
$\Sigma$ dans $T^1M$ par $\iota_\perp$.

Si $M'$ est une hypersurface immerg\'ee lisse strictement convexe de
$M$, nous noterons $N_+:M'\ra T^1M$ (resp.~$N_-:M'\ra T^1M$) {\it
  l'application de Gauss} sortante (resp.~rentrante) de $M'$, qui \`a
$x\in M'$ associe le vecteur normal sortant (respectivement rentrant)
\`a $M'$ en l'image dans $M$ de $x$, et $N_+M'$ (respectivement
$N_-M'$) la sous-vari\'et\'e immerg\'ee dans $T^1M$, image de
l'application de Gauss correspondante.

\bigskip
Nous commen\c{c}ons par un lemme de courbure n\'egative.

\blemm\label{lem:projection} Soit $X$ une vari\'et\'e riemannienne
lisse compl\`ete simplement connexe, \`a courbure sectionnelle au plus
$-1$.  Soit $K$ un convexe fermé de $X$ de bord $\partial K$ une
sous-vari\'et\'e de classe ${\rm C}^1$. Soit $v$ un vecteur tangent
unitaire d'origine un point de $\partial K$ et d'angle au plus
$\epsilon\in[0,\frac{\pi}{2}[$ avec une normale sortante en ce
point. Alors, pour $t>0$, le point $y\in\partial K$ le plus proche de
$\ga_v(t)$ est \`a distance au plus $\artanh(\sin\epsilon)$ du point
$x=\ga_v(0)$. De plus, si $t'=d(y,\ga_v(t))$, si $v'$ est le vecteur
tangent unitaire normal sortant de $K$ en $y$, alors l'angle $\alpha$
entre les vecteurs $\dot{\ga_{v'}}(t')$ et $\dot{\ga_{v}}(t)$ au point
$\ga_v(t)$ v\'erifie $\sin\alpha\leq \frac{\sinh(\artanh(\sin
  \epsilon))}{\sinh t}$.  \elemm

\dem Supposons que $x\ne y$ et notons $z= \ga_v(t)$. Considérons le
triangle $\ov x,\ov y,\ov z$ de comparaison dans $\hdr$ du triangle de
$X$ de sommets $x,y,z$, dont l'angle $\theta$ en $\ov x$ est au moins
$\frac{\pi}{2}-\epsilon$ et celui en $\ov y$ au moins $\frac{\pi}{2}$.
Notons $L$ la g\'eod\'esique de $\hdr$ passant par $\ov x$ et $\ov
y$. Alors, le point $\ov y'\in L$ le plus proche de $\ov z$ satisfait
$d(x,y)=d(\ov x,\ov y)\leq d(\ov x,\ov y')=\ell$. Notons $t'=d(\ov
z,\ov y')$. Alors par les formules des triangles hyperboliques
\cite[p.~145]{Beardon83}, nous avons $\cosh t=\cosh t'\cosh\ell$ et
$$
\cosh t'=\cosh t \cosh \ell -\sinh t\sinh \ell\cos\theta\geq \cosh
t(\cosh \ell-\sinh \ell \sin\epsilon)\;.
$$
Donc $\cosh\ell(\cosh \ell-\sinh \ell \sin\epsilon)\leq 1$ et
$\ell\leq \artanh(\sin\epsilon)$, d'où le premier résultat.  Le second
en d\'ecoule, toujours par comparaison, par la formule des sinus
hyperboliques \cite[p.~148]{Beardon83}.  \cqfd

\medskip Le r\'esultat suivant fournit, dans le fibr\'e unitaire
tangent, un th\'eor\`eme d'\'equidistribution, quand $t$ tend vers
$+\infty$, de l'image par le flot g\'eod\'esique au temps $t$ du fibré
normal unitaire d'une sous-vari\'et\'e immerg\'ee de volume fini
totalement g\'eod\'esique.

\btheo \label{theo:equidistributionenhaut} Soit $M$ une vari\'et\'e
hyperbolique connexe de volume fini, de dimension au moins $2$. Soit
$C$ une sous-vari\'et\'e immerg\'ee de $M$, connexe, totalement
géodésique et de volume fini. Soit $\Sigma$ une union de composantes
connexes du fibré normal unitaire de $C$. Alors la moyenne
riemannienne sur $g^t \Sigma$ converge vaguement vers la moyenne de
Liouville quand $t\ra+\infty$~:
$$
\forall\;\varphi\in \C_c(T^1M),\;\;\;\;
\lim_{t\ra+\infty}\;\oint\; \varphi\;\;d\vol_{g^t \Sigma}\;=
\;\oint\; \varphi\;\;d\vol_{T^1M}
\;.
$$  
\etheo

\dem 
Le volume de la sous-variété immergée $g^t \Sigma$ de $T^1M$, pour la
métrique riemannienne induite, est fini, car $C$ est de volume fini et
l'application canonique $g^t \Sigma\ra C$ contracte l'\'el\'ement de
volume d'au plus une constante ne d\'ependant que de $t,k,n$, par
homog\'en\'eit\'e. Nous pouvons donc consid\'erer la moyenne
riemannienne sur $g^t \Sigma$.

Soient $\wt M\ra M$ et $\wt C\ra C$ des revêtements universels
riemanniens.  Le relevé $\wt \iota:\wt C\ra \wt M$ de l'immersion
$\iota:C\ra M$ est un plongement lisse sur un sous-espace hyperbolique
encore noté $\wt C$ de l'espace hyperbolique $\wt M$. Notons $G_{\wt
  C}$ le groupe des isométries de $\wt M$ préservant $\wt C$. Il agit
transitivement sur le fibré normal unitaire $\nu^1\wt C$, et commute
avec l'action du flot géodésique sur $T^1\wt M$.

Pour tout $\eta\in\, ]0,\frac{\pi}{2}[ \,$, notons $\wt V_\eta$
l'ouvert $G_{\wt C}$-invariant de $T^1\wt M$ formé des vecteurs $w\in
T^1\wt M$ tels qu'il existe $s\in\,]0,\eta[$ et $v\in \nu^1\wt C$ tels
que les vecteurs $w$ et $\dot{\ga}_v(s)$ aient même origine et soient
d'angle strictement inférieur à $\eta$. L'application $\wt P:w\mapsto
v$ ainsi définie est une fibration lisse $G_{\wt C}$-équivariante de
$\wt V_\eta$ sur $\nu^1\wt C$. Notons $\Ga_C$ l'intersection de
$G_{\wt C}$ et du groupe de revêtement de $\wt M\ra M$. Par passage au
quotient, $\wt P$ induit une fibration $P$ de $\Ga_C\bs \wt V_{\eta}$
sur $\nu^1C$.  Notons $\wt U_\eta$ la réunion de composantes connexes
de $\wt V_\eta$ telle que $U_\eta= \Ga_C\bs \wt U_{\eta}$ soit
l'ouvert immerg\'e $P^{-1}(\Sigma)$ de $T^1M$. La mesure riemannienne
de $\wt U_\eta$ se désintègre par $\wt P$ par rapport à la mesure
riemannienne de $\nu^1\wt C$, de mesures conditionnelles sur les
fibres finies et invariantes par $G_{\wt C}$. Notons
$\chi_{U_\eta}:T^1M\ra \NN$ la fonction caract\'eristique (comptant
avec multiplicit\'e) de l'image de $U_\eta$ dans $T^1M$.
Puisque $C$ est de volume fini, l'intégrale $\int_{T^1M}
\chi_{U_\eta}(v)\,d \vol_{T^1M}(v)= \vol(U_\eta)$ est finie.

Montrons que pour tout $\eta'>0$, il existe $\eta>0$ tel que pour tout
$t>0$,
$$
g^tU_{\eta}\subset \V_{\eta'}g^t\Sigma\;.
$$ 
En effet, par passage au quotient, il suffit de démontrer que $g^t\wt
U_{\eta}\subset \V_{\eta'}g^t\wt\Sigma$, où $\wt\Sigma$ est la
préimage de $\Sigma$ dans $\nu^1\wt C$.  Soient $w\in \wt U_{\eta}$,
$s\in\,]0,\eta[$ et $v\in \wt\Sigma$ tels que $\angle\big( w,
\dot{\ga}_v(s)\big)<\eta$. Notons $K$ l'adhérence de $\V_s\wt C$, et
$v'\in \nu^1\wt C$ le vecteur tangent en l'origine du rayon géodésique
normal à $\wt C$ passant par $\ga_w(t)$, qui passe à l'instant $s$ par
le point $y$ de $\partial K$ le plus proche de $\ga_w(t)$. Remarquons
que $v'$ appartient aussi à $\wt \Sigma$. Par le lemme
\ref{lem:projection}, si $\eta$ est suffisamment petit, alors $g^{t}w$
fait un angle petit avec $g^{t'+s} v'$ o\`u $|t-t'-s|\leq
d(y,\ga_v(0))+\eta$ est petit, uniform\'ement en $t$ et $w$, donc est
uniform\'ement pr\`es de $g^{t}\Sigma$.

Soit $\varphi$ une fonction continue \`a support compact sur $T^1M$,
qui est donc uniform\'ement continue. Si $\eta$ est assez petit, par
homogénéité, la moyenne de $\varphi$ sur $g^t \Sigma$ est donc proche
de sa moyenne sur l'ouvert immerg\'e $g^{t}U_{\eta}$, uniform\'ement
en $t>0$. Cette derni\`ere moyenne est \'egale \`a
$$
\frac{\int_{T^1M} \chi_{U_{\eta}}(v)\varphi(g^{-t}v)\,
  d\vol_{T^1M}(v)}{\int_{T^1M} \chi_{U_{\eta}}(v)\,d \vol_{T^1M}(v)}
$$ 
par invariance de la mesure de Liouville par le flot
g\'eod\'esique. Elle converge donc quand $t$ tend vers $+\infty$, par
la propri\'et\'e de m\'elange du flot g\'eod\'esique pour la mesure de
Liouville, vers
$$
\frac{\int_{T^1M} \varphi(v)\, 
d \vol_{T^1M}(v)}{\operatorname{Vol}(T^1M)}\;,
$$ 
qui est la moyenne riemannienne de $\varphi$ sur $T^1M$.  \cqfd

\bigskip Le th\'eor\`eme pr\'ec\'edent est un r\'esultat
d'\'equidistribution au niveau du fibr\'e unitaire tangent
$T^1M$. Lorsque l'on veut en d\'eduire des propri\'et\'es
d'\'equidistribution dans la vari\'et\'e $M$ elle-m\^eme, des
ph\'enom\`enes surprenants apparaissent, lorsque la courbure est
variable.

Nous n'avons pas trouv\'e de r\'ef\'erence pour le r\'esultat
g\'en\'eral ci-dessous de g\'eom\'etrie riemannienne, connu des
sp\'ecialistes des sous-vari\'et\'es. Il g\'en\'eralise la remarque
suivante~: le cercle de centre $0$ et de rayon $r$ dans le plan
euclidien $\RR^2$ s'envoie par la normale sortante dans
$T\RR^2=\RR^2\times\RR^2$ par l'application $\psi:re^{i\theta}\mapsto
(re^{i\theta},e^{i\theta})$ et l'\'el\'ement de longueur
$\sqrt{r^2+1}\;d\theta$ sur l'image du cercle diff\`ere de l'image par
$\psi$ de l'\'el\'ement de longueur $r\,d\theta$ du cercle par
$\sqrt{1+\kappa^2}$ o\`u $\kappa=1/r$ est la courbure du cercle.

\bprop\label{prop:pansu} Soient $M$ une vari\'et\'e riemannienne lisse
de dimension $n+1$, $L$ une hypersurface lisse de $M$, $N$ un
champ de vecteurs unitaires lisse orthogonal le long de $L$,
$\mu$ la mesure riemannienne sur $L$, et $\nu$ la mesure
riemannienne sur l'image du plongement lisse de Gauss $\psi:L\ra
T^1M$ (d\'efini par $\psi(x)=N(x)$ pour tout $x$ dans $L$).
Alors $\psi_* \mu$ et $\nu$ sont absolument continues l'une par
rapport \`a l'autre et
$$
\forall x\in L\;,\;\;\;\frac{d\nu}{d\psi_*\mu}(\psi(x))=
\prod_{i=1}^n\big(1+\kappa_i(x)^2\big)^{1/2}\;,
$$
o\`u $\kappa_1(x),\dots,\kappa_n(x)$ sont les courbures principales de
l'hypersurface $L$ en $x$.  
\eprop

Rappelons que les courbures principales sont bien d\'efinies \`a
l'ordre pr\`es, et qu'elles sont chang\'ees en leurs oppos\'es si $N$
est remplacé par $-N$. Le produit ci-dessus est bien
invariant par permutation et par changement de signe des
$\kappa_i(x)$. L'absolue continuit\'e des mesures est imm\'ediate,
elles sont dans la classe de Lebesgue de la sous-vari\'et\'e lisse
$\psi(L)$.

\medskip \dem Le probl\`eme \'etant local, nous pouvons supposer que
$M$ est un voisinage suffisamment petit d'un point de $L$. Le
champ de vecteurs unitaires $N$ le long de $L$ s'\'etend donc en un champ
de vecteurs unitaires lisse sur $M$, que nous noterons encore $N:M\ra
TM$. Notons $p:TM\ra M$ la projection canonique et $\nabla$ la
connection riemannienne de $M$. Rappelons que le fibr\'e tangent de
$TM$ est somme directe $TTM=V\oplus H$ de deux sous-fibr\'es, dits
vertical et horizontal, orthogonaux pour la m\'etrique riemannienne de
$TM$, tels que, si $\pi_V:TTM\ra V$ est la projection orthogonale sur
le facteur vertical, alors pour tout $v\in TM$,

$\bullet$~~ $Tp_{\mid H_v}:H_v\ra T_{p(v)}M$, la restriction \`a $H_v$
de l'application tangente de $p$, est un isomorphisme lin\'eaire
isom\'etrique,

$\bullet$~~ $V_v=\operatorname{Ker} T_vp=T_v(T_{p(v)}M)=T_{p(v)}M$, et
le produit scalaire induit sur $V_v$ par la m\'etrique riemannienne de
$TM$ est \'egal au produit scalaire sur $T_{p(v)}M$,

$\bullet$~~pour tout champ de vecteurs $X:M\ra TM$ sur $M$, 
$$
\nabla_v X=\pi_V\circ TX(v)\;.
$$
Rappelons que l'op\'erateur de seconde forme fondamentale (pour le
choix de la normale $N$) de $L$ est l'endomorphisme $S$ du fibr\'e
$TL$ d\'efini (avec des conventions de signes variant suivant les
références) par $S(v)=\nabla_v N$ pour tout $v\in TL$. Comme $\psi$
est la restriction de $N$ \`a $L$, et $T\psi$ celle de $TN$ \`a $TL$,
nous avons donc $ \pi_V\circ T\psi=S$.  Comme $p\circ\psi={\rm id}_L$,
nous avons $Tp\circ T\psi={\rm id}_{TL}$.  Donc en notant $g^{TM}$ et
$g^L$ les tenseurs de m\'etriques riemanniennes sur $TM$ et $L$, nous
avons la formule bien connue
$$
\psi^*g^{TM}=g^L+S^*g^L\;.
$$
Par d\'efinition des courbures principales en un point $x$ de $L$,
l'op\'erateur $S$ est, sur $T_xL$, un endomorphisme sym\'etrique
$S(x)$ diagonalisable en base orthonorm\'ee (pour le produit scalaire
$g^L_x$) de valeurs propres \'egales aux courbures principales. Donc,
pour tout $x$ dans $L$, puisqu'en coordonnées locales, l'élément de
volume est la racine carrée du déterminant de la métrique,
$$
d\operatorname{vol}_{\psi^*g^{TM}}(x)= \sqrt{\det({\rm
    id}+S(x)^*S(x))}\;d\operatorname{vol}_{g^L}(x)=
\Big(\prod_{i=1}^n\big(1+\kappa_i(x)^2\big)\Big)^{1/2}
\;d\operatorname{vol}_{g^L}(x)
  \;.$$ 
Le r\'esultat en d\'ecoule.  
\cqfd

\medskip Soit $M$ une vari\'et\'e riemannienne lisse de dimension $n\geq
2$ et \`a courbure sectionnelle strictement n\'egative. Pour tout $v\in
T^1M$, notons $\F_{su}(v)$ le germe de feuille en $v$ du feuilletage
fortement instable de $T^1M$ pour le flot g\'eod\'esique. Notons
$\kappa_1(v),\dots, \kappa_{n-1}(v)$ les courbures principales en
$\pi(v)$ du germe d'hypersurface $\pi(\F_{su}(v))$ de $M$. Notons
$\kappa=\kappa_M:T^1M\ra\; ]0,+\infty[$ l'application d\'efinie par
$$
\kappa(v)=\prod_{i=1}^{n-1} \big(1+\kappa_i(v)^2\big)
^{\frac12}\;.
$$
Rappelons que si $\wt\pi:\widetilde M\ra M$ est un rev\^etement
riemannien universel de $M$, alors les feuilles du feuilletage
fortement instable de $T^1M$ sont les images par $T\wt \pi$ des
sous-vari\'et\'es $N_+H$ de $T^1\wt M$ lorsque $H$ parcourt l'ensemble
des horosph\`eres de $\widetilde M$. De plus, si $X$ est un espace
sym\'etrique \`a courbure strictement n\'egative, pour toute
horosph\`ere $H$ de $X$, le groupe des isom\'etries de $X$ fixant le
point \`a l'infini de $H$ agit transitivement sur $H$. Donc si $M$ est
localement sym\'etrique, alors l'application $\kappa:T^1M\ra\;
]0,+\infty[$ est constante. Il serait d'ailleurs int\'eressant de
savoir si les seules vari\'et\'es riemanniennes $M$ lisses compactes
\`a courbure sectionnelle strictement n\'egative, dont la fonction
$\kappa_M$ est constante, sont localement sym\'etriques.

\btheo\label{theo:equidistributionenbas} Soit $M$ une vari\'et\'e
riemannienne lisse compl\`ete, de volume fini, de dimension au moins
$2$ et de courbure sectionnelle pinc\'ee au plus $-1$. Soit $\Sigma'$
une hypersurface lisse proprement immerg\'ee strictement convexe de
volume fini, telle que la moyenne riemannienne sur $g^t N_+\Sigma'$
converge vaguement vers la moyenne de Liouville.  Pour toute fonction
$\varphi$ continue \`a support compact dans $M$, la $(\kappa\circ
N_+)$-moyenne de $\varphi$ sur $\pi(g^t N_+\Sigma')$ converge vers la
moyenne riemannienne de $\varphi$ sur $M$ quand $t\ra+\infty$.  
\etheo

Notons que lorsque $\kappa$ est constante (par exemple lorsque $M$ est
localement sym\'etrique), la $(\kappa\circ N_+)$-moyenne est \'egale
\`a la moyenne riemannienne; dans ce cas, la conclusion du théorème
dit que la moyenne riemannienne de $\varphi$ sur $\pi(g^t N_+\Sigma')$
converge vers la moyenne riemannienne de $\varphi$ sur $M$ quand
$t\ra+\infty$.  

Par le théorème \ref{theo:equidistributionenhaut}, le th\'eor\`eme
\ref{theo:introequidistrib} de l'introduction en d\'ecoule.  En effet,
si $C$ et $\Sigma$ vérifient les hypothèses du théorème
\ref{theo:equidistributionenhaut}, posons $\Sigma'=\pi(g^1\Sigma)$,
qui est une hypersurface immergée lisse strictement convexe de volume
fini.  Puisque $C$ est totalement géodésique, le rayon d'injectivité
d'un point de $C$ dans $C$ est supérieur ou égal à celui dans $M$.
Donc l'immersion de $C$ dans $M$ est propre. Par conséquent,
l'immersion lisse de $\Sigma'$ dans $T^1M$ est propre.  Le résultat
découle alors du fait que $g^t\Sigma=g^{t-1}N_+\Sigma'$ pour tout
$t\geq 1$.

\medskip \dem Pour tout $t>0$, notons $\Sigma(t)=\pi(g^t N_+\Sigma')$,
qui est une hypersurface lisse proprement immerg\'ee strictement
convexe de $M$ de volume fini. En particulier, les moyennes
riemanniennes sur $\Sigma(t)$ et $g^t N_+\Sigma'$ sont bien
d\'efinies. Pour tout $t>1$, notons $\kappa'_1(x),\dots,$
$\kappa'_{n-1}(x)$ les courbures principales en un point $x$ de
$\Sigma(t)$ (pour la normale sortante).  Posons alors
$$
\kappa'_t(x)=\prod_{i=1}^{n-1}\big(1+{\kappa'_i(x)}^2\big)^{1/2}\;.
$$

\medskip
\noindent{\bf \'Etape 1. } Montrons tout d'abord que, pour tout
$\varphi\in\C_c(M)$,  la $\kappa'_t$-moyenne de $\varphi$ sur
$\Sigma(t)$ converge vers la moyenne riemannienne de $\varphi$ sur $M$
quand $t\ra+\infty$.

\medskip Soit $\psi:x\mapsto N_+(x)$ l'application de Gauss de
$\Sigma(t)$ d\'efinie par la normale sortante $N_+$. Puisque
$\pi\circ\psi= \operatorname{id} _{\Sigma(t)}$, pour tout $x$ dans
$\Sigma(t)$, nous avons, par la proposition \ref{prop:pansu},
$$
\frac{d\pi_*\!\vol_{g^tN_+\Sigma'}}{d\vol_{\Sigma(t)}}(x)=
\frac{d\pi_*\!\vol_{g^tN_+\Sigma'}}{d\pi_*\psi_*\!\vol_{\Sigma(t)}}
(\pi\circ\psi(x))=
\frac{d\vol_{g^tN_+\Sigma'}}{d\psi_*\;\vol_{\Sigma(t)}}(\psi(x))=
\kappa'(x)\;.
$$
Donc, 
\begin{align*}
\int_{\Sigma(t)} \;\varphi\;\kappa'_t \;d\vol_{\Sigma(t)}& =
\int_{\Sigma(t)}
\;\varphi\;\frac{d\pi_*\!\vol_{g^tN_+\Sigma'}}{d\vol_{\Sigma(t)}}
\;d\vol_{\Sigma(t)}=
\int_{\Sigma(t)}
\;\varphi\;d\pi_*\!\vol_{g^tN_+\Sigma'}\\
& =\int_{g^tN_+\Sigma'}\;\varphi\circ\pi\;d\vol_{g^tN_+\Sigma'}\;.
\end{align*}
En particulier, en prenant une suite exhaustive $(K_n)_{n\in\NN}$ de
compacts de $M$, et $\varphi:M\ra[0,1] $ une application continue
constante \'egale \`a $1$ sur $K_n$ et nulle sur $^c K_{n+1}$, et en
faisant tendre $n$ vers $+\infty$, la convergence \'etant assur\'ee
par l'hypoth\`ese de pincement de la courbure, la finitude du volume
de $\Sigma(t)$ et le théorème de convergence dominée de Lebesque, nous
avons
$$
\int_{\Sigma(t)} \;\kappa'_t \;d\vol_{\Sigma(t)} =\Vol(g^tN_+\Sigma')\;.
$$
Or par l'hypoth\`ese,
\begin{align*}
\lim_{t\ra+\infty}
\oint_{g^tN_+\Sigma'}\varphi\circ\pi\;d\vol_{g^tN_+\Sigma'} &=
\oint_{T^1M} \varphi\circ \pi \;d\vol_{T^1M}
= \frac{1}{\operatorname{Vol}(M)}\int_{M} \varphi
\;d\vol_{M}\;,
\end{align*}
car les mesures conditionnelles de $d\vol_{T^1M}$ sur les fibres de
$T^1M\ra M$ ont une masse totale constante. Ceci conclut la premi\`ere
\'etape de cette d\'emonstration.

\medskip
\noindent{\bf \'Etape 2. } Montrons maintenant que, pour tout
$\varphi\in\C_c(M)$, la $(\kappa\circ N_+)$-moyenne de $\varphi$ sur
$\Sigma(t)$ converge aussi vers la moyenne riemannienne de $\varphi$
sur $M$ quand $t\ra+\infty$.

\medskip Notons que $\kappa\circ N_+$ est mesurable born\'e, car $M$ est
à courbure pincée. Donc la $(\kappa\circ N_+)$-moyenne de $\Sigma(t)$
est bien définie. Quitte \`a d\'ecomposer $\varphi$ en partie positive
et partie n\'egative, nous pouvons supposer que $\varphi\geq
0$. Soient $\epsilon\in\;]0,1[$ et $K$ un compact de $M$, contenant le
support de $\varphi$, tels que $\Vol(\,^cK)\leq
\frac{\epsilon}{4}\Vol(M)$.

\medskip Soit $\widetilde{x}$ un point d'un rev\^etement riemannien
universel $\widetilde{M}$ de $M$. Pour tout $t>0$, soit
$\widetilde{C_t}$ une hypersurface lisse strictement convexe plong\'ee
dans $\widetilde{M}$ telle que $\widetilde{x}\in\partial
\V_t(\widetilde{C_t})$.  Notons $N(\widetilde{x},t)$ le vecteur
unitaire normal sortant de $\V_t(\widetilde{C_t})$ en
$\widetilde{x}$. Soit $H_t$ l'horosph\`ere de $\widetilde{M}$ passant
par $\widetilde{x}$ et de vecteur unitaire normal sortant en
$\widetilde{x}$ \'egal \`a $N(\widetilde{x},t)$. Alors $H_t$ et
$\partial \N_t(\widetilde{C_t})$ ont quand $t$ tend vers $+\infty$ une
tangence \`a l'ordre deux, uniform\'ement quand $\widetilde{x}$
parcourt un compact de $\widetilde{M}$.

Par ce qui pr\'ec\`ede, si $t$ est assez grand, alors pour tout $x$
dans $K\cap \Sigma(t)$, nous avons
\begin{equation}\label{eq:courburesproches}
|\kappa\circ N_+(x) - \kappa'_t(x)| \leq \epsilon\;.
\end{equation} 
En particulier,
\begin{align*}
\Big|\int  \varphi\;\kappa\circ N_+(x)\; d\vol_{\Sigma(t)}-
\int  \varphi\;\kappa'_t\; d\vol_{\Sigma(t)}\Big| & \leq 
\epsilon \int  \varphi\; d\vol_{\Sigma(t)}
\leq \epsilon \int \varphi\; \kappa'_t(x)\;d\vol_{\Sigma(t)}
\;,
\end{align*}
car $\kappa'_t\geq 1$, et donc
\begin{equation}\label{eq:passkpakcircNp}
(1-\epsilon)\int  \varphi\;\kappa'_t(x)\;d\vol_{\Sigma(t)}
\leq \int  \varphi\;\kappa\circ N_+\; d\vol_{\Sigma(t)}\leq 
(1+\epsilon)\int  \varphi\; \kappa'_t(x)\;d\vol_{\Sigma(t)}
\;.
\end{equation}

Puisque la vari\'et\'e $M$ n'est pas compacte, nous ne pouvons pas
prendre $\varphi=1$ dans la formule pr\'ec\'edente pour estimer
(uniform\'ement en $t$) l'int\'egrale $\int \;\kappa\circ N_+\;
d\vol_{\Sigma(t)}$, et nous avons besoin de montrer qu'il n'y a pas de
ph\'enom\`ene de perte de masse \`a l'infini.

Fixons $\varphi_K$ une fonction de $M$ dans $[0,1]$, continue \`a
support compact contenu dans $K$, telle que
$\int\varphi_K\;d\vol_M\geq\vol(K)-\frac{\epsilon}{4}\Vol(M)$. En
particulier,
$$
\oint\varphi_K\;d\vol_M\geq \frac{\Vol(K)}{\Vol(M)}-
\frac{\epsilon}{4}\geq 1-\frac{\epsilon}{2}\;.
$$
Par l'\'etape 1 appliqu\'ee \`a $\varphi_K$, pour tout $t$ assez grand, 
nous avons
$$
\frac{\int\;\varphi_K\;\kappa'_t\;d\vol_{\Sigma(t)}}
{\int\;\kappa'_t\;d\vol_{\Sigma(t)}}\geq 
\oint\varphi_K\;d\vol_M\; -\frac{\epsilon}{2}\geq 1-\epsilon\;.
$$
Donc
\begin{equation}\label{eq:papertmass}
\frac{\int_{\,^cK}\;\kappa'_t\;d\vol_{\Sigma(t)}}
{\int\;\kappa'_t\;d\vol_{\Sigma(t)}}=1-
\frac{\int_{K}\;\kappa'_t\;d\vol_{\Sigma(t)}}
{\int\;\kappa'_t\;d\vol_{\Sigma(t)}}\leq 1-
\frac{\int\;\varphi_K\;\kappa'_t\;d\vol_{\Sigma(t)}}
{\int\;\kappa'_t\;d\vol_{\Sigma(t)}}\leq \epsilon\;.
\end{equation}
Puisque $K$ et $^cK$ sont disjoints de r\'eunion $M$, nous avons donc,
par les inégalités \eqref{eq:papertmass} et
\eqref{eq:courburesproches}, et puisque $\kappa\geq 1$,
\begin{align*}
(1-\epsilon)\;\int\;\kappa'_t\;d\vol_{\Sigma(t)}&\leq
\int_{K}\;\kappa'_t\;d\vol_{\Sigma(t)}\leq \int_{K}\;\kappa\circ
N_+\;d\vol_{\Sigma(t)}+\;\epsilon\;\int_K d\vol_{\Sigma(t)}\\ &\leq
(1+\epsilon)\;\int_{K}\;\kappa\circ N_+\;d\vol_{\Sigma(t)} \leq
(1+\epsilon)\;\int\;\kappa\circ N_+\;d\vol_{\Sigma(t)}\;,
\end{align*}
c'est-\`a-dire
\begin{equation}\label{eq:passkcircNakpinf}
\frac{1-\epsilon}{1+\epsilon}\;\int\;\kappa'_t\;d\vol_{\Sigma(t)}
\leq \int\;\kappa\circ N_+\;d\vol_{\Sigma(t)}\;.
\end{equation}

Montrons une comparaison analogue dans l'autre sens.
Puisque $\kappa'_t\geq 1$,  par l'in\'egalit\'e 
\eqref{eq:papertmass}, nous avons
\begin{equation}\label{eq:papertmassdeux}
\frac{\int_{\,^cK}\;d\vol_{\Sigma(t)}}
{\int\;\kappa'_t\;d\vol_{\Sigma(t)}}\leq
\frac{\int_{\,^cK}\;\kappa'_t\;d\vol_{\Sigma(t)}}
{\int\;\kappa'_t\;d\vol_{\Sigma(t)}}\leq \epsilon\;.
\end{equation}
Puisque la courbure sectionnelle de $M$ est pinc\'ee, et par comparaison
avec le cas \`a courbure constante, o\`u les courbures principales des
horosph\`eres sont des constantes strictement positives, les courbures
principales des horosph\`eres de $\wt M$ sont born\'ees. Il existe donc
$a\in[1,+\infty[$ tel que $1\leq \kappa\leq a$.

D'o\`u, respectivement par la positivit\'e de $\kappa'_t$, par
l'encadrement \eqref{eq:courburesproches}, puisque $\kappa\geq 1$,
puisque $M=K\sqcup \,^cK$, puisque $\kappa\leq a$, et par
l'in\'egalit\'e \eqref{eq:papertmassdeux}, nous avons
\begin{align*}
\int\;\kappa'_t\;d\vol_{\Sigma(t)}&
\geq \int_K\;\kappa'_t\;d\vol_{\Sigma(t)}
\geq \int_K\;\kappa\circ N_+\;d\vol_{\Sigma(t)}-
\;\epsilon\;\int_K\;\;d\vol_{\Sigma(t)}\\
&\geq (1-\epsilon)\int_K\;\kappa\circ N_+\;d\vol_{\Sigma(t)}\\
& = (1-\epsilon)\int\;\kappa\circ N_+\;d\vol_{\Sigma(t)}-
(1-\epsilon)\int_{\,^cK}\;\kappa\circ N_+\;d\vol_{\Sigma(t)}\\
& \geq(1-\epsilon)\int\;\kappa\circ N_+\;d\vol_{\Sigma(t)}-
(1-\epsilon)a\int_{\,^cK}\;d\vol_{\Sigma(t)}\\
&\geq (1-\epsilon)\int\;\kappa\circ N_+\;d\vol_{\Sigma(t)}-
(1-\epsilon)a\epsilon \int\;\kappa'_t\;d\vol_{\Sigma(t)}\;.
\end{align*}
Par cons\'equent,
\begin{equation}\label{eq:passkcircNakpsup}
\int\;\kappa\circ N_+\;d\vol_{\Sigma(t)}\leq
\frac{1+(1-\epsilon)a\epsilon}{1-\epsilon}\;\int\;\kappa'_t\;
d\vol_{\Sigma(t)}\;.
\end{equation}
La seconde \'etape d\'ecoule alors, outre de l'\'etape 1, de
l'encadrement \eqref{eq:passkpakcircNp} et des in\'egalit\'es
\eqref{eq:passkcircNakpinf} et \eqref{eq:passkcircNakpsup}. \cqfd

\medskip Remarquons que dans ce r\'esultat
\ref{theo:equidistributionenbas}, la r\'egularit\'e de $\Sigma'$ peut
sans doute \^etre affaiblie en ${\rm C}^1$. En effet, celle-ci permet
de d\'efinir l'application de Gauss (sortante) $N_+$ continue et la
moyenne riemannienne de $\Sigma'$, et si $t>0$, l'hypersurface
$\Sigma(t)=\pi(g^t N_+\Sigma')$ est ${\rm C}^{1,1}$ (voir par exemple
\cite{Walter76}), donc ses courbures principales sont d\'efinies
presque partout (pour la classe de mesure de Lebesgue). De plus, $g^t
N_+\Sigma'$ est une sous-vari\'et\'e lipschitzienne, et sa moyenne
riemannienne est encore d\'efinie.

\section{Comptage d'arcs perpendiculaires communs}
\label{sect:comptage}

Le but de cette partie est de donner un \'equivalent pr\'ecis, quand
$t$ tend vers $+\infty$, du nombre de rayons (localement)
g\'eod\'esiques partant perpendiculairement d'une sous-vari\'et\'e
totalement g\'eod\'esique et convergeant vers une pointe d'une vari\'et\'e
hyperbolique de volume fini, de longueur (naturellement normalis\'ee)
au plus $t$.

\medskip Soit $M$ une vari\'et\'e riemannienne lisse connexe,
compl\`ete, de volume fini, de dimension $n$ au moins $2$ et de
courbure sectionnelle au plus $-1$. Nous notons $\pi:T^1M\ra M$ son
fibr\'e unitaire tangent et $(g^t)_{t\in\RR}$ le flot g\'eod\'esique
sur $T^1M$. Nous munissons $T^1M$ de la m\'etrique riemannienne de
Sasaki.

Soit $C$ une sous-vari\'et\'e immerg\'ee lisse de volume fini dans
$M$, totalement g\'eod\'esique, de dimension $k\in [0,n-1]$.  Notons
que $C$ n'est pas suppos\'ee connexe. Par exemple, $C$ peut \^etre une
union finie de points ($k=0$), ou une union finie de g\'eod\'esiques
ferm\'ees ($k=1$), pas forc\'ement simples, et \'eventuellement
parcourues plusieurs fois.

\medskip Une mani\`ere naturelle de normaliser la longueur d'un rayon
g\'eod\'esique sortant de tout compact de $M$ est de ne compter sa
longueur que jusqu'au dernier instant o\`u il entre dans un petit
voisinage prescrit d'un bout de $M$. Il existe une mani\`ere canonique
de d\'efinir de tels voisinages (voir par exemple
\cite{BusKar91,BalGroSch85}). Nous aurons besoin de leur d\'efinition
pour n'importe quelle vari\'et\'e $M'$ riemannienne lisse connexe,
compl\`ete, de dimension $n$ au moins $2$ et de courbure sectionnelle
au plus $-1$, qui est g\'eom\'etriquement finie (voir par exemple
\cite{Bowditch95}). Rappelons qu'une {\it pointe} ({\it cusp} en
anglais) de $M'$ est une classe de rayons g\'eod\'esiques minimisants
asymptotes, le long desquels le rayon d'injectivit\'e tend vers
$0$. Lorsque $M'$ est de volume fini, l'ensemble des pointes de $M'$
est en bijection avec l'ensemble des bouts de l'espace topologique
localement compact $M'$, par l'application qui \`a un repr\'esentant
d'une pointe associe le bout vers lequel il converge. Soit $\wt\pi:\wt
M'\ra M'$ un rev\^etement universel riemannien de $M$, de groupe de
rev\^etement $\Ga$. Si $e$ est une pointe de $M$, si $\rho_e$ est un
rayon g\'eod\'esique minimisant dans la classe $e$, puisque $M'$ est
g\'eom\'etriquement finie, il existe (voir par exemple
\cite{Bowditch95}) une horoboule ouverte $H_e$ dans $\wt M'$,
centr\'ee au point \`a l'infini $\xi_e$ d'un relev\'e de $\rho_e$ dans
$\wt M'$, telle que $H_e$ et $\ga H_e$ soient disjointes si
$\ga\in\Ga$ ne fixe pas $\xi_e$.  En particulier, le quotient de $H_e$
par son stabilisateur dans $\Ga$ s'injecte dans $M$ par l'application
induite par $\wt \pi_{\mid H_e}$. L'image $V_e$ de $H_e$ dans $M'$
par $\wt \pi$ est appel\'ee un {\it voisinage cuspidal} (ou {\it de
  Margulis}) de la pointe $e$.

\medskip Notons $\H$ une union finie non vide de voisinages cuspidaux
de $M$, et $\partial\H$ la r\'eunion correspondante des bords de ces
voisinages cuspidaux. Munissons chacun de ces bord de la m\'etrique
riemannienne induite. Nous ne supposons pas que ces bords sont deux
\`a deux disjoints, ni qu'ils soient disjoints de $C$, ni que toutes
les pointes de $M$ interviennent. Nous noterons $\Vol(\H)$ la somme
des volumes riemanniens de ces voisinages cuspidaux et
$\Vol(\partial\H)$ celle de leurs bords. Le r\'esultat suivant est un
calcul \'el\'ementaire bien connu.

\blemm\label{lem:volumecuspvsbord}
Si $M$ est hyperbolique, alors $\Vol(\H)=\frac{\Vol(\partial\H)}{n-1}$.
\elemm

\dem Il suffit, par additivit\'e, de montrer le r\'esultat lorque $\H$
ne comporte qu'un seul voisinage cuspidal. Notons
$$
\HH^n_\RR=\big(\{x=(x_1,\dots,x_n)\in\RR^n\;:\;x_n>0\},\;
\frac{dx_1^2+\dots+dx_n^2}{x_n^2}\;\big)
$$
le mod\`ele du demi-espace sup\'erieur de l'espace hyperbolique r\'eel
de dimension $n$. Choisissons un rev\^etement riemannien universel
$\wt \pi:\HH^n_\RR\ra M$ tel que l'horosph\`ere d'\'equation $x_n=1$,
que nous noterons $\wt{\partial \H}$, s'envoie sur $\partial
\H$. Notons $\Omega$ une partie mesurable de $\wt{\partial \H}$, telle
que $\wt \pi_{\mid \Omega}:\Omega\ra \partial \H$ soit une bijection.
Remarquons que l'\'el\'ement de volume hyperbolique s'\'ecrit
$d\vol_{\HH^n_\RR}=d\vol_{\wt{\partial \H}}\frac{dx_n}{x_n^n}$. Par
d\'efinition d'un voisinage cuspidal, nous avons donc
$$
\Vol(\H)=\int_{\Omega\times[1,+\infty[}d\vol_{\HH^n_\RR}= 
\Vol(\partial\H)\int_1^{+\infty}\frac{dt}{t^n}=
\frac{\Vol(\partial\H)}{n-1}\;.\;\;\;\mbox{\cqfd}
$$

\bigskip Pour tout $t\geq 0$, notons $\N(t)=\N_{M,C,\H}(t)$ le nombre
compt\'e avec multiplicit\'e (voir ci-dessous) de segments
(localement) g\'eod\'esiques $\delta$ entre un point de $C$ et un
point de $\partial\H$, partant perpendiculairement de $C$ (si $k>0$)
et arrivant perpendiculairement en $\partial\H$, de longueur
alg\'ebrique au plus $t$. La longueur de $\delta$ est compt\'ee
positivement s'il sort de $\H$ en son point de d\'epart sur $\partial
\H$, et n\'egativement sinon (cette convention n'a pas d'influence
pour notre \'equivalent donn\'e dans le th\'eor\`eme
\ref{theo:comptage}, car il n'y a qu'un nombre fini de tels $\delta$
de longueur plus petite qu'une valeur donn\'ee). En particulier, si
$C$ et $\H$ sont disjoints, alors les longueurs alg\'ebriques des
$\delta$ sont toutes strictement positives. Un point de tangence entre
$\partial\H$ et $C$ est un tel segment $\delta$ (de longueur nulle). Si
$\iota:C\ra M$ est l'immersion consid\'er\'ee de $C$ dans $M$, si $v$
est le vecteur tangent de $\delta$ en son origine notée $x$, alors
la {\it multiplicit\'e} de $\delta$ est
\begin{equation}\label{eq:defmulpvariete}
m(\delta)=\operatorname{Card}\{x'\in C\;:\;
\iota(x')=x\;\;{\rm et}\;\; T_{x'}\,\iota\,(T_{x'}C)\perp v\}\;.
\end{equation}

Dans la formule ci-dessous, $\SSS_{m}$ d\'esigne la sph\`ere unit\'e de
$\RR^{m+1}$, munie de sa m\'etrique riemannienne usuelle. 

\btheo \label{theo:comptage} Soit $M$ une vari\'et\'e hyperbolique
connexe, de volume fini, de dimension $n\geq 2$.  Soit $C$ une
sous-vari\'et\'e lisse immerg\'ee de volume fini dans $M$, totalement
g\'eod\'esique, de dimension $k<n$. Soit $\H$ une union finie non vide
de voisinages cuspidaux de $M$. Alors quand $t$ tend vers $+\infty$,
$$
\N_{M,C,\H}(t)\sim 
\frac{\Vol(\SSS_{n-k-1})\Vol(\H)\Vol(C)}{\Vol(\SSS_{n-1})\Vol(M)}
\;e^{(n-1)t}\;.
$$
\etheo

Nous renvoyons \`a la fin de cette partie pour des remarques
sur cet \'enonc\'e. Avant de commencer la preuve du th\'eor\`eme
\ref{theo:comptage}, nous d\'emontrons trois r\'esultats d'int\'er\^et
ind\'ependant en g\'eom\'etrie hyperbolique.

La proposition suivante est un exercice de g\'eom\'etrie hyperbolique
(voir \cite{Buser92} pour $n=2$ et \cite{Basmajian94} pour $k=n-1$).

\bprop \label{prop:computvol} Soient $\wh M$ une vari\'et\'e
hyperbolique de dimension $n\geq 2$, $\wh C$ une sous-vari\'et\'e
(plong\'ee) de volume fini totalement g\'eod\'esique de dimension $k<n$,
telle que l'inclusion de $\wh C$ dans $\wh M$ soit une \'equivalence
d'homotopie.  Quand $t$ tend vers $+\infty$, nous avons l'\'equivalent
$$
\Vol(\V_t\wh C)\sim 
\frac{\Vol(\SSS_{n-k-1})\Vol(\wh C)}{(n-1)2^{n-1}}\;\;e^{(n-1)t}\;.
$$
\eprop

Pour usage ult\'erieur, remarquons qu'il d\'ecoule de ce r\'esultat que,
pour tout $\epsilon\in\RR$,
\begin{equation}\label{eq:equivvolplusepsilon}
\Vol(\V_{t+\epsilon}\wh C)\sim e^{(n-1)\epsilon}\Vol(\V_t\wh C)
\end{equation}
Pour une valeur exacte du volume du $t$-voisinage de $\wh C$, voir la 
formule \eqref{eq:formulevolum}, qui lorsque $k=1$ devient
$$
\Vol(\V_t\wh C)=
\frac{\Vol(\SSS_{n-2})\Vol(\wh C)}{(n-1)}\;(\sinh t)^{n-1}\;.
$$

\dem Par additivité, nous pouvons supposer que $\wh C$ (et donc $\wh
M$) est connexe. Consid\'erons de nouveau $\HH^n_\RR$, le mod\`ele du
demi-espace sup\'erieur de l'espace hyperbolique r\'eel de dimension
$n$. Si $k=0$, alors $\wh C$ est un point, et le résultat est
connu. Donc nous supposons que $k\geq 1$. Notons $\HH^k_\RR$ son
sous-espace hyperbolique d'\'equations $x_1=0,\dots,
x_{n-k}=0$. Choisissons un rev\^etement riemannien universel $\wt
\pi:\HH^n_\RR\ra \wh M$ tel que $\wt \pi(\HH^k_\RR)=\wh C$.  Notons
$\Omega$ une partie mesurable de $\HH^k_\RR$, telle que $\wt \pi_{\mid
  \Omega}:\Omega\ra \wh C$ soit une bijection. En particulier,

$$
\int_\Omega d\vol_{\HH^k_\RR}=\Vol(\wh C)\;.
$$
Pour tout $x=(x_1,\dots,x_n)\in\HH^n_\RR$, notons $y$ la derni\`ere
coordonn\'ee de la projection orthogonale hyperbolique de $x$ sur
$\HH^k_\RR$, notons $s$ la distance hyperbolique de $x$ \`a $\HH^k_\RR$,
notons $r=(\sum_{i=1}^{n-k}x_i^2)^{1/2}$, et notons $\theta$ l'angle
(euclidien) de l'arc de cercle entre la projection de $x$ et $x$ port\'e
par la g\'eod\'esique hyperbolique entre ces deux points, s'ils sont
distincts (voir la figure ci-dessous).

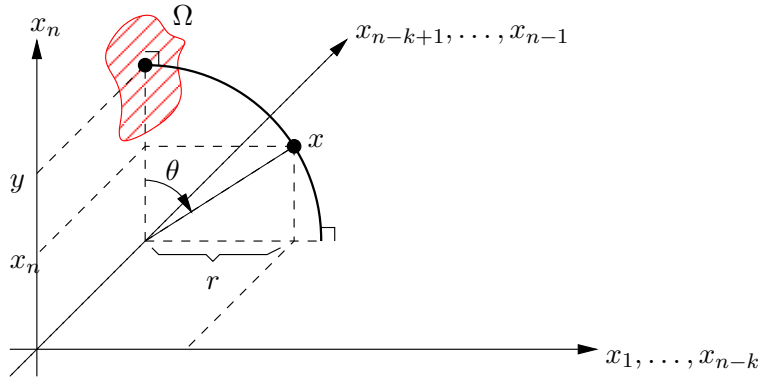
\begin{center}
\begin{picture}(0,0)%
\includegraphics{fig_calcvoltvois.pstex}%
\end{picture}%
\setlength{\unitlength}{3729sp}%
\begingroup\makeatletter\ifx\SetFigFontNFSS\undefined%
\gdef\SetFigFontNFSS#1#2#3#4#5{%
  \reset@font\fontsize{#1}{#2pt}%
  \fontfamily{#3}\fontseries{#4}\fontshape{#5}%
  \selectfont}%
\fi\endgroup%
\begin{picture}(3990,2535)(886,-2053)
\put(2206,-1456){\makebox(0,0)[lb]{\smash{{\SetFigFontNFSS{11}{13.2}{\rmdefault}{\mddefault}{\updefault}{\color[rgb]{0,0,0}$r$}%
}}}}
\put(1936,-736){\makebox(0,0)[lb]{\smash{{\SetFigFontNFSS{11}{13.2}{\rmdefault}{\mddefault}{\updefault}{\color[rgb]{0,0,0}$\theta$}%
}}}}
\put(1981,299){\makebox(0,0)[lb]{\smash{{\SetFigFontNFSS{11}{13.2}{\rmdefault}{\mddefault}{\updefault}{\color[rgb]{0,0,0}$\Omega$}%
}}}}
\put(1036,254){\makebox(0,0)[lb]{\smash{{\SetFigFontNFSS{11}{13.2}{\rmdefault}{\mddefault}{\updefault}{\color[rgb]{0,0,0}$x_n$}%
}}}}
\put(3196,209){\makebox(0,0)[lb]{\smash{{\SetFigFontNFSS{11}{13.2}{\rmdefault}{\mddefault}{\updefault}{\color[rgb]{0,0,0}$x_{n-k+1},\dots, x_{n-1}$}%
}}}}
\put(4861,-1951){\makebox(0,0)[lb]{\smash{{\SetFigFontNFSS{11}{13.2}{\rmdefault}{\mddefault}{\updefault}{\color[rgb]{0,0,0}$x_1,\dots, x_{n-k}$}%
}}}}
\put(2881,-511){\makebox(0,0)[lb]{\smash{{\SetFigFontNFSS{11}{13.2}{\rmdefault}{\mddefault}{\updefault}{\color[rgb]{0,0,0}$x$}%
}}}}
\put(901,-781){\makebox(0,0)[lb]{\smash{{\SetFigFontNFSS{11}{13.2}{\rmdefault}{\mddefault}{\updefault}{\color[rgb]{0,0,0}$y$}%
}}}}
\put(901,-1321){\makebox(0,0)[lb]{\smash{{\SetFigFontNFSS{11}{13.2}{\rmdefault}{\mddefault}{\updefault}{\color[rgb]{0,0,0}$x_n$}%
}}}}
\end{picture}%

\end{center}

Une formule classique (voir par exemple \cite[page 145]{Beardon83})
montre que $\sinh s=\tan\theta$. Nous avons $x_n=y\cos\theta$ et
$r=y\sin\theta$, de sorte que $dr\wedge dx_n= y\, d\theta\wedge dy$.
La mesure de Lebesgue de $\RR^{n-k}$ s'\'ecrit en coordonn\'ees
sph\'eriques $dx_1\dots dx_{n-k}=r^{n-k-1}dr\;
d\vol_{\SSS_{n-k-1}}$. L'\'el\'ement de volume hyperbolique s'\'ecrit
donc
\begin{align*}
d\vol_{\HH^n_\RR}&=\frac{dx_1\dots dx_n}{x_n^n}=
\frac{r^{n-k-1}}{x_n^n}\;dr\;d\vol_{\SSS_{n-k-1}}dx_{n-k+1}\dots dx_n
\\ & =\frac{\sin^{n-k-1}\theta}{y^k\cos^n\theta}\;
d\theta \;d\vol_{\SSS_{n-k-1}}dx_{n-k+1}\dots dx_{n-1}dy\;.
\end{align*}

Pour tout $t>0$, posons $\theta_t=\arctan(\sinh t)$. Puisque $\wh M$
se r\'etracte par d\'eformation sur $\wh C$, l'application $\wt\pi$
envoie le domaine $\D=\{x\in\HH^n_\RR\;:\; y\in\Omega, \theta\in
[0,\theta_t]\}$ sur $\V_t\wh C$ en pr\'eservant la mesure. Donc
\begin{align}
\Vol(\V_t\wh C)&=\int_\D d\vol_{\HH^n_\RR}=
\int_{\SSS_{n-k-1}} d\vol_{\SSS_{n-k-1}} \int_0^{\theta_t} 
\frac{\sin^{n-k-1}\theta}{\cos^n\theta} \;d\theta
\int_\Omega \frac{dx_{n-k+1}\dots dx_{n-1}dy}{y^k}\notag
\\ & =\Vol(\SSS_{n-k-1})\Vol(\wh C)
\int_0^{t}\tanh^{n-k-1}s\;\cosh^{n-1} s \;ds\;,
\label{eq:formulevolum}
\end{align}
en faisant le changement de variable $\theta=\arctan(\sinh s)$. La
remarque suivant l'\'enonc\'e de la proposition \ref{prop:computvol}
en d\'ecoule. Un argument classique de d\'ecoupage en deux de
l'int\'egrale et les \'equivalents $\tanh s\sim 1$ et $\cosh s\sim
e^s/2$ quand $s$ tend vers $+\infty$ montrent facilement la
proposition \ref{prop:computvol}. \cqfd

\medskip Le second r\'esultat de g\'eom\'etrie hyperbolique est sans
doute aussi bien connu, et une preuve nous a \'et\'e communiqu\'ee par
P.~Pansu.  Si $X=\HH^n_\RR$ est l'espace hyperbolique r\'eel de
dimension $n\geq 2$, nous appellerons {\it horosph\`ere instable}
(resp.~{\it stable}) le relev\'e, dans le fibr\'e unitaire tangent
$T^1X$, par la normale sortante (resp.~rentrante) d'une horosph\`ere
de $X$.

\bprop \label{prop:pansujacobi} Dans $T^1\HH^n_\RR$, une horosph\`ere
instable et une horosph\`ere stable, qui se rencontrent, sont
orthogonales en leur point d'intersection pour la m\'etrique de
Sasaki.  \eprop

\dem La preuve utilise l'isomorphisme lin\'eaire bien connu, pour tout
vecteur tangent unitaire $v$ d'une vari\'et\'e riemannienne compl\`ete
$M$, entre l'espace vectoriel $T_vTM$ et celui des champs de Jacobi le
long de la g\'eod\'esique d\'efinie par $v$.

Notons $M=\HH^n_\RR$ et $p:TM\ra M$ la projection canonique dans
cette preuve. Rappelons (voir le d\'ebut de la d\'emonstration de la
proposition \ref{prop:pansu}) que la m\'etrique de Sasaki sur $TM$ est
d\'efinie en posant, pour tous $v\in T^1M$, et $X,Y\in T_vTM$
\begin{equation}\label{eq:formuleSasaki}
\langle X,Y\rangle_{TM}=\langle \nabla_v (Z_X),\nabla_v (Z_Y)\rangle_M+
\langle T_vp(X),T_vp(Y)\rangle_M\;,
\end{equation}
o\`u $Z_X:M\ra TM$ est n'importe quel champ de vecteurs lisse tels que
$Z_X(p(v))=v$ et $TZ_X(v)=X$.

Pour tous $v\in T^1M$ et $X\in T_vTM$, notons $t\mapsto J_X(t)$ le
champ de Jacobi le long de la g\'eod\'esique $\ga_v:t\mapsto p(g^tv)$
tel que $J_X(0)=T_vp(X)$ et $\dot{J}_X(0)=\nabla_v Z_X$. Par transport
parall\`ele le long de $\ga_v$, nous consid\'erons que $J_X(t)\in
T_{p(v)}M$ pour tout $t\in \RR$. Remarquons que l'espace tangent \`a
l'horosph\`ere stable (resp.~instable) passant par $v$ est l'ensemble
des $X\in T_vTM$ tels que $J_X$ soit orthogonal le long de $\ga_v$ et
$\lim_{t\ra+\infty}J_X(t)=0$ (resp.~$\lim_{t\ra-\infty}J_X(t)=0$).

Puisque $M$ est \`a courbure constante $-1$, l'\'equation des champs
de Jacobi est $\ddot{J}(t)-J(t)=0$. Soit $X\in T_vTM$ (resp.~$Y\in
T_vTM$) un vecteur tangent \`a l'horosph\`ere stable (resp.~instable)
passant par $v$. Alors $J_X(t)=e^{-t}J_X(0)$, donc
$\dot{J}_X(0)=-J_X(0)$. De m\^eme, nous avons $\dot{J}_Y(0)=J_Y(0)$.
Donc, par la formule \eqref{eq:formuleSasaki}, nous avons
$\langle X,Y\rangle_{TM}=\langle \dot{J}_X(0),\dot{J}_Y(0)\rangle+ 
\langle {J}_X(0),{J}_Y(0)\rangle=0\;,$
ce qui montre le r\'esultat.
\cqfd

\medskip Le troisi\`eme r\'esultat est en fait valable plus
g\'en\'eralement que pour l'espace hyperbolique r\'eel $X=\HH^n_\RR$.

\blemm\label{lem:procherentranthoro} Pour tout $\epsilon>0$, il existe
$\delta>0$ tel que pour toute hypersurface ${\rm C}^1$ strictement
convexe $S$ d'une vari\'et\'e riemannienne $X$ compl\`ete, simplement
connexe, de dimension $n\geq 2$, \`a courbure sectionnelle pinc\'ee au
plus $-1$, pour toute horoboule ferm\'ee $H$ disjointe de $S$, pour tout
vecteur $v\in N_+S$, pour tout $s\geq 1$, si $g^sv\in
\V_\delta\big(N_-(\partial H)\big)$, alors il existe $v'\in N_+S$ et
$s'\geq 0$ tel que $g^{s'}v'\in N_-(\partial H)\cap
B_{T^1X}(g^sv,\epsilon)$.  \elemm

\dem Soit $(X,S,H,v,s)$ comme dans l'\'enonc\'e. Notons $s'=d(S,H)$
(qui est strictement positive), et $v'$ le vecteur tangent \`a
l'origine du segment perpendiculaire commun allant de $S$ \`a
$H$. Nous avons en particulier $g^{s'}v'\in N_-(\partial H)$. Nous
pouvons supposer que $v'\neq v$. Pour $\pi:T^1X\ra X$ la projection
canonique, d\'efinissons $x=\pi(v)$, $y=\pi(g^sv)$, $x'=\pi(v')$,
$y'=\pi(g^{s'}v')$ (voir la figure ci-dessous \`a gauche). Montrons
que $g^{s'}v'$ est proche de $g^sv$ si $g^sv$ est suffisamment proche
d'un \'el\'ement de $N_-(\partial H)$ (avec constantes uniformes).

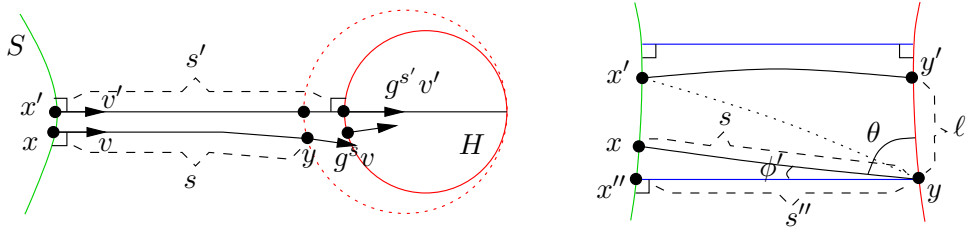
\begin{center}
\begin{picture}(0,0)%
\includegraphics{fig_quadrilatere.pstex}%
\end{picture}%
\setlength{\unitlength}{3729sp}%
\begingroup\makeatletter\ifx\SetFigFontNFSS\undefined%
\gdef\SetFigFontNFSS#1#2#3#4#5{%
  \reset@font\fontsize{#1}{#2pt}%
  \fontfamily{#3}\fontseries{#4}\fontshape{#5}%
  \selectfont}%
\fi\endgroup%
\begin{picture}(6328,1570)(256,-1339)
\put(271,-130){\makebox(0,0)[lb]{\smash{{\SetFigFontNFSS{11}{13.2}{\rmdefault}{\mddefault}{\updefault}{\color[rgb]{0,0,0}$S$}%
}}}}
\put(379,-760){\makebox(0,0)[lb]{\smash{{\SetFigFontNFSS{11}{13.2}{\rmdefault}{\mddefault}{\updefault}{\color[rgb]{0,0,0}$x$}%
}}}}
\put(883,-769){\makebox(0,0)[lb]{\smash{{\SetFigFontNFSS{11}{13.2}{\rmdefault}{\mddefault}{\updefault}{\color[rgb]{0,0,0}$v$}%
}}}}
\put(1450,-1015){\makebox(0,0)[lb]{\smash{{\SetFigFontNFSS{11}{13.2}{\rmdefault}{\mddefault}{\updefault}{\color[rgb]{0,0,0}$s$}%
}}}}
\put(2449,-874){\makebox(0,0)[lb]{\smash{{\SetFigFontNFSS{11}{13.2}{\rmdefault}{\mddefault}{\updefault}{\color[rgb]{0,0,0}$g^sv$}%
}}}}
\put(2788,-388){\makebox(0,0)[lb]{\smash{{\SetFigFontNFSS{11}{13.2}{\rmdefault}{\mddefault}{\updefault}{\color[rgb]{0,0,0}$g^{s'}v'$}%
}}}}
\put(3280,-790){\makebox(0,0)[lb]{\smash{{\SetFigFontNFSS{11}{13.2}{\rmdefault}{\mddefault}{\updefault}{\color[rgb]{0,0,0}$H$}%
}}}}
\put(364,-529){\makebox(0,0)[lb]{\smash{{\SetFigFontNFSS{11}{13.2}{\rmdefault}{\mddefault}{\updefault}{\color[rgb]{0,0,0}$x'$}%
}}}}
\put(1465,-211){\makebox(0,0)[lb]{\smash{{\SetFigFontNFSS{11}{13.2}{\rmdefault}{\mddefault}{\updefault}{\color[rgb]{0,0,0}$s'$}%
}}}}
\put(2236,-844){\makebox(0,0)[lb]{\smash{{\SetFigFontNFSS{11}{13.2}{\rmdefault}{\mddefault}{\updefault}{\color[rgb]{0,0,0}$y$}%
}}}}
\put(904,-484){\makebox(0,0)[lb]{\smash{{\SetFigFontNFSS{11}{13.2}{\rmdefault}{\mddefault}{\updefault}{\color[rgb]{0,0,0}$v'$}%
}}}}
\put(5279,-936){\makebox(0,0)[lb]{\smash{{\SetFigFontNFSS{11}{13.2}{\rmdefault}{\mddefault}{\updefault}{\color[rgb]{0,0,0}$\phi'$}%
}}}}
\put(5990,-714){\makebox(0,0)[lb]{\smash{{\SetFigFontNFSS{11}{13.2}{\rmdefault}{\mddefault}{\updefault}{\color[rgb]{0,0,0}$\theta$}%
}}}}
\put(4263,-316){\makebox(0,0)[lb]{\smash{{\SetFigFontNFSS{11}{13.2}{\rmdefault}{\mddefault}{\updefault}{\color[rgb]{0,0,0}$x'$}%
}}}}
\put(4258,-770){\makebox(0,0)[lb]{\smash{{\SetFigFontNFSS{11}{13.2}{\rmdefault}{\mddefault}{\updefault}{\color[rgb]{0,0,0}$x$}%
}}}}
\put(4195,-1082){\makebox(0,0)[lb]{\smash{{\SetFigFontNFSS{11}{13.2}{\rmdefault}{\mddefault}{\updefault}{\color[rgb]{0,0,0}$x''$}%
}}}}
\put(6393,-1078){\makebox(0,0)[lb]{\smash{{\SetFigFontNFSS{11}{13.2}{\rmdefault}{\mddefault}{\updefault}{\color[rgb]{0,0,0}$y$}%
}}}}
\put(6569,-672){\makebox(0,0)[lb]{\smash{{\SetFigFontNFSS{11}{13.2}{\rmdefault}{\mddefault}{\updefault}{\color[rgb]{0,0,0}$\ell$}%
}}}}
\put(6348,-230){\makebox(0,0)[lb]{\smash{{\SetFigFontNFSS{11}{13.2}{\rmdefault}{\mddefault}{\updefault}{\color[rgb]{0,0,0}$y'$}%
}}}}
\put(5005,-595){\makebox(0,0)[lb]{\smash{{\SetFigFontNFSS{11}{13.2}{\rmdefault}{\mddefault}{\updefault}{\color[rgb]{0,0,0}$s$}%
}}}}
\put(5449,-1270){\makebox(0,0)[lb]{\smash{{\SetFigFontNFSS{11}{13.2}{\rmdefault}{\mddefault}{\updefault}{\color[rgb]{0,0,0}$s''$}%
}}}}
\end{picture}%

\end{center}

Soit $\epsilon'>0$. La projection canonique $\pi$ est
$1$-lipschitzienne, et isom\'etrique sur chaque orbite du flot
g\'eod\'esique. Donc quitte \`a remplacer $H$ par une horoboule
ferm\'ee concentrique et $s'$ par $s'\pm d(y,\partial H)$, nous
pouvons supposer que $y\in \partial H$, et que l'angle en $y$ entre
$g^sv$ et le vecteur unitaire normal rentrant dans $H$ en $y$ est au
plus $\epsilon'$.

Consid\'erons le quadrilat\`ere g\'eod\'esique $Q$ de sommets
$x,x',y',y$. Par convexit\'e, ses angles en $x,x',y'$ sont au moins
$\frac{\pi}{2}$, et l'angle en $y$, not\'e $\theta$, appartient \`a
$[\frac{\pi}{2}-\epsilon',\frac{\pi}{2}[$.

Montrons que $\ell=d(y,y')$ tend uniform\'ement vers $0$ quand
$\epsilon'$ tend vers $0$. Ceci implique le r\'esultat, car puisque
les horosph\`eres ont des courbures principales born\'ees par
l'hypoth\`ese de pincement de la courbure de $X$, les vecteurs
unitaires normaux rentrant dans $H$ en $y$ et $y'$ seront
uniform\'ement proches.

Par comparaison, nous pouvons supposer que la courbure est constante,
\'egale \`a $-1$. Puisque remplacer le triangle g\'eod\'esique de
sommets $x',y,y'$ par le triangle g\'eod\'esique isom\'etrique, de
sommets $x',y,y''$, contenu dans le plan hyperbolique contenant
$x,y,x'$, et de l'autre c\^ot\'e de la g\'eod\'esique passant par
$x',y$, ne fait qu'augmenter les angles en $x'$ et $y$ sans changer
les autres angles ni les longueurs des c\^ot\'es, nous pouvons
supposer que $Q$ est contenu dans un plan hyperbolique. Quitte \`a
remplacer le c\^ot\'e $[x',y']$ par le segment g\'eod\'esique
perpendiculaire commun aux droites passant par $x,x'$ et $y,y'$, ce
qui ne peut qu'augmenter $\ell=d(y,y')$ et laisse $s=d(x,y)$
inchang\'e, nous pouvons supposer que les angles de $Q$ en $x'$ et
$y'$ sont \'egaux \`a $\frac{\pi}{2}$.

Notons $x''$ la projection orthogonale de $y$ sur la g\'eod\'esique
passant par $x,x'$, notons $\phi'\in[0,\epsilon']$ l'angle en $y$ du
triangle g\'eodesique de sommets $x'',y,x$, et posons $s''=d(x'',y)$
(voir la figure ci-dessus \`a droite).  Par la formule 7.11.2(iii) de
\cite[page 147]{Beardon83}, nous avons
$$
\tanh s''=\tanh s\;\cos\phi'\geq \tanh 1\;\cos\epsilon'\;,
$$
car $s\geq 1$. Donc en utilisant les formules de \cite[page
157]{Beardon83} pour le quadrilatère de sommets $x'',x',y',y$ dont
l'angle en $y$ est au moins $\theta$ et les autres angles sont droits,
nous avons
$$
\tanh \ell\leq\frac{\cos\theta}{\tanh s''}\leq
\frac{\cos(\frac{\pi}{2}-\epsilon')} {\tanh 1
  \cos\epsilon'}=\frac{\tan \epsilon'}{\tanh 1} \;,
$$
qui tend uniform\'ement vers $0$ quand $\epsilon'$ tend vers $0$. Le
r\'esultat en d\'ecoule. \cqfd

\bigskip \noindent{\bf D\'emonstration du th\'eor\`eme
  \ref{theo:comptage}. }
Fixons $\epsilon_0>0$ suffisamment petit.

\medskip Nous effectuons tout d'abord quelques r\'eductions pour
simplifier la situation g\'eom\'etrique de l'\'enonc\'e.

Si $C$ est r\'eunion de deux sous-vari\'et\'es immerg\'ees
de volume fini totalement g\'eod\'esiques $C'$ et $C''$, alors les
deux membres de l'\'equivalent du th\'eor\`eme sont les sommes des
termes correspondants \`a $C'$ et $C''$. Nous pouvons donc supposer
que $C$ est connexe.

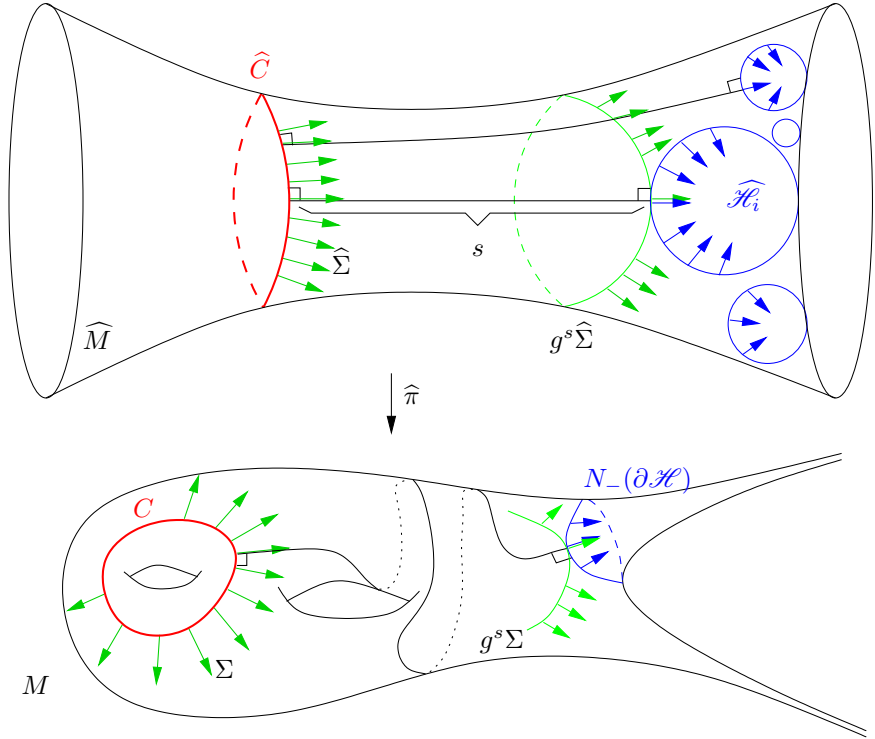
\begin{center}
\begin{picture}(0,0)%
\includegraphics{fig_revetement.pstex}%
\end{picture}%
\setlength{\unitlength}{3398sp}%
\begingroup\makeatletter\ifx\SetFigFontNFSS\undefined%
\gdef\SetFigFontNFSS#1#2#3#4#5{%
  \reset@font\fontsize{#1}{#2pt}%
  \fontfamily{#3}\fontseries{#4}\fontshape{#5}%
  \selectfont}%
\fi\endgroup%
\begin{picture}(6316,5379)(623,-8038)
\put(1171,-5191){\makebox(0,0)[lb]{\smash{{\SetFigFontNFSS{10}{12.0}{\rmdefault}{\mddefault}{\updefault}{\color[rgb]{0,0,0}$\wh M$}%
}}}}
\put(3511,-5596){\makebox(0,0)[lb]{\smash{{\SetFigFontNFSS{10}{12.0}{\rmdefault}{\mddefault}{\updefault}{\color[rgb]{0,0,0}$\wh \pi$}%
}}}}
\put(4006,-4516){\makebox(0,0)[lb]{\smash{{\SetFigFontNFSS{10}{12.0}{\rmdefault}{\mddefault}{\updefault}{\color[rgb]{0,0,0}$s$}%
}}}}
\put(2133,-7598){\makebox(0,0)[lb]{\smash{{\SetFigFontNFSS{10}{12.0}{\rmdefault}{\mddefault}{\updefault}{\color[rgb]{0,0,0}$\Sigma$}%
}}}}
\put(4822,-6201){\makebox(0,0)[lb]{\smash{{\SetFigFontNFSS{10}{12.0}{\rmdefault}{\mddefault}{\updefault}{\color[rgb]{0,0,1}$N_-(\partial\H)$}%
}}}}
\put(4576,-5183){\makebox(0,0)[lb]{\smash{{\SetFigFontNFSS{10}{12.0}{\rmdefault}{\mddefault}{\updefault}{\color[rgb]{0,0,0}$g^s\wh \Sigma$}%
}}}}
\put(2986,-4644){\makebox(0,0)[lb]{\smash{{\SetFigFontNFSS{10}{12.0}{\rmdefault}{\mddefault}{\updefault}{\color[rgb]{0,0,0}$\wh \Sigma$}%
}}}}
\put(2386,-3215){\makebox(0,0)[lb]{\smash{{\SetFigFontNFSS{10}{12.0}{\rmdefault}{\mddefault}{\updefault}{\color[rgb]{1,0,0}$\wh C$}%
}}}}
\put(5866,-4160){\makebox(0,0)[lb]{\smash{{\SetFigFontNFSS{10}{12.0}{\rmdefault}{\mddefault}{\updefault}{\color[rgb]{0,0,1}${\wh \H}_i$}%
}}}}
\put(1538,-6425){\makebox(0,0)[lb]{\smash{{\SetFigFontNFSS{10}{12.0}{\rmdefault}{\mddefault}{\updefault}{\color[rgb]{1,0,0}$C$}%
}}}}
\put(4066,-7378){\makebox(0,0)[lb]{\smash{{\SetFigFontNFSS{10}{12.0}{\rmdefault}{\mddefault}{\updefault}{\color[rgb]{0,0,0}$g^s\Sigma$}%
}}}}
\put(721,-7711){\makebox(0,0)[lb]{\smash{{\SetFigFontNFSS{10}{12.0}{\rmdefault}{\mddefault}{\updefault}{\color[rgb]{0,0,0}$M$}%
}}}}
\end{picture}%

\end{center}

Pour tout $h>0$, soit $\H_h$ la r\'eunion des voisinages cuspidaux des
m\^emes bouts que $\H$, contenue dans $\H$, et telle que $\partial
\H_h$ soit une hypersurface immerg\'ee \`a \'equidistance $h$ de
$\partial \H$.  Avec notre convention concernant la longueur
alg\'ebrique des segments perpendiculaires communs, nous avons
$\N_{M,C,\H_h}(t)=\N_{M,C,\H}(t+h)$ puisque les g\'eod\'esiques de
$\HH^n_\RR$ perpendiculaires \`a une horosph\`ere sont
perpendiculaires \`a toutes les horosph\`eres concentriques. Puisque
le flot g\'eod\'esique contracte par $e^{-t}$ la distance riemannienne
induite des feuilles stables, $\Vol(\H_h)=e^{-(n-1)h}
\Vol(\H)$. Quitte \`a remplacer $\H$ par $\H_h$ avec $h$ assez grand,
nous pouvons donc supposer que les voisinages cuspidaux qui composent
$\H$ sont \`a distance strictement sup\'erieure \`a $2\epsilon_0$ deux
\`a deux, ainsi que de $C$ priv\'e d'un voisinage cuspidal dans $C$ de
chaque pointe de $C$ repr\'esent\'ee par un rayon g\'eod\'esique
minimizant contenu dans $\H$.

Nous commen\c{c}ons la preuve du th\'eor\`eme \ref{theo:comptage} en
introduisant divers objets g\'eom\'etriques qui nous serons utiles.

Notons $\wh \pi:\wh M\ra M$ le rev\^etement riemannien de $M$ d\'efini
par $C$~: fixons un point base dans $C$ et prenons pour point base
dans $M$ son image par l'immersion $\iota$ de $C$ dans $M$; alors, en
notant $\iota_*:\pi_1C\ra \pi_1M$ l'application induite sur les
groupes fondamentaux par $\iota$, le rev\^etement $\wh \pi:\wh M\ra M$
est le rev\^etement connexe de $M$ d\'efini par le sous-groupe
$\iota_*\pi_1C$ de $\pi_1M$. Nous munissons bien s\^ur $\wh M$ de la
m\'etrique riemannienne relev\'ee. Par le th\'eor\`eme du
rel\`evement, l'immersion de $C$ dans $M$ se rel\`eve en un plongement
de $C$ dans $\wh M$, dont nous notons l'image $\wh C$.  Nous noterons
encore $(g^t)_{t\in\RR}$ le flot g\'eod\'esique sur $T^1\wh M$. Notons
$\wh \Sigma=\nu^1\wh C$ (voir la figure ci-dessus), et d\'efinissons
$$
\wh \H=\wh\pi^{-1}(\H)\;\;\;{\rm et}\;\;\;
\wh \Omega_t=\bigcup_{s\in[0,t]}\;g^s \wh \Sigma\;.
$$
Les propri\'et\'es de ces objets g\'eom\'etriques dont nous aurons
besoin sont les suivantes.

\begin{itemize}
\item[$\bullet$] La partie $\wh C$ est une sous-vari\'et\'e
  (plong\'ee) de volume fini, connexe, totalement g\'eod\'esique, de
  dimension $k$, dans la vari\'et\'e hyperbolique $\wh M$, et $\wh
  \Sigma$ est une sous-variété (plongée) de $T^1\wh M$.
\item[$\bullet$] La vari\'et\'e $\wh M$ est g\'eom\'etriquement finie, et
  ses pointes correspondent aux pointes de $\wh C$. Elle se r\'etracte
  par d\'eformation forte sur $\wh C$ (le long des rayons
  g\'eod\'esiques perpendiculaires \`a $\wh C$), de sorte que chaque
  voisinage cuspidal de $\wh M$ soit pr\'eserv\'e lors de la
  d\'eformation et s'envoit sur un voisinage cuspidal de $\wh C$.
\item[$\bullet$] L'application $\wh \pi$ envoie $\wh C$ sur (l'image
  dans $M$ de) $C$ en pr\'eservant les mesures riemanniennes; en
  particulier
\begin{equation}\label{volCChat}
\Vol(\wh C)=\Vol(C)\;.
\end{equation}
\item[$\bullet$] L'image de $\wh \Sigma$ par l'application tangente
  $T\wh \pi$ est l'image dans $T^1M$ du fibré normal unitaire de
  $C$. Puisque $\wh\pi$ est un rev\^etement riemannien, $T\wh \pi$
  envoie $g^t\wh\Sigma$ sur $g^t\Sigma$ en pr\'eservant les
  mesures riemanniennes~:
\begin{equation}\label{eq:mememes} 
(T\wh\pi)_*(d\vol_{g^t\wh \Sigma})=d\vol_{g^t\Sigma}\;.
\end{equation}
\item[$\bullet$] La partie $\wh \Omega_t$ est une sous-vari\'et\'e \`a
  bord, lisse, de volume fini, de dimension $n$ dans la vari\'et\'e
  $T^1\wh M$ de dimension $2n-1$, qui est la réunion disjointe des
  $g^s \wh \Sigma$ pour $s\in[0,t]$.
\item[$\bullet$] L'application de $\wh \Omega_t$ dans $[0,t]$, de
  fibre $g^s \wh \Sigma$ au dessus de $s\in[0,t]$, est une fibration
  (diff\'erentiablement) triviale; puisque les g\'eod\'esiques de
  vecteur origine dans $\wh \Sigma$ sont perpendiculaires aux
  hypersurfaces \'equidistantes \`a $\wh \Sigma$, la mesure
  riemannienne de $\wh \Omega_t$ se d\'esint\`egre par cette fibration
  en la mesure de Lebesgue sur $[0,t]$, de mesure conditionnelle
  au-dessus du point $s$ la mesure riemannienne sur $g^s\wh \Sigma$:
\begin{equation}\label{eq:prepacesaro} 
d\vol_{\Omega_t}=\int_{0}^t d\vol_{g^s\wh \Sigma}\;ds\;.
\end{equation}
\item[$\bullet$] Les courbures principales des hypersurfaces
  $\partial\V_t\wh C$ convergent uniform\'ement, quand $t$ tend vers
  $+\infty$, vers $1$ (la courbure principale des horosph\`eres dans
  $\HH^n_\RR$).  Par le point pr\'ec\'edent et la proposition
  \ref{prop:pansu} (appliqu\'ee \`a la vari\'et\'e $\wh M$ et \`a
  l'hypersurface $\partial\V_t\wh C$), nous avons donc, quand $t$ tend
  vers $+\infty$, l'\'equivalent
\begin{equation}\label{eq:VolOmegat}
\Vol(\wh \Omega_t)
\sim 2^{\frac{n-1}{2}}\Vol (\V_{t}\wh C)\;.
\end{equation}
\item[$\bullet$] La partie $\wh \H$ est la r\'eunion d'une union finie
  $V$ de voisinages cuspidaux de pointes de $\wh M$ \`a distance
  strictement sup\'erieure \`a $2\epsilon_0$ deux \`a deux, et d'une
  union d'horoboules isom\'etriquement plong\'ees dans $\wh M$, ces
  horoboules plong\'ees \'etant \`a distance strictement sup\'erieure
  \`a $2\epsilon_0$ deux \`a deux, ainsi que de $\wh C$ et de $V$. De
  plus, $\partial (\wh \H-V)=\wh\pi^{-1}(\partial \H)-\partial V$, que
  nous noterons $\partial^* \wh \H$ pour simplifier les notations dans
  la suite, est une sous-vari\'et\'e (plong\'ee) strictement convexe
  lisse de dimension $n-1$ dans $\wh M$, donc $N_-(\partial^*\wh \H)$
  est une sous-vari\'et\'e (plong\'ee) lisse de dimension $n-1$ dans
  $T^1\wh M$.
\end{itemize}

\medskip L'int\'er\^et d'introduire ces objets g\'eom\'etriques pour
notre probl\`eme de comptage r\'eside dans le r\'esultat suivant.

\bprop \label{prop:relatcomptgeom}
Pour tout $t\geq 0$, l'ensemble $I_t=\wh \Omega_t\cap N_-(\partial
\wh \H)$ est fini, de cardinal \'egal \`a $\N(t)$.  
\eprop

\dem Par stricte convexit\'e des horosph\`eres dans $\HH^n_\RR$ (et
additivit\'e des dimensions), chaque composante connexe de
$N_-(\partial^* \wh \H)$ rencontre $\wh \Omega_t$ en au plus un point
(et l'intersection est transverse). Comme aucune composante de
$N_-(\partial V)$ ne rencontre $\wh \Omega_t$, nous avons $I_t=\wh
\Omega_t\cap N_-(\partial^* \wh \H)$.  Pour usage ult\'erieur, notons
que deux éléments de $I_t$ sont donc \`a distance strictement
sup\'erieure \`a $2\epsilon_0$, car leurs points bases le sont.


Pour tout $t\geq 0$, tout segment de g\'eod\'esique $\ga$ de longueur
$s\leq t$, partant perpendiculairement en
un point de $C$ et arrivant perpendiculairement en un point de
$\partial \H$, de multiplicit\'e $m$, se rel\`eve par le rev\^etement
$\wh \pi$ en exactement $m$ segments g\'eod\'esiques de longueurs $s$
(d'origines deux \`a deux distinctes), partant perpendiculairement de
$\wh C$, et arrivant perpendiculairement en des points
$x_{\ga,1}\dots, x_{\ga,m}$ (deux \`a deux distincts) de $\partial^*
\wh \H$. Les vecteurs unitaires $N_-(x_{\ga,1}),\dots, N_-(x_{\ga,m})$
rentrants normaux \`a $\partial^* \wh \H$ appartiennent \`a $g^{s} \wh
\Sigma\cap N_-(\partial^* \wh \H)$, donc \`a $I_t$. Le r\'esultat en
d\'ecoule facilement. \cqfd

\begin{center}
\begin{picture}(0,0)%
\includegraphics{fig_interOmegwhH.pstex}%
\end{picture}%
\setlength{\unitlength}{3398sp}%
\begingroup\makeatletter\ifx\SetFigFontNFSS\undefined%
\gdef\SetFigFontNFSS#1#2#3#4#5{%
  \reset@font\fontsize{#1}{#2pt}%
  \fontfamily{#3}\fontseries{#4}\fontshape{#5}%
  \selectfont}%
\fi\endgroup%
\begin{picture}(4710,2904)(2326,-2053)
\put(2341,-1366){\makebox(0,0)[lb]{\smash{{\SetFigFontNFSS{10}{12.0}{\rmdefault}{\mddefault}{\updefault}{\color[rgb]{0,0,0}$\wh \Sigma$}%
}}}}
\put(7021,-511){\makebox(0,0)[lb]{\smash{{\SetFigFontNFSS{10}{12.0}{\rmdefault}{\mddefault}{\updefault}{\color[rgb]{0,0,0}$g^t\wh \Sigma$}%
}}}}
\put(4051,-1906){\makebox(0,0)[lb]{\smash{{\SetFigFontNFSS{10}{12.0}{\rmdefault}{\mddefault}{\updefault}{\color[rgb]{0,0,0}$\wh\Omega_t\subset T^1\wh M$}%
}}}}
\put(5597,-755){\makebox(0,0)[lb]{\smash{{\SetFigFontNFSS{10}{12.0}{\rmdefault}{\mddefault}{\updefault}{\color[rgb]{0,0,0}$U_v$}%
}}}}
\put(5344,-659){\makebox(0,0)[lb]{\smash{{\SetFigFontNFSS{10}{12.0}{\rmdefault}{\mddefault}{\updefault}{\color[rgb]{0,0,0}$v$}%
}}}}
\end{picture}%

\end{center}

\medskip Pour continuer la preuve du th\'eor\`eme \ref{theo:comptage},
nous introduisons maintenant l'objet analytique qui nous sera
utile. Montrons que pour tout $\eta>0$ suffisamment petit, il existe
une application $\psi_\eta\in\C_c(T^1M)$, \`a support dans
$\V_\eta(N_-(\partial\H))$, telle que les mesures
$\psi_\eta\;d\vol_{T^1M}$ convergent vaguement vers la mesure
riemannienne $d\vol_{N_-(\partial\H)}$ de $N_-(\partial\H)$ quand
$\eta$ tend vers $0$.  Ce r\'esultat est classique, mais nous en donnons
une preuve car nous aurons besoin d'une version assez pr\'ecise, et de
ses notations, par la suite.

Notons $\nu\H$ le fibr\'e normal de la sous-vari\'et\'e
$N_-(\partial\H)$ de la vari\'et\'e riemannienne $T^1M$, muni de la
m\'etrique fibrée induite par la métrique de $T^1M$. Par compacit\'e
de $N_-(\partial\H)$ et par les propri\'et\'es de l'exponentielle
riemannienne, il existe $\eta_0>0$ tel que pour tout $\eta\in\;
]0,\eta_0]$, l'application exponentielle normale de $\nu\H$ dans
$T^1M$ soit un diff\'eomorphisme lisse du $\eta$-voisinage de la
section nulle dans $\nu\H$ sur le $\eta$-voisinage de
$N_-(\partial\H)$ dans $T^1M$.  Pour tout $v$ dans $N_-(\partial\H)$,
nous noterons $D_\perp(v,\eta)$ l'image par cette application de la
boule de centre $0$ et de rayon $\eta$ dans la fibre de $\nu\H$
au-dessus de $v$. Pour tout $v$ dans $N_-(\partial\H)$, prenons sur
$D_\perp(v,\eta_0)$ une fonction continue $\psi_{v,\eta}$ \`a support
compact contenu dans $D_\perp(v,\eta)$ telle que
$\psi_{v,\eta}\;d\lambda$, o\`u $d\lambda$ est la mesure riemannienne
sur $D_\perp(v,\eta_0)$, converge vers la masse de Dirac en $v$ quand
$\eta\ra 0$. Par un argument de trivialit\'e locale du fibr\'e $\nu\H$
et de partition finie de l'unit\'e, nous obtenons donc une application
$\psi_\eta$ cherch\'ee, par orthogonalit\'e en $v\in N_-(\partial\H)$
de $D_\perp(v,\eta_0)$ et de $T_vN_-(\partial\H)$.

\bigskip Le coeur de la preuve du th\'eor\`eme \ref{theo:comptage} est
maintenant d'estimer de deux mani\`eres diff\'erentes l'int\'egrale
$\int \psi_\eta\circ T\wh \pi\;d\vol_{\wh \Omega_t}$. Fixons
$\epsilon>0$.

\medskip Donnons tout d'abord des informations sur le support de la
mesure $\psi_\eta\circ T\wh \pi\;d\vol_{\wh \Omega_t}$.

Notons $\wh\Omega_\infty=\bigcup_{t\in[0,+\infty[}\wh\Omega_t$ et
$I_\infty= \bigcup_{t\in[0,+\infty[} I_t$, qui sont des unions
croissantes. Remarquons que $\wh\Omega_\infty$ est une
sous-vari\'et\'e lisse de $T^1\wh M$, et que $I_\infty$ est un
ensemble de points $(2\epsilon_0)$-s\'epar\'es de $\wh\Omega_\infty$,
comme mentionn\'e dans la preuve de la proposition
\ref{prop:relatcomptgeom}. Remarquons aussi que si $\eta$ est assez
petit, alors $\wh\Omega_\infty$ ne rencontre pas le $\eta$-voisinage
de $N_-(\partial V)$, car l'angle en un point de $\partial V$ entre un
vecteur unitaire normal rentrant dans $V$ et le vecteur tangent d'un
rayon géodésique perpendiculaire à $\wh C$, en ce point, est au moins
$\frac{\pi}{2}$.

Pour tout $v$ dans $I_\infty$, notons $s_v$ l'unique $s\in[0,+\infty[$
tel que $v\in g^s\wh\Sigma$, et $U_v$ l'intersection triple du
support de $\psi_\eta\circ T\wh \pi$, de la boule de centre $v$ et de
rayon $\epsilon_0$, et de $\wh\Omega_\infty$. Comme il vient d'\^etre
dit, les $U_v$ sont deux \`a deux disjoints. Les propri\'et\'es des
$U_v$ que nous utiliserons sont les suivantes (voir la figure
pr\'ec\'edente).

\blemm\label{lem:technique}
(1) Pour tout $\eta>0$ suffisamment petit, pour tout $t\geq \epsilon_0$,
l'intersection de $\wh\Omega_{t-\epsilon_0}$ et du support de
$\psi_\eta\circ T\wh \pi$ est contenue dans la r\'eunion disjointe des
$U_v$ pour $v\in I_t$. De m\^eme, l'intersection de
$\wh\Omega_{t+\epsilon_0}$ et du support de $\psi_\eta\circ T\wh \pi$
contient la r\'eunion disjointe des $U_v$ pour $v\in I_t$.

(2) Pour tout $\eta>0$ suffisamment petit, pour tout $v\in I_\infty$, 
pour tout $t\geq s_v+\epsilon_0$, 
\begin{equation}\label{eq:majopetitdisc}
 1-\epsilon\leq \int_{U_v}\psi_\eta\circ T\wh\pi \;d\vol_{\Omega_t}
\leq 1+\epsilon\;.
\end{equation}
\elemm

\dem Tout d'abord, par le lemme \ref{lem:procherentranthoro}, si
$\eta$ est assez petit, tout point de l'intersection du support de
$\psi_\eta\circ T\wh \pi$ et de $\wh\Omega_\infty$ est contenu dans la
boule de centre $v$ et de rayon $\min\{\epsilon_0,\eta_0\}$ pour un
$v\in I_\infty$.  La premi\`ere affirmation en d\'ecoule.

Pour $s$ assez grand, la sous-vari\'et\'e $g^s\wh \Sigma$ de $T^1\wh
M$ est, au voisinage de chacun de ses points, proche d'un germe de
feuille instable, donc par la proposition \ref{prop:pansujacobi}, est
presque perpendiculaire au germe de feuille stable passant par ce
point.  En tout point de $T^1M$, l'orbite du flot g\'eod\'esique est
orthogonale aux feuilles stables et instables. Donc si $\eta_0$ est
suffisamment petit, il existe $t_0\geq 0$ tel que si $v\in I_\infty$
v\'erifie $s_v\geq t_0$, alors $\wh\Omega_\infty\cap B(v,\eta_0)$ est
(diff\'erentiablement) proche de $D_\perp(v,\eta_0)$.  Comme
l'int\'egrale riemannienne de $\psi_\eta\circ T\wh \pi$ sur
$D_\perp(v,\eta_0)$ est \'egale \`a $1$ (puisque $T\wh \pi$ est un
rev\^etement riemannien), ceci montre la seconde affirmation lorsque
$s_v\geq t_0$.  Comme la restriction de $\psi_\eta$ \`a
$D_\perp(v,\eta_0)$ converge vers la masse de Dirac unit\'e en $v$,
quitte \`a choisir $\eta$ suffisamment petit, nous pouvons supposer
que cette estimation soit aussi vraie pour l'ensemble fini des $v\in
I_\infty$ tels que $s_v\in[0,t_0]$.  \cqfd

\medskip Proc\'edons maintenant \`a la double estimation de
l'int\'egrale globale $\int \psi_\eta\circ T\wh \pi\;d\vol_{\wh
  \Omega_t}$.

Respectivement par le lemme \ref{lem:technique} (1), par le lemme
\ref{lem:technique} (2), et par la proposition
\ref{prop:relatcomptgeom}, nous avons
$$
\int \psi_\eta\circ T\wh \pi\;d\vol_{\wh \Omega_{t-\epsilon_0}}\leq
\sum_{v\in I_t}\int_{U_v}\psi\circ T\wh \pi\;d\vol_{\wh \Omega_{t}}\leq 
(1+\epsilon)\operatorname{Card} I_t=(1+\epsilon)\N(t)\;.
$$
Avec un argument similaire pour la minoration, nous obtenons
l'encadrement suivant, reliant la fonction de comptage \`a l'int\'egrale
$\int \psi_\eta\circ T\wh \pi\;d\vol_{\wh \Omega_{t}}$~: pour tout
$t\geq \epsilon_0$,
\begin{equation}\label{eq:encadrun}
\frac{1}{1+\epsilon}\int \psi_\eta\circ T\wh \pi\;
d\vol_{\wh \Omega_{t-\epsilon_0}}\leq \N(t)\leq
\frac{1}{1-\epsilon}\int \psi_\eta\circ T\wh \pi\;
d\vol_{\wh \Omega_{t+\epsilon_0}}\;.
\end{equation}

Nous allons maintenant utiliser un argument d'\'equidistribution pour
donner un \'equivalent de l'int\'egrale $\int \psi_\eta\circ T\wh
\pi\;d\vol_{\wh \Omega_{t}}$ quand $t$ tend vers $+\infty$.

Par la formule \eqref{eq:mememes} et le th\'eor\`eme
d'\'equidistribution \ref{theo:equidistributionenhaut}, nous avons
$$
\lim_{t\ra+\infty}\oint\psi_\eta\circ T\wh \pi\;
d\vol_{g^t\wh \Sigma}=
\lim_{t\ra+\infty}\oint\psi_\eta\;
d\vol_{g^t\Sigma}=\oint \psi_\eta\;
d\vol_{T^1M}\;.
$$
Par la formule \eqref{eq:prepacesaro} et un argument de moyenne de
C\'esaro, nous avons donc
\begin{equation}\label{eq:limequidistr}
\lim_{t\ra+\infty}\oint\psi_\eta\circ T\wh \pi\;
d\vol_{\wh \Omega_t}=\oint \psi_\eta\;
d\vol_{T^1M}\;.
\end{equation}
Or par d\'efinition de l'application $\psi_\eta$,  nous avons
$$
\lim_{\eta\ra 0} \oint \psi_\eta\;
d\vol_{T^1M}=\frac{\Vol(N_-(\partial\H))}{\Vol(T^1M)} >0\;.
$$
Fixons $\eta>0$ suffisamment petit tel que
\begin{center}
  $(1-\epsilon)\frac{\Vol(N_-(\partial\H))}{\Vol(T^1M)}\leq
  {\displaystyle \oint \psi_\eta\; d\vol_{T^1M}}\leq
  (1+\epsilon)\frac{\Vol(N_-(\partial\H))}{\Vol(T^1M)}$.
\end{center}
De ceci, par les formules \eqref{eq:encadrun} et
\eqref{eq:limequidistr}, il d\'ecoule que, pour tout $t$ assez grand,
nous avons
\begin{center}
$
\frac{(1-\epsilon)^2}{1+\epsilon}\;
\frac{\Vol(N_-(\partial\H))}{\Vol(T^1M)}\Vol(\wh\Omega_{t-\epsilon_0})
\leq\N(t)\leq 
\frac{(1+\epsilon)^2}{1-\epsilon}\;
\frac{\Vol(N_-(\partial\H))}{\Vol(T^1M)}\Vol(\wh\Omega_{t+\epsilon_0})
$.
\end{center}
Donc, par les formules \eqref{eq:VolOmegat} et
\eqref{eq:equivvolplusepsilon}, et puisque $\epsilon$ et $\epsilon_0$
peuvent \^etre pris arbitrairement petits, nous avons l'\'equivalent
quand $t$ tend vers $+\infty$,
$$
\N(t)\sim\frac{2^{\frac{n-1}{2}}\Vol(N_-(\partial\H))}{\Vol(T^1M)}
\Vol(\V_t\wh C)\;.
$$
Par la proposition \ref{prop:pansu} et par le lemme
\ref{lem:volumecuspvsbord}, nous avons 
$$
\Vol(N_-(\partial\H))=
2^{\frac{n-1}{2}}\Vol(\partial\H)=(n-1)\;2^{\frac{n-1}{2}}\Vol(\H)\;.
$$
Le th\'eor\`eme \ref{theo:comptage} d\'ecoule alors de la formule
\eqref{eq:volTunMM} et de la proposition \ref{prop:computvol}
d'estimation de volumes de voisinages tubulaires.  \cqfd

\medskip
Nous concluons cette partie par diverses remarques sur le
th\'eor\`eme \ref{theo:comptage}.

\medskip
\noindent{\bf Remarques. } (1) Si $C$ est de codimension $1$ et
transversalement orientable, nous pouvons compter les perpendiculaires
communes \`a $C$ et \`a $\partial\H$ qui partent d'un c\^ot\'e donn\'e
de $C$. Fixons une des deux composantes connexes $\Sigma$ de $\nu^1
C$. Si $\iota:C\ra M$ est l'immersion consid\'er\'ee de $C$ dans $M$,
si $\alpha$ est un segment perpendiculaire commun \`a $C$ et \`a
$\partial\H$, si $v$ est le vecteur tangent (unitaire) de $\alpha$ en
son origine $x$ dans $C$, alors appelons {\it multiplicit\'e positive}
de $\alpha$ l'entier
$$
m_+(\alpha)=\operatorname{Card}\{x'\in C\;:\;
\iota(x')=x,\;\; T_{x'}\,\iota\,(T_{x'}C)\perp v, 
\;\;{\rm et}\;\;v\in\Sigma\}\;.
$$
La fonction de comptage associ\'ee est d\'efinie, pour tout $t\geq 0$,
en posant $\N_+(t)$ le nombre compt\'e avec multiplicit\'e positive de
perpendiculaires communes entre $C$ et $\partial\H$, de longueur
alg\'ebrique au plus $t$.  Comme toute composante connexe de $\nu^1 C$
s'\'equidistribue quand on la pousse par le flot g\'eod\'esique, on
peut v\'erifier que la m\'ethode utilis\'ee dans la preuve du
th\'eor\`eme \ref{theo:comptage} s'applique (en remplaçant $\wh
\Sigma$ par la composante connexe de $\nu^1 \wt C$ qui revêt
$\Sigma$). Le volume de la partie du $t$-voisinage de $\wh C$ d'un
c\^ot\'e donn\'e de $\wh C$ \'etant la moiti\'e du volume de tout ce
$t$-voisinage, nous en d\'eduisons donc le r\'esultat suivant.

\bcoro 
Soit $(M,C,\H)$ comme dans l'\'enonc\'e du th\'eor\`eme
\ref{theo:comptage}, avec $C$ de codimension $1$ et transversalement
orientable.  Alors quand $t$ tend vers $+\infty$,
$$
\N_+(t)\sim 
\frac{\Vol(\SSS_{n-k-1})\Vol(\H)\Vol(C)}{2\Vol(\SSS_{n-1})\Vol(M)}
\;e^{(n-1)t}\;.\;\;\;\mbox{\cqfd}
$$ 
\ecoro

\medskip(2) Soit $(M,C,\H)$ comme dans l'\'enonc\'e
du th\'eor\`eme \ref{theo:comptage}. Soient $m\in\NN-\{0\}$ et
$p:M'\ra M$ un rev\^etement de $M$ \`a $m$ feuillets. Si $\iota:C\ra
M$ est l'immersion de $C$ dans $M$, notons $p':C'\ra C$ le
rev\^etement riemannien de $C$ \`a $m$ feuillets, et $\iota':C'\ra M'$
l'immersion (totalement g\'eod\'esique) tels que le diagramme suivant
commute:
$$
\begin{array}{ccc}
C'&\stackrel{\iota'}{\longrightarrow}& M'\\
_{p'}\downarrow\;\;\;\; & & \;\;\downarrow{}_p\\
C&\stackrel{\iota}{\longrightarrow}& M
\end{array}
$$
Posons $\H'=p^{-1}(\H)$. Alors $\H'$ est une r\'eunion de voisinage
cuspidaux de pointes de $M'$. Il est facile de v\'erifier que
$\Vol(M')=m\Vol(M)$, $\Vol(\H')=m\Vol(\H)$, $\Vol(C')=m\Vol(C)$, et,
par d\'efinition des multiplicit\'es, que
$\N_{M',C',\H'}=m\,\N_{M,C,\H}$. Ceci est bien s\^ur compatible avec
la formule asymptotique du th\'eor\`eme \ref{theo:comptage}, et peut
permettre d'\'etendre ce r\'esultat aux cas des (bons) orbifolds, mais
nous pr\'ef\'erons une approche globale, plus utile pour nos
applications.

\medskip (3) Soient $n\geq 2$ et $\Ga$ un groupe discret de covolume
fini d'isom\'etries de $\HH^n_\RR$ (ayant \'eventuellement de la
torsion). Soient $p\geq 1$ et $C_1,\dots, C_p$ des sous-espaces
hyperboliques de dimension $k<n$ de $\HH^n_\RR$, de stabilisateurs
$\Ga_{C_1},\dots, \Ga_{C_p}$ de covolume fini.  Soient $q\geq 1$ et
$\H_1,\dots,\H_q$ des horoboules ouvertes de $\HH^n_\RR$ {\it
  pr\'ecis\'ement invariantes} par $\Ga$ (c'est-à-dire qu'en notant
$\Ga_{\H_j}$ le stabilisateur de $\H_j$, si $\ga\in\Ga$ et
$\H_j\cap\ga \H_j \neq\emptyset$ alors $\ga\in\Ga_{\H_j}$).

Pour tous $i\in\{1,\dots, p\}, j\in \{1,\dots,q\}$ et $\alpha,
\beta\in\Ga$, notons $\delta =\delta_{i,j,\alpha,\beta}$ l'arc
g\'eod\'esique de $\HH^n_\RR$, perpendiculaire commun entre $\alpha
C_i$ et $\beta\partial\H_j$ lorsqu'il existe (c'est-\`a-dire lorsque
le point \`a l'infini de $\beta\H_j$ n'est pas un point \`a l'infini
de $\alpha C_i$). Le groupe $\Ga$ agit naturellement sur l'ensemble
des tels arcs. Nous compterons la longueur $\ell(\delta)$ de $\delta$
positivement s'il sort de $\beta\H_j$ en son point de d\'epart sur
$\beta\partial \H_j$, et n\'egativement sinon. Appelons {\it
  multiplicit\'e} de $\delta$ (lorsqu'il existe) le rationnel
\begin{equation}\label{eq:defmultorbi}
m(\delta)=\frac{1}{\operatorname{Card}\big(\alpha\Ga_{C_i}\alpha^{-1}
\cap\beta\Ga_{\H_j}\beta^{-1}\big)}\;,
\end{equation}
inverse du cardinal du stabilisateur dans $\Ga$ de $\alpha
C_i\cup\beta \H_j$. (Cette multiplicit\'e est naturelle~: dans tout
probl\`eme de comptage d'objets ayant d'\'eventuelles sym\'etries, la
bonne fonction de comptage consiste \`a donner pour multiplicit\'e \`a
un objet l'inverse du cardinal de son groupe de sym\'etries.)

Puisque le stabilisateur dans $\Ga$ de $\alpha C_i\cup\beta \H_j$ fixe
point par point $\delta$, et par discr\'etude de $\Ga$, la
multiplicit\'e $m(\delta)$ est bien définie, et invariante par $\Ga$.
Par convention, posons $m(\delta_{i,j,\alpha,\beta})=0$ et
$\ell(\delta_{i,j,e,\ga})=-\infty$ si $\alpha C_i$ et
$\beta\partial\H_j$ n'ont pas d'arc perpendiculaire commun. Notons que
pour tous les $\ga\in\Ga$, $\alpha_i\in \Ga_{C_i}$ et
$\beta_j\in\Ga_{\H_j}$, l'arc $\delta_{i,j,\alpha,\beta}$ existe si et
seulement si l'arc $\delta_{i,j,\ga\alpha\alpha_i,\ga\beta\beta_j}$
existe, et alors $\ga\delta_{i,j,\alpha,\beta} =
\delta_{i,j,\ga\alpha\alpha_i,\ga\beta\beta_j}$.

Pour tout $t\geq 0$, notons $\N(t)=\N_{\Ga,(C_i),(\H_j)}(t)$ le
nombre, compt\'e avec multiplicit\'e, de classes modulo $\Ga$ d'arcs
g\'eod\'esiques $\delta_{i,j,\alpha,\beta}$ de longueur au plus $t$.
C'est aussi le nombre, compt\'e avec multiplicit\'e, de classes modulo
$\Ga_{\H_j}$ de perpendiculaires communes entre $\partial\H_j$ et les
$\alpha C_i$ de longueur au plus $t$, ou encore le nombre, compt\'e
avec multiplicit\'e, de classes modulo $\Ga_{\C_i}$ de
perpendiculaires communes entre $C_i$ et les $\beta\partial\H_j$ de
longueur au plus $t$~:
$$
\N(t)=\sum_{1\leq i\leq p,\;1\leq j \leq q,\;
  [\ga]\in\Ga_{\H_j}\backslash \Ga/\Ga_{\ga C_i},
  \;\ell(\delta_{i,j,\ga,e})\leq t} m(\delta_{i,j,\ga,e})\;.
$$

\bcoro\label{coro:comptagedoubleclasse} 
Pour tout $t\geq 0$, $\N(t)$
est bien d\'efini, est fini, et quand $t$ tend vers $+\infty$,
$$
\N(t)\sim\frac{\Vol(\SSS_{n-k-1})
\Big(\sum_{j=1}^q\Vol(\Ga_{\H_j}\backslash\H_j)\Big)
\Big(\sum_{i=1}^p\Vol(\Ga_{C_i}\backslash C_i)\Big)}
{\Vol(\SSS_{n-1})\Vol(\Ga\backslash\HH^n_\RR)} \;e^{(n-1)t}\;.
$$
\ecoro

\dem Par additivit\'e des deux membres de l'\'equivalent, nous
pouvons supposer que $p=q=1$, et nous notons $\H_\infty=\H_1$,
$\Ga_\infty=\Ga_{\H_1}$, $C_0=C_1$, $\Ga_0=\Ga_{C_1}$. Par le lemme
de Selberg, $\Ga$ \'etant de type fini (et $\RR$ de caract\'eristique
nulle), il existe un sous-groupe $\Ga'$ sans torsion, distingu\'e,
d'indice fini $m$ dans $\Ga$. 

Notons $M'=\Ga' \backslash \HH^n_\RR$, qui est une vari\'et\'e
hyperbolique de volume fini de dimension $n$, et $\pi':\HH^n_\RR\ra
M'$ la projection canonique. Puisque le lieu de ramification du
rev\^etement branch\'e $M'\ra \Ga \backslash \HH^n_\RR$ est de mesure
nulle, nous avons
$$
\Vol(M')=m\Vol(M)\;.
$$

Notons $\ga_1\H_\infty,\dots,\ga_{m_{\infty}}\H_\infty$ un syst\`eme
de repr\'esentants modulo $\Ga'$ des images par $\Ga$ de $\H_\infty$.
L'image de $\ga_j\H_\infty$ dans $M'$ par $\pi'$ est un voisinage
cuspidal d'une pointe $e_j$ de $M'$, car $\H_\infty$ est
pr\'ecis\'ement invariante sous $\Ga$. De plus, $e_j\neq e_{j'}$ si
$j\neq j'$, par d\'efinition du syst\`eme de repr\'esentants.  Notons
$\H'$ l'image dans $M'$ de la r\'eunion disjointe
$\bigcup_{j=1}^{m_\infty} \ga_j\H_\infty$, qui est ouverte. Par le
m\^eme argument que pr\'ec\'edemment, nous avons
$$
\Vol(\H')=m\Vol(\Ga_\infty\backslash\H_\infty)\;.
$$

Le groupe $\Ga_0$ agit \`a gauche par translations \`a droite $(g,
\Ga'\ga) \mapsto \Ga'\ga g^{-1}$ sur l'ensemble fini $\Ga' \backslash
\Ga$ d'ordre $m$. Notons $\{\Ga'\ga_1,\dots, \Ga'\ga_{m_0}\}$ un
syst\`eme de repr\'esentants des orbites~: puisque $\Ga_0$ est le
stabilisateur de $C_0$ dans $\Ga$, $\ga_1 C_0, \dots, \ga_{m_{0}}C_0$
est un syst\`eme de repr\'esentants modulo $\Ga'$ des images par $\Ga$
de $C_0$. Le stabilisateur de $\Ga'\ga_i$ pour cette action de $\Ga_0$
est $\Ga_0\cap\ga_i^{-1}\Ga' \ga_i$. Par la formule des classes, nous
avons
$$
m=\sum_{i=1}^{m_0}\operatorname{Card}\big(
\Ga_0/(\Ga_0\cap\ga_i^{-1}\Ga'\ga_i)\big)\;.
$$
Notons $C'=\coprod_{i=1}^{m_0} \Ga'_i\backslash (\ga_iC_0)$, o\`u
$\Ga'_i= \operatorname{Stab}_{\Ga'}(\ga_i C_0)=
\Ga'\cap\ga_i\Ga_0\ga_i^{-1}$, qui est une vari\'et\'e hyperbolique
de volume fini de dimension $k<n$. Les inclusions des $\ga_i C_0$ dans
$\HH^n_\RR$ pour $i=1,\dots, m_0$ induisent une immersion totalement
g\'eod\'esique de $C'$ dans $M'$. De plus
\begin{align*}
\Vol(C')=\sum_{i=1}^{m_0}\Vol\big(\Ga'_i\backslash (\ga_iC_0)\big)
&=\sum_{i=1}^{m_0}\;[\ga_i\Ga_0\ga_i^{-1}:\Ga'_i]\;
\Vol\big((\ga_i\Ga_0\ga_i^{-1})\backslash (\ga_iC_0)\big)\\
&=m\Vol(\Ga_0\backslash C_0)\;,
\end{align*}
par la formule des classes ci-dessus.

Pour toute vari\'et\'e lisse $N$, notons $\G_k(TN)\ra N$ le fibr\'e
lisse des sous-espaces vectoriels de dimension $k$ des espaces
tangents aux points de $N$.  Pour tout $t\geq 0$, notons $\E(t)$
l'ensemble des images dans $2^{M'}\times\G_k(TM')$, par l'application
induite par $\pi'$, des couples $(\delta,E)$ pour lesquels il existe
$\alpha,\beta\in\Ga$ tels que $\delta$ soit la perpendiculaire commune
entre $\alpha C_0$ et $\beta\partial\H_\infty$, de longueur au plus
$t$, et $E$ soit l'espace tangent \`a $\alpha C_0$ au point initial de
$\delta$ sur $\alpha C_0$. Par d\'efinition des multiplicit\'es (voir
la formule \eqref{eq:defmulpvariete}), et puisque deux sous-variétés
totalement géodésiques connexes d'une variété hyperbolique, qui ont
même espace tangent en un point commun, coïncident, nous avons
$$
 \operatorname{Card}\big(\E(t)\big)=\N_{M',C',\H'}(t)\;.
$$
Le groupe fini $\Ga'\backslash\Ga$ de cardinal $m$ agit sur $M'$, donc
sur $\E(t)$. Puisque $\Ga'$ agit librement sur $\HH^n_\RR$, le
cardinal du stabilisateur dans ce groupe de l'image d'un $(\delta,E)$
comme ci-dessus est \'egal au cardinal du stabilisateur dans $\Ga$ de
$(\delta,E)$, donc de $\alpha C_0\cup\beta \H_\infty$, puisque $C_0$
est totalement g\'eod\'esique et que toute g\'eod\'esique
perpendiculaire \`a une horoboule passe par son point \`a l'infini. Donc
le cardinal de l'orbite par $\Ga'\backslash\Ga$ de l'image de
$(\delta,E)$ est \'egal \`a
$\frac{m}{\operatorname{Card}\big(\alpha\Ga_{C_i}\alpha^{-1}
  \cap\beta\Ga_{\H_j}\beta^{-1}\big)}$. Par d\'efinition des
multiplicit\'es (voir la formule \eqref{eq:defmultorbi}), nous avons
donc
$$
 \operatorname{Card}\big(\E(t)\big)=m\;\N_{\Ga',(C_0),(\H_\infty)}(t)\;.
$$
Le r\'esultat d\'ecoule alors du th\'eor\`eme \ref{theo:comptage}
appliqu\'e \`a $M',C',\H'$, et des calculs pr\'ec\'edents des
volumes. \cqfd

\medskip (4) La derni\`ere remarque est une formule de comptage de
points fixes hyperboliques (dans une orbite donn\'ee) d'un groupe
discret de covolume fini d'isom\'etries de $\HH^n_\RR$, la
complexit\'e \'etant (l'inverse de) la distance entre les deux points
fixes de l'isom\'etrie hyperbolique. Nous nous en servirons pour les
applications arithm\'etiques qui suivent.  Nous laissons au lecteur le
soin d'\'enoncer le r\'esultat lorsque l'on remplace $\ga_0$ par un
nombre fini de tels \'el\'ements (deux \`a deux sans puissances
conjugu\'ees).
Si $x$ est l'un des deux points fixes d'une isom\'etrie hyperbolique
de $\HH^n_\RR$, nous noterons $x^\sigma$ l'autre point fixe.

\bcoro\label{coro:comptageorbitepointfixehyperbolique} Soit $n\geq 2$,
et soit $G$ un groupe discret de covolume fini d'isom\'etries du
mod\`ele du demi-espace sup\'erieur de $\HH^n_\RR$ ayant $\infty$
comme point fixe parabolique. Soit $\ga_0$ un \'el\'ement hyperbolique
de $G$, notons $x_0\in\partial_\infty\HH^n_\RR$ un de ses deux points
fixes, $\ell_0$ sa distance de translation, et $n_{G,0}$ l'indice de
$\ga_0^\ZZ$ dans le stabilisateur dans $G$ de son axe de
translation. Notons $\A_\infty$ le volume euclidien d'un domaine
fondamental mesurable de l'hyperplan affine de hauteur $1$ sous
l'action du stabilisateur $G_\infty$ du point $\infty$ dans $G$, et
$||\cdot||$ la norme euclidienne sur $\partial_\infty\HH^n_\RR -
\{\infty\} =\RR^{n-1}$. Alors quand $t$ tend vers $+\infty$,
$$
\sum_{x\in G_\infty\backslash G \{x_0,x_0^\sigma\},\;
||x-x^\sigma||\geq\frac{1}{t}}\;{\textstyle 
\frac{1}{\operatorname{Card}
(G_\infty\cap\operatorname{Stab}_G x)}}
\sim\frac{2^n\;\ell_0\;\Vol(\SSS_{n-2})\;\A_\infty}
{(n-1)\;n_{G,0}\;\Vol(\SSS_{n-1})\;\Vol(G\backslash\HH^n_\RR)}
\;t^{n-1}\;.
$$
\ecoro

\dem Appliquons le corollaire \ref{coro:comptagedoubleclasse} avec
$k=1$, $\Ga=G$, une seule horoboule $\H_\infty$ constitu\'ee des
points de $\HH^n_\RR$ dont la derni\`ere coordonn\'ee est strictement
sup\'erieure \`a $h$, o\`u $h$ est choisi assez grand pour que $\H_\infty$
soit pr\'ecis\'ement invariante, et une seule g\'eod\'esique $C_0$
\'egale \`a l'axe de translation de $\ga_0$. Nous avons
$$
\Vol(G_\infty\backslash\H_\infty)=\frac{\A_\infty}{(n-1)h^{n-1}}\;,
$$
par le calcul du lemme \ref{lem:volumecuspvsbord}.  En notant $G_0$ le
stabilisateur de $C_0$, nous avons $\Vol(G_0\backslash C_0) =
\frac{\ell_0}{n_{G,0}}$. Si $x$ et $y$ sont les extr\'emit\'es,
diff\'erentes du point $\infty$, d'une g\'eod\'esique $C_*$, alors par
un calcul imm\'ediat de distance hyperbolique, toujours en comptant
alg\'ebriquement les distances \`a $\partial\H_\infty$, nous avons
$d(\partial\H_\infty, C_*)=\log\frac{2h}{||y-x||}$. En particulier,
$||y-x||\geq 1/t$ si et seulement si $d(\partial\H_\infty, C_*)\leq
\log(2ht)$.

Pour tout $\ga\in G$, en notant $x$ et $x^\sigma$  les points \`a
l'infini de $\ga C_0$, nous avons $\operatorname{Card}\big(G_\infty
\cap \operatorname{Stab}_G(x)\big)= \frac{1}{2}\operatorname{Card}
\big(G_\infty\cap \ga G_0\ga^{-1}\big)$ s'il existe un \'el\'ement de
$G_\infty$ \'echangeant $x$ et $x^\sigma$, et $\operatorname{Card}
\big(G_\infty\cap \operatorname{Stab}_\Ga(x)\big) =\operatorname{Card}
\big(G_\infty\cap \ga G_0\ga^{-1}\big)$ sinon, auquel cas $x$ et
$x^\sigma$ d\'efinissent deux \'el\'ements distincts de $G_\infty
\backslash G \{x_0,x_0^\sigma\}$. Ainsi,
$$
\sum_{x\in G_\infty\backslash G \{x_0,x_0^\sigma\},\;
||x-x^\sigma||\geq\frac{1}{t}}\;
\frac{1}{\operatorname{Card}
(G_\infty\cap\operatorname{Stab}_Gx)}$$
$$=
2\sum_{[\ga]\in G_\infty\backslash G/ G_0,\;
d(\partial\H_\infty, \ga C_0)\leq \log(2ht)}\;
\frac{1}{\operatorname{Card}
(G_\infty\cap\ga G_0\ga^{-1})}
$$
Le r\'esultat d\'ecoule alors du corollaire
\ref{coro:comptagedoubleclasse}. \cqfd

\medskip La plupart de nos applications arithm\'etiques d\'ecouleront du
corollaire \ref{coro:comptageorbitepointfixehyperbolique} de comptage
de points fixes hyperboliques, en prenant $n=2,3$ ou $5$.

\section{Applications arithm\'etiques}

Nous renvoyons par exemple \`a \cite{Samuel67, Narkiewicz90} pour
toutes les informations sur les corps de nombres que nous utiliserons,
et \`a \cite{Landau58, Buell89, BucVol07} pour celles sur les formes
quadratiques binaires.

\subsection{Irrationnels quadratiques et formes quadratiques binaires}
\label{subsec:justifcomplexiteh}

Nous rappelons dans cette partie la bijection bien connue entre
l'ensemble des r\'eels irrationnels quadratiques $\alpha$ et
l'ensemble des formes quadratiques binaires $Q$ enti\`eres primitives
ind\'efinies ne repr\'esentant pas $0$ sur $\QQ$. Comme promis dans
\cite{ParPau10MZ}, nous donnons des interpr\'etations alg\'ebriques de
la complexit\'e $h(\alpha)=\frac{2}{|\alpha-\alpha^\sigma|}$, en
l'exprimant en fonction du discriminant de sa forme quadratique
binaire associ\'ee, et, lorsque $\alpha$ est un entier alg\'ebrique,
en fonction de sa trace et de sa norme.

\bigskip Notons $K$ le corps des rationnels $\QQ$ (respectivement une
extension quadratique imaginaire de $\QQ$), et $\widehat K=\RR$
(respectivement $\wh K=\CC$). Soit $\OOO_K$ l'anneau des entiers de $K$.
Pour tout $\alpha\in\wh K$ irrationnel quadratique sur $K$, notons
$\alpha^\sigma$ son conjugu\'e de Galois sur $K$,
$\Tr(\alpha)=\alpha+\alpha^\sigma$ sa trace (relative),
$N(\alpha)=\alpha\alpha^\sigma$ sa norme (relative), et
$$
h(\alpha)=\frac{2}{|\alpha-\alpha^\sigma|}\;.
$$ 
Le groupe $\operatorname{PSL}_2(\widehat K)$ agit sur $\PP^1(\widehat
K)= \widehat K\cup \{\infty\}$ par homographies. Son sous-groupe
$\Ga_K=\operatorname{PSL}_2(\OOO_K)$ pr\'eserve l'ensemble des \'el\'ements
de $\widehat K$ irrationnels quadratiques sur $K$, et son action
commute avec l'action de l'automorphisme de Galois.

Remarquons que $h(\alpha)$ est fini, non nul, et que
$h(\alpha^\sigma)=h(\alpha)$. Nous avons montr\'e dans
\cite{ParPau10MZ} que pour tout $\beta$ appartenant au stabilisateur
$\Ga_{K,\infty}$ du point $\infty$ dans $\Ga_K$, nous avons
$h(\beta\alpha)=h(\alpha)$, et qu'il n'existe qu'un nombre fini
d'orbites modulo $\Ga_{K,\infty}$ d'\'el\'ements $\alpha$ de $\widehat
K$ irrationnels quadratiques sur $K$ tels que $h(\alpha)\leq c$, pour
toute constante $c>0$.

Pour tout $\alpha$ dans $\wh K$, il est bien connu (voir par exemple
\cite[Lem.~6.2]{ParPau10MZ}) que $\alpha$ est irrationnel
quadratique sur $K$ si et seulement s'il existe un \'el\'ement
hyperbolique du r\'eseau $\Ga_K$ de $\operatorname{PSL}_2(\widehat K)$
ayant $\alpha$ pour point fixe, l'autre point fixe \'etant alors
$\alpha^\sigma$.

\bigskip Apr\`es ces rappels sur les irrationnels quadratiques, passons
aux rappels sur les formes quadratiques binaires.

Soit $ Q(X,Y)=AX^2+BXY+CY^2 $ une forme quadratique binaire enti\`ere
(sur $\QQ$).  Nous noterons $A=A_Q$ si pr\'eciser $Q$ est
n\'ecessaire.  La forme $Q$ est dite {\em primitive} si ses
coefficients $A$, $B$, $C$ sont premiers entre eux. Notons $
\discr=\discr_Q=B^2-4AC $ le {\it discriminant} de $Q$, qui est congru
\`a $0$ ou $1$ modulo $4$. La forme quadratique r\'eelle $Q$ est
ind\'efinie si et seulement si $\discr>0$. La forme $Q$ est dite {\it
  irr\'eductible} si elle ne v\'erifie pas l'une des assertions
\'equivalentes suivantes~:
\begin{itemize}
\item[$\bullet$] $Q$ est, à multiple rationnel près, le produit de
  deux formes lin\'eaires enti\`eres;
\item[$\bullet$] $Q$ repr\'esente $0$ sur $\QQ$;
\item[$\bullet$] $\discr$ est un carr\'e.
\end{itemize}
Remarquons que $A\neq 0$ si $\discr$ n'est pas un carr\'e.  Le groupe
$\SL 2\ZZ$ agit \`a droite sur l'ensemble des formes quadratiques
binaires enti\`eres par $(\ga,Q)\mapsto Q\circ\ga$, en pr\'eservant le
discriminant~: $ \forall\;\ga\in\SL 2\ZZ,\;D_{Q\circ\ga}=D_Q $.  Deux
formes quadratiques binaires enti\`eres sont dites {\it \'equivalentes}
si elles sont dans la m\^eme orbite sous $\SL 2\ZZ$.

Le r\'esultat suivant est bien connu (voir par exemple
\cite[Theo.~4.5.8 (2)]{BucVol07}), sa preuve (non directe) est un
pr\'etexte pour introduire des objets qui nous servirons
ult\'erieurement. Posons
$$
\alpha_Q=\frac{-B+\sqrt{\discr}}{2A}\;.
$$ 

\bprop\label{prop:bijqaudbinirrquad} L'application $Q\mapsto\alpha_Q$
est une bijection de l'ensemble des formes quadratiques binaires
enti\`eres primitives ind\'efinies irr\'eductibles dans l'ensemble des
r\'eels irrationnels quadratiques, telle que, pour tout $\ga\in \SL
2\ZZ$,
$$
\alpha_{-Q}=(\alpha_Q)^{\sigma}\;\;\;{\rm et}\;\;\;
\alpha_{Q\circ\ga}=\ga\alpha_Q\;.
$$
\eprop

\dem Puisque $Q$ est primitive, le sous-groupe $\autom(Q,\ZZ)$ de $\SL
2\ZZ$ des \'el\'ements fixant $Q(X,Y)=AX^2+BXY+CY^2$ est (voir par
exemple \cite[Theo.~202]{Landau58})
$$
\autom(Q,\ZZ)=\Big\{\ga_{Q,t,u}=
\begin{pmatrix}
\dfrac{t-Bu}2 & -Cu\\ Au & \dfrac{t+Bu}2
\end{pmatrix}:\;t,u\in\ZZ,\;\; t^2-\discr u^2=4\Big\}\;.
$$
Notons $(t_Q,u_Q)$ la solution fondamentale de l'\'equation de
Pell-Fermat $t^2-\discr u^2=4$.

L'image $\PSO=\autom(Q,\ZZ)/\{\pm \operatorname{id}\}$ de
$\autom(Q,\ZZ)$ dans $\operatorname{PSL}_2(\ZZ)$ est infinie cyclique,
engendr\'ee par $\pm\ga_Q=\pm\ga_{Q,t_Q,u_Q}$.  L'application
$Q\mapsto\pm\ga_Q$ est une bijection de l'ensemble des formes
quadratiques binaires enti\`eres primitives ind\'efinies
irr\'eductibles dans l'ensemble des \'el\'ements hyperboliques
primitifs de $\operatorname{PSL}_2(\ZZ)$ (voir
\cite[Prop.~1.4]{Sarnak82}), telle que, pour tout $\ga\in \SL 2\ZZ$,
$$
\ga_{-Q}=(\ga_Q)^{-1}\;\;\;{\rm et}\;\;\;
\ga_{Q\circ\ga}=\ga\ga_Q\ga^{-1}\;.
$$
Il est imm\'ediat de montrer que l'irrationnel quadratique $\alpha_Q$
est le point fixe attractif de $\pm\ga_Q$ dans $\PP^1(\RR)$, l'autre
point fixe \'etant le conjugu\'e de Galois $\alpha_Q^\sigma$. Remarquons
que $\PSO$ est le stabilisateur de $\alpha_Q$ dans
$\operatorname{PSL}_2(\ZZ)$.  

Comme l'application de l'ensemble des \'el\'ements hyperboliques
primitifs $\ga$ de $\operatorname{PSL}_2(\ZZ)$ dans l'ensemble des
r\'eels irrationnels quadratiques, qui \`a $\ga$ associe son point fixe
attractif, est une bijection, le r\'esultat en d\'ecoule.  \cqfd

\medskip
La complexit\'e $h(\alpha)$ d'un r\'eel irrationnel quadratique $\alpha$
v\'erifie alors les propri\'et\'es suivantes.

\blemm \label{lem:relathdiscr} Soit $Q$ une forme quadratique binaire
enti\`ere ind\'efinie irr\'eductible.

\noindent (1) Nous avons
$$
h(\alpha_Q)=\frac{2|A_Q|}{\sqrt{\discr_Q}}\;.
$$
(2) Pour tout $\ga=\pm\left(\begin{array}{cc}a & b \\ c & d
\end{array}\right)$ dans $\operatorname{PSL}_2(\ZZ)$, nous avons
$$
h(\ga\alpha_Q)=\frac{|Q(d,-c)|}{|A_Q|}\;h(\alpha_Q)=
\frac{2}{\sqrt{\discr_Q}}\;|Q(d,-c)|\;. 
$$
(3) Si $\alpha$ est un entier quadratique r\'eel, alors
$$
h(\alpha)=\frac{2}{\sqrt{\Tr(\alpha )^2-4 N(\alpha)}}\;.
$$
\elemm

\dem (1) La premi\`ere assertion d\'ecoule imm\'ediatement des
d\'efinitions de $h(\cdot)$ et de $\alpha_Q$.

(2) Notons pour simplifier $\alpha=\alpha_Q$. Remarquons tout d'abord
que $N(\alpha)=C/A$, et ${\rm Tr}(\alpha)=-B/A$, donc
\begin{equation}\label{eq:relatformbinnorm}
Q(X,Y)=A\big(X^2-{\rm Tr}(\alpha)XY+N(\alpha)Y^2\big)\;.
\end{equation}
Or
\begin{equation}\label{eq:computactgah}
\ga\alpha-\ga\alpha^\sigma  =
\frac{a\alpha+b}{c\alpha+d}-
\frac{a\alpha^\sigma+b}{c\alpha^\sigma+d} = 
\frac{\alpha-\alpha^\sigma}{c^2N(\alpha)
     +cd\Tr(\alpha)+d^2}\;.
\end{equation}
Par (1),  l'assertion (2) en d\'ecoule.

(3) La forme quadratique binaire $Q_\alpha(X,Y)= X^2-\Tr(\alpha)XY +
N(\alpha)Y^2$ (qui v\'erifie $Q_\alpha=Q_{\alpha^\sigma}$) est
enti\`ere, primitive, et n'est pas produit de formes lin\'eaires
enti\`eres, car $Q_\alpha(X,Y)=N(X-\alpha Y)$ sur $\ZZ^2$ et $\alpha$
est irrationnel. Son discriminant $D_{Q_\alpha}= \Tr(\alpha)^2 -4
N(\alpha)$ est strictement positif, car $\alpha$ est r\'eel
irrationnel. Il est imm\'ediat que $\alpha_{Q_\alpha}$ est \'egal \`a
$\alpha$ ou \`a $\alpha^\sigma$. Le r\'esultat d\'ecoule donc de (1).
\cqfd

\medskip \rem Si $Q$ est une forme quadratique binaire enti\`ere
ind\'efinie primitive irr\'eductible, alors $\alpha_Q$
est un entier alg\'ebrique si et seulement si $A_Q=1$.  En effet, par la
formule \eqref{eq:relatformbinnorm}, nous avons $A_Q=1$ si et
seulement si $Q(X,Y)= X^2-{\rm Tr}(\alpha)XY+N(\alpha)Y^2$, et un r\'eel
quadratique est un entier quadratique si et seulement si sa norme et
sa trace sont enti\`eres.

\medskip Nous dirons indiff\'eremment qu'une forme quadratique binaire
enti\`ere $Q$, qu'un r\'eel irrationnel quadratique $\alpha$, un
\'el\'ement hyperbolique $\gamma$ de $\operatorname{PSL}_2(\ZZ)$ est
{\it r\'eciproque} s'il est, respectivement, \'equivalent \`a $-Q$,
dans la m\^eme orbite que son conjugu\'e de Galois sous
$\operatorname{PSL}_2(\ZZ)$, conjugu\'e \`a son inverse dans
$\operatorname{PSL}_2(\ZZ)$.

\bprop \label{prop:equivreciproque}
Soit $Q$ une forme quadratique binaire enti\`ere primitive
ind\'efinie de discriminant non carr\'e. Les conditions suivantes sont
\'equivalentes:
\begin{enumerate}
\item[(1)] $Q$ est r\'eciproque;
\item[(2)] $\alpha_Q$ est r\'eciproque;
\item[(3)] la p\'eriode du d\'eveloppement en fraction continue de
  $\alpha_Q$ est un palindrome \`a permutation cyclique pr\`es;
\item[(4)] $\pm\gamma_Q$ est r\'eciproque;
\item[(5)] $\gamma_Q$ est conjugu\'e \`a une matrice sym\'etrique.
\end{enumerate}

Lorsque $\alpha_Q$ est un entier alg\'ebrique, la condition (2) est
\'equivalente \`a demander que l'ordre $\ZZ+\alpha_Q\ZZ$ contienne une
unit\'e de norme $-1$.  
\eprop

Voir aussi \cite{PolRud01,Long02,Burger05,Sarnak07} pour d'autres
caract\'erisations.

\medskip \dem L'\'equivalence entre (1), (2) et (4) d\'ecoule de la
proposition \ref{prop:bijqaudbinirrquad} et de sa preuve, le point
fixe attractif de l'inverse d'un \'el\'ement hyperbolique \'etant son
point fixe r\'epulsif. Pour l'\'equivalence entre (2) et (3), comme
d\'emontr\'e par exemple dans \cite[\S 23]{Perron13}, la partie
p\'eriodique du d\'eveloppement en fraction continue de
$\alpha^\sigma$ est, \`a permutation cyclique pr\`es, obtenue en
renversant l'ordre de la partie p\'eriodique du d\'eveloppement en
fraction continue de $\alpha$; de plus, deux nombres r\'eels sont dans
la m\^eme orbite par homographies enti\`eres si et seulement si les
queues de leur d\'eveloppement en fraction continue
co\"{\i}ncident. Pour l'\'equivalence entre (4) et (5), voir
\cite{Sarnak07}.

Pour montrer la derni\`ere assertion, remarquons d'une part que si
$\alpha$ est un r\'eel quadratique irrationnel, si $\ga=
\pm\left(\begin{array}{cc}a & b \\ c & d \end{array}\right)\in
\operatorname{PSL}_2(\ZZ)$ envoie $\alpha$ sur $\alpha^\sigma$, alors
par l'\'equation \eqref{eq:computactgah}, nous avons $N(d+\alpha_Q c)
=-1$. R\'eciproquement, si $\alpha$ est un entier quadratique
irrationnel, supposons que $\ZZ+\alpha\ZZ$ contienne un \'el\'ement de
norme $-1$.  Soient $a$ et $c$ dans $\ZZ$ tels que $N(c\alpha-a)=-1$,
c'est-\`a-dire que $c^2N(\alpha)-ac\Tr(\alpha)+a^2=-1$. En particulier,
$\alpha$ \'etant un entier quadratique, $c$ divise $a^2+1$. Posons
$d=-a$ et $b=-\frac{a^2+1}{c}$. Nous avons alors $ad-bc =1$ et
$cN(\alpha)-a\Tr(\alpha)-b=0$, c'est-\`a-dire que $\alpha^\sigma
(c\alpha+d)=a\alpha+b$. Donc $\alpha$ est r\'eciproque. \cqfd

\subsection{Comptage d'irrationnels quadratiques r\'eels 
dans une orbite de sous-groupes de congruences du 
groupe modulaire}
\label{subsec:comptirraquad}

Nous donnons dans cette partie un th\'eor\`eme de comptage
d'irrationnels quadratiques dans une orbite par homographies d'un
sous-groupe d'indice fini du groupe modulaire (\'etendu par
l'automorphisme de Galois).

\btheo\label{theo:appliarithneg2} Soit $\alpha_0$ un nombre r\'eel,
irrationnel quadratique sur $\QQ$. Soient $G$ un sous-groupe d'indice
fini du groupe modulaire $\Ga={\rm PSL}_2(\ZZ)$, et $\cdot$ son action
par homographies sur $\RR\cup\{\infty\}$.  Alors quand $s$ tend vers
$+\infty$,
$$
\operatorname{Card}\{\alpha\in G \cdot \{\alpha_0,\alpha_0^\sigma\} 
\!\!\!\!\mod q_G\;:\; h(\alpha)\leq s\} 
\sim \;\frac{12 \;q_G \;\log\big(\frac{|\operatorname{tr}\ga_0|+
\sqrt{\operatorname{tr}^2\ga_0-4}}{2}\big)}{\pi^2\;[\Ga:G]\;n_0}\;
   \;s\;,
$$
o\`u $q_G$ est le plus petit entier strictement positif $q$ tel que
l'homographie $z\mapsto z+q$ appartienne \`a $G$; $\ga_0$ est un
\'el\'ement non trivial de $G$ fixant $\alpha_0$; et $n_0$ est
l'indice de ${\ga_0}^\ZZ$ dans le stabilisateur de
$\{\alpha_0,\alpha_o^\sigma\}$ dans $G$.  
\etheo

Avant de donner la preuve de ce r\'esultat, voici quelques valeurs
possibles pour $q_G$ et $[\Ga:G]$. Soit $p\in\NN-\{0,1\}$. Si
$G=\Ga(p)$ est le sous-groupe de congruence principal modulo $p$,
c'est-\`a-dire le noyau du morphisme ${\rm PSL}_2(\ZZ)\ra{\rm PSL}_2
(\ZZ/p\ZZ)$ de r\'eduction modulo $p$ des coefficients, alors (voir
par exemple \cite[Theo.~5.5.4]{Katok92}, \cite[p.~22]{Shimura71}))
\begin{equation}\label{eq:indexprincipal}
  (q_G,[\Ga:G])= \begin{cases} (2,6)\;\;{\rm si}\;\;p=2,\\
    \Big(p,\;\frac{p^3}{2}{\displaystyle 
    \prod_{\ppp \;{\rm premier},\; \ppp | p}} 
    (1-\frac{1}{\ppp^2})\;\Big)\;\;{\rm si} \;\;p>2\;.
\end{cases}
\end{equation}
Et si $G=\Ga_0(p)=\Big\{\pm\left(\begin{array}{cc}a & b \\ c & d 
\end{array}\right)\in {\rm PSL}_2(\ZZ)\;:\;c\equiv 0\mod p\Big\}$ 
est la pr\'eimage du sous-groupe triangulaire sup\'erieur par le
morphisme de r\'eduction modulo $p$ des coefficients, alors (voir par
exemple \cite[p.~24]{Shimura71})
\begin{equation}\label{eq:indexnonprincipal}
(q_G,[\Ga:G])= \Big(1,p\prod_{\ppp \;{\rm premier},\; \ppp | p} 
(1+\frac{1}{\ppp})\Big)\;.
\end{equation}

\dem Remarquons tout d'abord que $q_G$ existe, car $z\mapsto z+1$
appartient \`a $\Ga$, et $G$ est d'indice fini dans $\Ga$. De m\^eme,
$\ga_0$ existe et est hyperbolique (donc v\'erifie
$|\operatorname{tr}\ga_0|>2$), par exemple par la démonstration de la
proposition \ref{prop:bijqaudbinirrquad} ou par
\cite[Lem.~6.2]{ParPau10MZ}, et $n_0$ est donc bien d\'efini.

Nous appliquons le corollaire
\ref{coro:comptageorbitepointfixehyperbolique} avec $n=2$, $\A_\infty
=q_G$ (tout segment de longueur euclidienne $q_G$ dans la droite
horizontale \`a hauteur $1$ \'etant un domaine fondamental de cette
droite sous l'action du stabilisateur $G_\infty$ du point $\infty$
dans $G$), et
$$
\Vol(G\backslash \HH^2_\RR)=[\Ga:G]\;\Vol(\PSL 2\ZZ
\backslash \HH^2_\RR)=\frac{\pi\; [\Ga:G]}{3}\;.
$$
Rappelons la formule 
\begin{equation}\label{eq:calclongtransl}
\cosh\frac{\ell_0}{2}=\frac{|\operatorname{tr}\ga_0|}{2}
\end{equation}
 bien connue (voir par exemple
\cite[page 173]{Beardon83}) de la distance de translation $\ell_0$ d'un
\'el\'ement hyperbolique $\ga_0$ de $\PSL 2\RR$.

Notons que $h(\alpha)\leq t$ si et seulement si $|\alpha-
\alpha^\sigma| \geq\frac{2}{t}$, et que $G_\infty$ \'etant sans
torsion, toutes les multiplicit\'es valent $1$. Le r\'esultat
d\'ecoule donc du corollaire
\ref{coro:comptageorbitepointfixehyperbolique}. \cqfd

\medskip Remarquons que si $\alpha_0$ est r\'eciproque, si $\ga_0$ est
primitif, et si $G=\Ga$, alors $q_G=1$ et $n_0=2$. Donc quand $s$ tend
vers $+\infty$,
$$
\operatorname{Card}\{\alpha\in \Ga \cdot \alpha_0 
\!\!\!\!\mod \ZZ\;:\; h(\alpha)\leq s\} 
\sim \;\frac{6 \;\log\big(\frac{|\operatorname{tr}\ga_0|+
\sqrt{\operatorname{tr}^2\ga_0-4}}{2}\big)}{\pi^2}\;
   \;s\;.
$$
Par exemple, le nombre d'or $\phi=(1+\sqrt 5)/2$ est un point fixe de
l'\'el\'ement hyperbolique primitif sym\'etrique $\begin{pmatrix} 2
  &1\\1& 1 \end{pmatrix}$ de $\Ga={\rm PSL}_2(\ZZ)$. Donc, quand $s$
tend vers $+\infty$,
$$
\operatorname{Card} \{\alpha\in \;\Ga\cdot \phi\!\!\mod\ZZ\;:\;
h(\alpha)\leq t\}\sim \frac{12\log\phi}{\pi^2}\;t\;.
$$

Pour conclure cette partie, montrons le corollaire
\ref{coro:introcompirrquad} de l'introduction.

\medskip\noindent{\bf D\'emonstration du corollaire
  \ref{coro:introcompirrquad}.} 
Si $\alpha_0$ est un entier quadratique irrationnel r\'eel, la forme
quadratique binaire $Q_\alpha(X,Y)= X^2-\Tr(\alpha)XY + N(\alpha)Y^2$
est, comme vu dans la preuve du lemme \ref{lem:relathdiscr} (3),
enti\`ere, primitive, irr\'eductible, de discriminant $\discr_Q$ \'egal
au discriminant $D=\Tr(\alpha_0)^2-4N(\alpha_0)$ de
$\OOO=\ZZ+\alpha_0\ZZ$. Par la preuve de la proposition
\ref{prop:bijqaudbinirrquad}, la valeur absolue de la trace de
l'\'el\'ement hyperbolique primitif $\ga_0=\pm\ga_{Q_{\alpha_0}}$ de
$\operatorname{PSL}_2(\ZZ)$ fixant $\alpha_0$ est $t_0$, o\`u, comme
d\'efini dans l'introduction, $(t_0,u_0)$ est la solution fondamentale
de l'\'equation de Pell-Fermat $t^2-Du^2=4$. Par la derni\`ere
assertion de la proposition \ref{prop:equivreciproque}, l'indice $n_0$
de $\ga_0^\ZZ$ dans le stabilisateur de
$\{\alpha_0,\alpha_0^\sigma\}$, qui vaut $2$ si $\alpha_0$ est
r\'eciproque, et $1$ sinon, est \'egal au $n_0$ d\'efini en
introduction. Le r\'esultat d\'ecoule alors du th\'eor\`eme
\ref{theo:appliarithneg2}. \cqfd

\medskip D'autres r\'esultats de comptage peuvent \^etre obtenus en
variant le sous-groupe arithm\'etique. Par exemple, en renvoyant \`a
\cite{Katok92,MacRei03} pour toute information sur les constructions
par les quaternions de groupes fuchsiens, soient
\begin{itemize}
\item[$\bullet$] $A$ une alg\`ebre de quaternion non ramifi\'ee sur
  $\QQ$, de norme $N$ et de trace $\operatorname{Tr}$, comme
  $\Big(\begin{array}{c}p,-1\\\hline\QQ\end{array}\Big)$ o\`u
  $p\in\NN-\{0\}$ v\'erifie $p\equiv 1\mod 4$, de base standard
  $1,i,j,k$ (voir par exemple \cite[page 88]{MacRei03}),
\item[$\bullet$] $\OOO$ un ordre de $A$, comme $\OOO=\ZZ+i\ZZ+j\ZZ+k\ZZ$
  dans l'exemple ci-dessus,
\item[$\bullet$] $\rho:A\otimes\RR \ra \M_2(\RR)$ un isomorphisme
  d'alg\`ebres de quaternions sur $\RR$, comme $a+ib+cj+dk\mapsto
  \Big(\begin{array}{cc}a+b\sqrt{p}& c+d\sqrt{p}\\ -c+d\sqrt{p}&
    a-b\sqrt{p}\end{array}\Big)$ dans l'exemple ci-dessus (voir
  par exemple \cite[page 114]{Katok92}).
\end{itemize}
Posons $\Ga(A,\OOO)=\rho(\{x\in\OOO\;:\;N(x)=1\})/\{\pm {\rm id}\}$, qui
est un r\'eseau non uniforme de ${\rm PSL}_2(\RR)$. Comme $\rho$
pr\'eserve norme et trace, les \'el\'ements hyperboliques de $\Ga(A,\OOO)$
sont les $\pm\rho(x)$ pour les $x\in\OOO$ tels que $|{\rm Tr}\;x|>2$ et
$N(x)=1$, ceux paraboliques (ou l'identit\'e) les $\pm\rho(x)$ pour les
$x\in\OOO$ tels que $|{\rm Tr}\;x|=2$ et $N(x)=1$. Le corollaire
\ref{coro:comptageorbitepointfixehyperbolique} s'applique encore pour
donner des r\'esultats de comptage asymptotique de points fixes
hyperboliques dans une orbite donn\'ee de $\Ga(A,\OOO)$. 

Dans notre exemple, si $(x,y)$ est une solution de l'\'equation de
Pell-Fermat $x^2-py^2=1$ et si $(x',y')$ est une solution de
l'\'equation de Pell-Fermat ${x'}^2-p{y'}^2=-1$ (qui existe car $p\equiv
1\mod 4$), alors $0$ est un point fixe de l'\'el\'ement hyperbolique
$\pm\Big(\begin{array}{cc} x+y\sqrt{p}& 0\\ 0&
  x-y\sqrt{p}\end{array}\Big)$ et $\frac{x'}{1-y'\sqrt{p}}$ est le
point fixe de l'\'el\'ement parabolique
$\pm\Big(\begin{array}{cc}1+x'\sqrt{p}& py'+\sqrt{p}\\ -py'+\sqrt{p}&
  1-x'\sqrt{p}\end{array}\Big)$, que nous pouvons envoyer en $\infty$
par une conjugaison d\'efinie sur $\QQ(\sqrt{p})$. Ainsi, l'ensemble des
points fixes hyperboliques de $\Ga(A,\OOO)$ est un joli sous-ensemble
arithm\'etiquement d\'efini de $\RR$ (form\'e d'\'el\'ements quadratiques
sur $\QQ(\sqrt{p})$), auxquel nous pouvons appliquer le r\'esultat de
comptage \ref{coro:comptageorbitepointfixehyperbolique}.

\subsection{Comptage de repr\'esentations de formes quadratiques 
binaires enti\`eres}
\label{subsec:comptrepquadbin}

Soit $\P$ le sous-ensemble de $\ZZ^2$ des \'el\'ements de
coordonn\'ees premi\`eres entre elles. Il est invariant par l'action
lin\'eaire de $\SL 2\ZZ$. Pour toute forme quadratique binaire
enti\`ere ind\'efinie irr\'eductible $Q$, notons
$$ 
\Psi_Q(s)=\card\big(\autom(Q,\ZZ)\bs\{x\in\P\;:\; |Q(x)|\le s\}\big)
$$
le nombre des repr\'esentations propres non équivalentes par $Q$ des
entiers de valeur absolue au plus $s$. Le dessin ci-dessous décrit,
lorsque $Q(X,Y)=X^2-XY-Y^2$ et $s=11$ les orbites de $\autom(Q,\ZZ)$ sur
les hyperboles d'équation $|Q(x)|= s$.

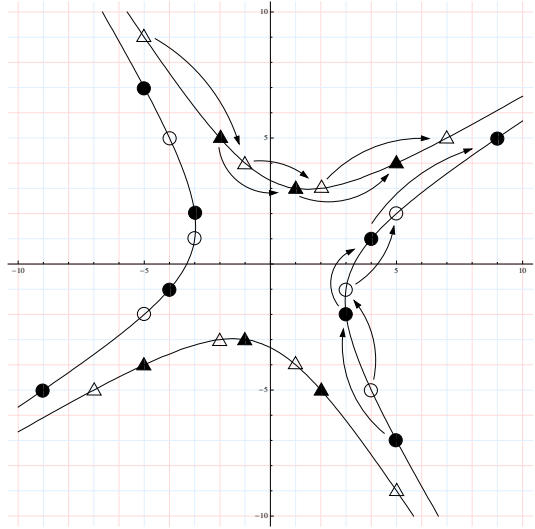
\begin{center}
\begin{picture}(0,0)%
\includegraphics{fig_goldenratio.pstex}%
\end{picture}%
\setlength{\unitlength}{1579sp}%
\begingroup\makeatletter\ifx\SetFigFont\undefined%
\gdef\SetFigFont#1#2#3#4#5{%
  \reset@font\fontsize{#1}{#2pt}%
  \fontfamily{#3}\fontseries{#4}\fontshape{#5}%
  \selectfont}%
\fi\endgroup%
\begin{picture}(8267,8267)(1201,-9228)
\end{picture}%

\end{center}

Notons $\discr'$ le quotient du discriminant $\discr_Q$ de $Q$ par le 
carr\'e du plus grand diviseur commun des coefficients de $Q$. 
Notons $(t_Q,u_Q)$ la solution fondamentale de
l'\'equation de Pell-Fermat 
$t^2-\discr'\, u^2=4$, et
$$
\epsilon_{Q}=\frac{t_Q+u_Q\sqrt{\discr'}}{2}\;\;\;{\rm
  et}\;\;\;R_Q=\log \epsilon_Q\;.
$$

\btheo\label{theo:formes} Soit $Q$ une forme quadratique binaire
enti\`ere ind\'efinie de discriminant non carr\'e. Alors quand $s$
tend vers $+\infty$,
$$
\Psi_Q(s)\sim
\frac{12 \;R_Q}{\pi^2\sqrt{\discr_Q}}\; s\;.
$$  
\etheo

Le corollaire \ref{coro:introformes} de l'introduction en découle.

\medskip 
\dem 
Puisque $\Psi_{kQ}(s)=\Psi_Q(\frac{s}{k})$, $D_{kQ}=k^2D_Q$, et
$\epsilon_{Q}= \epsilon_{kQ}$, nous pouvons supposer que $Q$ est
primitive.

Nous appliquons le corollaire \ref{coro:comptagedoubleclasse} avec
$n=2$, $k=1$, $\Ga=\operatorname{PSL}_2(\ZZ)$, une seule horoboule
$\H_\infty$ (de stabilisateur not\'e $\Ga_\infty$) constitu\'ee des
points du demi-plan supérieur de coordonn\'ee verticale au moins $1$,
et une seule g\'eod\'esique $C_0$ \'egale \`a l'axe de translation de
$\pm\ga_Q$. Il est imm\'ediat que
$$
\Vol(\Ga_\infty\backslash \H_\infty)= 1\;\;\;{\rm et}
\;\;\;\Vol(\Ga\backslash \HH^2_\RR)=\frac{\pi}{3}\;.
$$
Notons $i_Q$ l'indice de $\PSO$ dans $\Ga_0=\Stab_\Ga C_0$ (qui vaut
$2$ ou $1$ suivant que $Q$ est r\'eciproque ou pas). Puisque
$\pm\ga_Q$ engendre $\PSO$, et comme la trace de $\ga_Q$ est $t_Q$,
nous avons, par la formule \eqref{eq:calclongtransl} et un petit
calcul,
$$
\Vol(\Ga_0\backslash C_0)=\frac{2\;\arcosh (\frac{t_Q}{2})}{i_Q}=
\frac{2\log \epsilon_Q}{i_Q}=
\frac{2R_Q}{i_Q}\;.
$$

Le groupe $\SL 2\ZZ$ agit \`a droite sur $\ZZ^2$ par $(\ga,x)\mapsto
\ga^{-1}x$, l'orbite de $(1,0)$ est exactement l'ensemble des couples
$(d,-c)$ d'entiers premiers entre eux, et le stabilisateur de $(1,0)$
est le sous-groupe unipotent triangulaire sup\'erieur $U$. Par le
lemme \ref{lem:relathdiscr} (2), $\Psi_Q(s)$ est donc le nombre de
doubles classes dans $U\backslash \SL 2\ZZ/\autom(Q,\ZZ)$ des
\'el\'ements $\ga$ de $\SL 2\ZZ$ tels que $h(\ga\alpha_Q)\leq
\frac{2}{\sqrt{D_Q}}\;s$.  Puisque $-\operatorname{id}$ appartient \`a
$\autom(Q,\ZZ)$, l'application canonique de $U\backslash \SL
2\ZZ/\autom(Q,\ZZ)$ dans $\Ga_\infty\backslash \Ga/\Ga_0$ est
surjective, chaque pr\'eimage est d'ordre $i_Q$, et la valeur de
$h(\ga\alpha_Q)$ est constante sur chaque pr\'eimage. Comme d\'ej\`a
vu dans la preuve du corollaire
\ref{coro:comptageorbitepointfixehyperbolique}, nous avons
$d(\H_\infty, \ga C_0) =\log h(\ga\alpha_Q)$. Donc, en utilisant le
corollaire \ref{coro:comptagedoubleclasse},
\begin{align*}
\Psi_Q(s)&= i_Q\operatorname{Card}\big\{[\ga]\in\Ga_{\infty}\backslash
\Ga/\Ga_{0}\;:\; \;d(\H_\infty, \ga C_0)\leq \log
\Big(\frac{2}{\sqrt{D_Q}}\;s\Big)\big\}\\
&\sim i_Q\;\frac{\Vol(\SSS_{0})
\Vol(\Ga_{\infty}\backslash\H_\infty)\Vol(\Ga_{0}\backslash C_0)}
{\Vol(\SSS_{1})\Vol(\Ga\backslash\HH^2_\RR)} 
\;\Big(\frac{2}{\sqrt{D_Q}}\;s\Big)
=\frac{12\, R_{Q}}{\pi^2\sqrt{\discr_Q}}\; s
\;.\;\;\mbox{\cqfd}
\end{align*}

\medskip\rem Soient $p\in\NN-\{0,1\}$ et $Q(X,Y)=AX^2+BXY+CY^2$ une
forme quadratique binaire enti\`ere ind\'efinie de discriminant non
carr\'e. Consid\'erons les sous-groupes de congruence $\Ga(p)$ et
$\Ga_0(p)$ d\'efinis dans la partie \ref{subsec:comptirraquad}. En
reprenant les notations de la preuve de la proposition \ref
{prop:bijqaudbinirrquad}, la preuve de \cite[Prop.~3.3]{Sarnak82}
montre, lorsque $Q$ est primitive, que $\ga_{Q,t,u}$ appartient \`a
$\Ga(p)$ si et seulement si $p$ divise $u$. De m\^eme, si $A\equiv 1\mod
p$, alors $\ga_{Q,t,u}$ appartient \`a $\Ga_0(p)$ si et seulement si $p$
divise $u$. Notons encore $\discr'$ le quotient du discriminant
$\discr_Q$ de $Q$ par le carr\'e du plus grand diviseur commun des
coefficients de $Q$.  Notons $(t_{Q,p},u_{Q,p})$ la solution de
l'\'equation de Pell-Fermat $t^2-\discr'\, u^2=4$ telle que
$t_{Q,p}>0,u_{Q,p}>0$, $p$ divise $u_{Q,p}$ et $u_{Q,p}$ est
minimale. Posons
$$
\epsilon_{Q,p}=\frac{t_{Q,p}+u_{Q,p}\sqrt{\discr'}}{2}\;\;\;{\rm
  et}\;\;\;R_{Q,p}=\log\epsilon_{Q,p}\;.
$$
Posons $\P(p)$ l'ensemble des $(X,Y)$ dans $\P$ tels que $X\equiv 1\mod
p$ et $Y\equiv 0\mod p$, qui est invariant par $\Ga(p)$, et $\P_{0}(p)$
l'ensemble des $(X,Y)$ dans $\P$ tels que $Y\equiv 0\mod p$, qui est
invariant par $\Ga_0(p)$. Notons
$$
\Psi_{Q,p}(s)=\card\Big(\big(\autom(Q,\ZZ)\cap\Ga(p)\big)\bs
\{(X,Y)\in\P(p)\;:\;
|Q(X,Y)|\le s\} \Big)
$$
et
$$
\Psi_{Q,p,0}(s)=\card\Big(\big(\autom(Q,\ZZ)\cap\Ga_0(p)\big)\bs
\{(X,Y)\in\P_{0}(p)\;:\; |Q(X,Y)|\le s\}\Big)\;.
$$
En appliquant le m\^eme raisonnement que dans la preuve du
th\'eor\`eme \ref{theo:formes} en rempla\c{c}ant le groupe
$\Ga=\operatorname{PSL}_2(\ZZ)$ par $\Ga(p)$ et $\Ga_0(p)$ et en
utilisant les formules \eqref{eq:indexprincipal} et
\eqref{eq:indexnonprincipal}, nous avons, quand $s$ tend vers
$+\infty$,
$$
\Psi_{Q,p}(s)\sim
\begin{cases}
\dfrac
{4\;R_{Q,p}}{\pi^2\sqrt{\discr_Q}}\;s
\quad\textrm{si}\;p=2,\\
\phantom{hhhhhhhh}\\
\dfrac{24\;R_{Q,p}}{\pi^2p^2\sqrt{\discr_Q}}  
\;\;{\displaystyle{\prod_{\ppp \;{\rm premier},\; \ppp | p} 
\big(1-\frac{1}{\ppp^2}\big)}^{-1}}\;\;s\;\;\;{\rm si} \;\;p>2\;,  
\end{cases}
$$
et, si $A\equiv 1\mod p$, 
$$
\Psi_{Q,p,0}(s)
\sim\dfrac{12\;R_{Q,p}}{\pi^2p\sqrt{\discr_Q}}  
\;\;{\displaystyle{\prod_{\ppp \;{\rm premier},\; \ppp | p} 
\big(1+\frac{1}{\ppp}\big)}^{-1}}\;\;s\;.
$$

\subsection{Comptage d'irrationnels quadratiques 
imaginaires dans une orbite de sous-groupes 
de congruences du groupe de Bianchi}
\label{subsec:comptirraquadimag}

\'Etant donn\'e un corps de nombres quadratique imaginaire $K$, le but
de cette partie est de donner un \'equivalent quand $s$ tend vers
$+\infty$ du nombre modulo translations enti\`eres sur $K$
d'irrationnels quadratiques sur $K$, appartenant \`a une orbite
donn\'ee par homographies et anti-homographies enti\`eres sur $K$, de
complexit\'e au plus $s$, o\`u la complexit\'e d'un irrationnel
quadratique sur $K$ est mesur\'ee par l'inverse de sa distance \`a son
conjugu\'e de Galois.

Si $K$ est un corps de nombres, rappelons que la fonction z\'eta (de
Dedekind) de $K$, d\'efinie sur $\{z\in\CC\;:\;\operatorname{Re} s>1\}$,
est
$$
\zeta_K(s)=\sum_{\aaa}\frac 1{N(\aaa)^s}\;,
$$
o\`u $\aaa$ parcourt les id\'eaux non nuls de l'anneau $\OOO$ des
entiers de $K$ et $N(\aaa)=\card(\OOO/\aaa)$.

\btheo\label{theo:appliarithdim3} Soit $D$ un entier strictement
n\'egatif, congru \`a $0$ ou \`a $1$ modulo $4$, sans facteur
carr\'e. Soit $\OOO_D$ l'anneau des entiers du corps de nombres
quadratique imaginaire $K_D$ de discriminant $D$. Soit $\alpha_0$ un
nombre complexe, irrationnel quadratique sur $K_D$. Soient $G$ un
sous-groupe d'indice fini de $\PSL 2{\OOO_D}$, $\cdot$ son action par
homographies sur $\CC\cup \{\infty\}$, et $\Lambda$ le r\'eseau des
$\lambda\in\OOO_D$ tels que $\pm\Big(\begin{array}{cc}1&\lambda
  \\0&1\end{array} \Big)\in G$.  Alors quand $s$ tend vers $+\infty$,
$$
\operatorname{Card}\{\alpha\in G \cdot
\{\alpha_0,\alpha_0^\sigma\}\!\!\!\mod \Lambda\;:\; h(\alpha)\leq s\}
\sim \;\frac{4 \pi^2\;n_{G,\infty}
  \;\big|\log\big|\frac{\operatorname{tr}\ga_0+
    \sqrt{\operatorname{tr}\ga_0^2-4}}{2}\big|\;\big|}
{n_G\;n_{G,0}\;\omega_D\;|D|\;\zeta_{K_D}(2)}\; \;s^2\;,
$$
o\`u $\omega_D$ est le nombre de racines de l'unit\'e dans $\OOO_D$;
$n_G$ est l'indice de $G$ dans $\PSL 2{\OOO_D}$; $n_{G,\infty}$ est
l'indice de $\Lambda$ dans le stabilisateur de $[1:0]$ dans $\PSL
2{\OOO_D}$; $\ga_0$ est un \'el\'ement de $G$ fixant $\alpha_0$ et de
modules de valeurs propres diff\'erents de $1$; et $n_{G,0}$ est
l'indice de ${\ga_0}^\ZZ$ dans le stabilisateur de
$\{\alpha_0,\alpha_0^\sigma\}$ dans $G$.  \etheo

\dem Le groupe $G$, d'indice fini dans le sous-groupe arithm\'etique
$\PSL 2{\OOO_D}$ de $\PSL 2\CC$, est un groupe discret de covolume fini
d'isom\'etries de $\htr$. Par une formule due \`a Humbert (voir par
exemple les parties 8.8 and 9.6 de \cite{ElsGruMen98}), nous avons
$$
\Vol(G\bs\htr)=n_G\Vol(\PSL 2{\OOO_D}\bs\htr)=
\frac{n_G}{4\pi^2}|D|^{3/2}\zeta_{K_D}(2)\;.
$$

Le point $\infty$ de $\partial_\infty\htr$ est le point fixe d'un
\'el\'ement parabolique (par exemple $z\mapsto z+1$) de $\PSL
2{\OOO_D}$, donc d'un \'el\'ement parabolique de son sous-groupe
d'indice fini $G$. Notons $G_\infty$ le stabilisateur de $\infty$ dans
$G$. Il contient $\Lambda$ avec indice fini $n'_{G,\infty}$, et il est
contenu dans le stabilisateur de $\infty$ dans $\PSL 2{\OOO_D}$ avec
indice fini $n''_{G,\infty}$. Nous avons donc $n_{G,\infty}=
n'_{G,\infty}n''_{G,\infty}$. Un domaine fondamental pour l'action sur
$\CC$ du stabilisateur de $\infty$ dans $\PSL 2{\OOO_D}$ est le
rectangle $[0,1]\times [0,\frac{\sqrt{|D|}}{2}]$ si $D\neq-3,-4$, le
rectangle $[0,1]\times[0,\frac{1}{2}]$ si $D=-4$, et le
parall\'elogramme de sommets $0,\frac{1}{2}-\frac{i}{2\sqrt{3}},
\frac{1}{2}+\frac{i}{2\sqrt{3}}, \frac{i}{\sqrt{3}}$ si $D= -3$, voir
par exemple \cite[page 318]{ElsGruMen98}. Notons que $\omega_D$ vaut
respectivement $2,4,6$ dans ces trois cas. Donc l'aire d'un domaine
fondamental mesurable pour l'action de $G_\infty$ sur l'horosph\`ere
de $\htr$ form\'e des points de hauteur euclidienne $1$ est
$$
\A_\infty=\frac{n''_{G,\infty}\sqrt{|D|}}{\omega_D}\;.
$$

Remarquons que $\ga_0$ existe bien, et est hyperbolique, voir par
exemple \cite[Lem.~6.2]{ParPau10MZ}, et $n_{G,0}$ est donc bien
d\'efini. Par un calcul classique, la distance de translation $\ell_0$
d'un \'el\'ement $\ga_0$ de $\PSL 2\CC$ v\'erifie
$$
\ell_0=
2\,\Big|\log\big|\frac{\operatorname{tr}\ga_0+
    \sqrt{\operatorname{tr}\ga_0^2-4}}{2}\big|\;\Big|\;,
$$
qui ne d\'epend pas du choix de la racine carr\'ee du nombre complexe
$\operatorname{tr}\ga_0^2-4$, ni du choix du repr\'esentant de $\ga_0$
dans $\SL 2\CC$ modulo $\{\pm\id\}$.

\medskip
Nous pouvons donc appliquer le corollaire
\ref{coro:comptageorbitepointfixehyperbolique}, avec $n=3$, $G=G$
et $t=\frac{s}{2}$, en remarquant que nous avons $\card(G_\infty \cap
\operatorname{Stab}_G\alpha)>1$ seulement pour un nombre fini
d'orbites modulo $G_\infty$ d'\'el\'ements $\alpha$ de $G\cdot
\{\alpha_0,\alpha_0^\sigma\}$, et que pour passer d'un comptage
modulo $G_\infty$ \`a un comptage modulo $\Lambda$, il faut
multiplier le premier nombre par $n'_{G,\infty}$.  
\cqfd

\medskip Voici un exemple d'application de ce th\'eor\`eme
\ref{theo:appliarithdim3}, dont le corollaire
\ref{coro:introorbianchi} de l'introduction d\'ecoule (voir la
remarque suivant la d\'emonstration).

\bcoro \label{coro:orbianchi}
Soit $\phi=\frac{1+\sqrt{5}}2$ le nombre d'or, soit $\aaa$ un
id\'eal non nul de $\OOO_D$, et soit $\Ga_{0}(\aaa)$ le sous-groupe de
congruence $\Big\{\pm\left(\begin{array}{cc}a & b \\ c & d \end{array}
\right)\in {\rm PSL}_2(\OOO_D) \;:\; c\in\aaa\Big\}$.  Alors quand $s$
tend vers $+\infty$, le cardinal de $\{\alpha\in \Ga_{0}(\aaa) \cdot
\{\phi,\phi^{\sigma}\} \!\!\!\mod \OOO_D\;:\; h(\alpha)\leq s\} $ est
\'equivalent \`a
$$
\frac{4\pi^2\;k_\aaa\;\log\phi} {|D|\;  \zeta_{K_D}(2)\;
  N(\aaa)\prod_{\ppp |\aaa} \big(1+\frac{1}{N(\ppp)}\big)}\; \;s^2\;,
$$
o\`u $k_\aaa$ est le plus petit entier $k\in\NN-\{0\}$ tel que le
$2k$-\`eme terme de la suite de Fibonacci appartienne \`a $\aaa$, si
$D\neq -4$ et si $\phi^\sigma$ n'est pas dans la m\^eme orbite que
$\phi$ par $\Ga_{0}(\aaa)$.  
\ecoro

\dem Appliquons le th\'eor\`eme \ref{theo:appliarithdim3} ci-dessus avec
$\alpha_0=\phi$ et $G=\Ga_{0}(\aaa)$. Alors $\Lambda=\OOO_D$,
$n_G=N(\aaa) \prod_{\ppp |\aaa} \big(1+\frac{1}{N(\ppp)}\big)$ par
\cite[page 117]{Newman72} (o\`u comme d'habitude, le produit est pris
sur les id\'eaux premiers divisant $\aaa$), et $n_{G,\infty}=
\frac{\omega_D}{2}$.

Soit $\ga_1=\Big(\begin{array}{cc}2&1 \\1&1\end{array}\Big)$ si $D\neq
-4$ et $\ga_1=\Big(\begin{array}{cc}i&i \\i&0\end{array}\Big)$
sinon. Alors $\pm \ga_1$ appartient \`a $\PSL 2{\OOO_D}$ et fixe $\phi$
(donc $\phi^\sigma$). Notons $k'_\aaa$ le plus petit entier
strictement positif $k$ tel que $\pm\ga_1^k$ appartienne \`a $\Ga_{0}
(\aaa)$, qui existe car $\Ga_{0} (\aaa)$ est d'indice fini dans $\PSL
2{\OOO_D}$. Posons $\ga_0=\pm\ga_1^{k'_\aaa}$, qui appartient \`a $\Ga_{0}
(\aaa)$. Remarquons que la distance de translation de $\ga_0$ vaut
$k'_\aaa$ fois celle de $\pm\ga_1$, donc
$$
\Big|\log\Big|\frac{\operatorname{tr}\ga_0+
  \sqrt{\operatorname{tr}\ga_0^2-4}}{2}\Big|\Big| =
\left\{\begin{array}{ll}
    2 k'_\aaa\log\phi& {\rm si}\;D\neq -4\\
    k'_\aaa\log\phi& {\rm sinon}\;. \end{array} \right.
$$

Notons $c(\ga)$ le coefficient $2$-$1$ d'un \'el\'ement $\ga\in\SL
2\CC$. Si $P= \Big(\begin{array}{cc} \phi&\phi^\sigma
  \\1&1\end{array}\Big)$, alors $\ga_1= P \Big(\begin{array}{cc}
  \phi^2&0 \\0&\phi^{-2}\end{array} \Big) P^{-1}$ si $D\neq -4$ et
$\ga_1= P \Big(\begin{array}{cc} i\phi&0 \\0&-i\phi^{-1}\end{array}
\Big) P^{-1}$ sinon.  Si $D\neq -4$, pour tout $k\in\NN-\{0\}$, nous
avons donc $c({\ga_1}^k)= \frac{\phi^{2k}-\phi^{-2k}}{\sqrt{5}}$. Ce
nombre est le $2k$-\`eme terme $x_{2k}$ de la suite de Fibonacci
$(x_n)_{n\in\NN}$, d\'efinie par la relation de r\'ecurrence
$x_{n+2}=x_{n+1}+x_n$, de valeurs initiales $x_0=0$ et
$x_1=1$. L'entier $k_\aaa$ d\'efini ci-dessus est donc \'egal \`a
$k'_\aaa$ si $D\neq -4$.

\medskip Calculons maintenant $n_{G,0}$. Notons $F_G$ (respectivement
$F_\Ga$) le fixateur de $\{\phi,\phi^\sigma\}$ dans $\Ga_{0} (\aaa)$
(respectivement $\PSL 2{\OOO_D}$). Pour tout $\ga\in\SL 2\CC$ tel que
$\pm\ga\in F_\Ga$, il existe $\rho=\rho(\ga)>0$ et
$\theta\in\RR/2\pi\ZZ$ tel que $\ga= P \Big(\begin{array}{cc} \rho
  e^{i\theta}&0 \\0& \rho^{-1} e^{-i\theta}\end{array} \Big) P^{-1}$;
donc $|c(\ga)|^2=\frac{\rho^2+\rho^{-2} -2\cos(2\theta)}{5}$. En
particulier, si $\rho<\phi$, alors $|c(\ga)|^2<1$. Quitte \`a
remplacer $\ga$ par son inverse, nous pouvons supposer que $\ga$ ne
translate pas de $\phi^\sigma$ vers $\phi$ sur la g\'eod\'esique entre
$\phi$ et $\phi^\sigma$. Rappelons que la norme d'un \'el\'ement de
$\OOO_D$, \'egale au carr\'e du module, est enti\`ere, et que l'ensemble
des distances de translations des \'el\'ements de $F_\Ga$ est un
sous-groupe discret de $\RR$. Par cons\'equent, si $\rho<\phi$, alors
$c(\ga)=0$ et donc $\pm\ga$, qui fixe $\infty,\phi,\phi^\sigma$, vaut
l'identit\'e. Nous ne pouvons avoir $\phi<\rho<\phi^2$, sinon
$0<\rho(\Big(\begin{array}{cc}2&1
  \\1&1\end{array}\Big)\ga^{-1})<\phi$, ce que nous venons
d'exclure. Si $\rho=\phi$, alors $|c(\ga)|^2\leq 1$, avec \'egalit\'e
si et seulement si $\theta=\frac{\pi}{2} \mod \pi$. L'\'el\'ement
$\ga$ est alors \'egal \`a $\pm\Big(\begin{array}{cc}i&i
  \\i&0\end{array}\Big)$, qui appartient \`a $\PSL 2{\OOO_D}$ si et
seulement si $D=-4$. Si $\rho=\phi^2$, alors $\theta=0\mod \pi$, sinon
$\rho(\Big(\begin{array}{cc}2&1
  \\1&1\end{array}\Big)\ga^{-1})=1<\phi$, ce que nous avons
exclu. Donc $F_\Ga$ est \'egal \`a $(\pm \ga_1)^\ZZ$. Par d\'efinition
de $k'_\aaa$, nous en d\'eduisons que $F_G$ est \'egal \`a
${\ga_0}^\ZZ$.  L'indice $n_{G,0}$ de ${\ga_0}^\ZZ$ dans le
stabilisateur de $\{\phi,\phi^\sigma\}$ dans $\Ga_{0}(\aaa)$ vaut donc
$2$ s'il existe un \'el\'ement de $\Ga_{0}(\aaa)$ envoyant $\phi$ sur
$\phi^\sigma$, et $1$ sinon.

Le r\'esultat d\'ecoule alors du th\'eor\`eme \ref{theo:appliarithdim3}. 
\cqfd

\medskip En notant $\delta_{\xi,\eta}$ le symbole de Kronecker (qui
vaut $1$ si $\xi=\eta$ et $0$ sinon), la d\'emonstration ci-dessus
montre que le corollaire \ref{coro:orbianchi} reste valable dans tous
les cas, en divisant l'asymptotique par $(2-\delta_{D,-4})
(1+\delta_{\Ga_{0}(\aaa)\cdot\phi,\Ga_{0}(\aaa)\cdot\phi^\sigma})$ et
en y rempla\c{c}ant $k_\aaa$ par l'entier $k'_\aaa$ d\'efini dans la
d\'emonstration ci-dessus. Puisque $c(P \Big(\begin{array}{cc} i\phi
  &0 \\0& -i\phi^{-1}\end{array} \Big) P^{-1})= i^n\,
\frac{\phi^n-(\phi^\sigma)^n}{\sqrt{5}}$, cet entier $k'_\aaa$ est le
plus petit entier $k$ tel que le $k$-\`eme terme de la suite de
Fibonacci appartienne \`a $\aaa$.

\medskip \rem Soient $K$ un corps de nombres quadratique imaginaire,
$\aaa$ un id\'eal non nul de l'anneau $\OOO_K$ de ses entiers, et $G$ un
sous-groupe d'indice fini de $\PSL 2{\OOO_K}$. Un irrationnel
quadratique $\alpha$ sur $K$ est dit {\it $G$-r\'eciproque} si
$\alpha$ et son conjugu\'e de Galois $\alpha^\sigma$ sur $K$ sont dans
la m\^eme orbite par homographies de $G$. Par exemple, les seuls
\'el\'ements de $\PSL 2{\OOO_K}$ envoyant $\phi$ sur $\phi^\sigma$ sont
les $\alpha_k=\Big(\begin{array}{cc}0&-1 \\1&0\end{array}\Big)\ga_1^k$
pour $k\in\NN$, o\`u $\ga_1$ est d\'efini dans la preuve ci-dessus. Si
$D\neq -4$, le coefficient $2$-$1$ de $\alpha_k$ est, au signe pr\`es,
le $(2k+1)$-\`eme terme de la suite de Fibonacci. Il est facile de
v\'erifier que celui-ci vaut $2$ si $k=1$, et n'est pas congru \`a $0$
modulo $3$. Si $D\neq -4$, nous en d\'eduisons donc que $\phi$ est
$\Ga_{0}(2)$-r\'eciproque, mais n'est pas $\Ga_{0}(3)$-r\'eciproque.

\bibliography{../biblio}

\bigskip
\noindent {\small 
\begin{tabular}{l}  
Department of Mathematics and Statistics, P.O. Box 35\\ 
40014 University of Jyv\"askyl\"a, FINLAND\\ 
{\it e-mail: parkkone@maths.jyu.fi} 
\end{tabular} 
\medskip

\noindent \begin{tabular}{l}
DMA, UMR 8553 CNRS\\
Ecole Normale Sup\'erieure, 45 rue d'Ulm\\
75230 PARIS Cedex 05, FRANCE\\
{\it e-mail: Frederic.Paulin@ens.fr}
\end{tabular} and \begin{tabular}{l}
D\'epartement de math\'ematique, B\^at.~425\\
Universit\'e Paris-Sud 11\\
91405 ORSAY Cedex, FRANCE\\
{\it e-mail: frederic.paulin@math.u-psud.fr}
\end{tabular} 
}

\end{document}